# The Congruences of a Finite Lattice

## A *Proof-by-Picture* Approach

## Third Edition

**George Grätzer**

*To László Fuchs,*
*my thesis advisor, my teacher,*
*who taught me to set the bar high;*

*and to my coauthors,*
*who helped to raise the bar;*
*especially to Gábor Czédli,*
*for his help in the last dozen years or so.*

# *Short Contents*

# *Contents*

# *Glossary of Notation*

| Symbol | Explanation | Page |
|---|---|---|
| $a_*$ $(a^*)$ | the unique lower (upper) cover of $a$ | 17 |
| $0_I$ and $1_I$ | the zero and unit of the interval $I$ | 12 |
| $\mathrm{Atom}(U)$ | set of atoms of the ideal $U$ | 106 |
| $\mathrm{Aut}\, L$ | automorphism group of $L$ | 11 |
| $\mathsf{B}_n$ | Boolean lattice with $n$ atoms | 4 |
| $\mathrm{C_l}(L), \mathrm{C_{ll}}(L), \mathrm{C_{ul}}(L)$ | left boundary chains of a planar lattice | 47 |
| $\mathsf{C}_n$ | $n$-element chain | 4 |
| $\mathrm{con}(a,b)$ | smallest congruence under which $a \equiv b$ | 15 |
| $\mathrm{con}(H)$ | smallest congruence collapsing $H$ | 15 |
| $\mathrm{con}(\mathfrak{p})$ | principal congruence generated by $\mathfrak{p}$ | 38 |
| $\mathrm{Con}\, L$ | congruence lattice of $L$ | 14, 70 |
| $\mathrm{J}(\mathrm{Con}\, L)$ | join-irreducible congruences of $L$ | 38 |
| $\mathrm{C_r}(L), \mathrm{C_{lr}}(L), \mathrm{C_{ur}}(L)$ | right boundary chains of a planar lattice $L$ | 47 |
| $\mathrm{Cube}\, K$ | cubic extension of $K$ | 93 |
| $\mathbf{D}$ | class (variety) of distributive lattices | 22 |
| $\mathrm{Diag}(K)$ | diagonal embedding of $K$ into $\mathrm{Cube}\, K$ | 93 |
| $\mathrm{Dn}\, P$ | down sets of the ordered set $P$ | 5, 9 |
| $\mathrm{ext}\colon \mathrm{Con}\, K \to \mathrm{Con}\, L$ | for $K \leq L$, extension map: $\boldsymbol{\alpha} \mapsto \mathrm{con}_L(\boldsymbol{\alpha})$ | 41 |
| $\mathrm{fil}(a)$ | filter generated by the element $a$ | 13 |
| $\mathrm{fil}(H)$ | filter generated by the set $H$ | 13 |
| $\mathrm{Flag}(c_3)$ | flag lattice | 345 |
| $\mathrm{Free}_{\mathbf{D}}(3)$ | free distributive lattice on three generators | 23 |
| $\mathrm{Free}_{\mathbf{K}}(H)$ | free lattice generated by $H$ in a variety $\mathbf{K}$ | 23 |
| $\mathrm{Free}_{\mathbf{M}}(3)$ | free modular lattice on three generators | 25 |
| $\mathrm{Fuse}(P, A)$ | Fusion of ordered sets | 315 |
| $\mathrm{hom}_{\{\vee,0\}}(X, Y)$ | $\{\vee, 0\}$-homomorphisms of $X$ into $Y$ | 285 |
| $\mathrm{id}(a)$ | ideal generated by the element $a$ | 13 |
| $\mathrm{id}(H)$ | ideal generated by the set $H$ | 13 |
| $\mathrm{Id}\, L$ | ideal lattice of $L$ | 13, 71 |
| (Id) | condition to define ideals | 13, 71 |

| Symbol | Explanation | Page |
|---|---|---|
| $\mathrm{J}(D)$ | ordered set of join-irreducible elements of $D$ | 17 |
| $\mathrm{J}^{+}(D)$ | $\mathsf{J}(D) \cup \{0, 1\}$ | 307 |
| $\mathsf{J}(a)$ | set of join-irreducible elements below $a$ | 17 |
| $\ker(\boldsymbol{\gamma})$ | congruence kernel of $\boldsymbol{\gamma}$ | 15 |
| $\mathbf{L}$ | class (variety) of all lattices | 22 |
| $L_{\mathrm{bottom}}, L_{\mathrm{top}}$ | bottom and top of a rectangular lattice $L$ | 55 |
| lc(L) | left corner of a rectangular lattice $L$ | 47 |
| $L_{\mathrm{left}}, L_{\mathrm{right}}$ | left and right boundary of $L$ | 47 |
| $\mathbf{M}$ | class (variety) of modular lattices | 22 |
| Max | maximal elements of an ordered set | 71 |
| $\mathbf{mcr}(n)$ | minimal congruence representation function | 126 |
| $\mathbf{mcr}(n, \mathbf{V})$ | **mcr** for a class $\mathbf{V}$ | 126 |
| $\mathsf{M}_3$ | five-element modular nondistributive lattice | 28 |
| $\mathsf{M}_3[L]$ | ordered set of Boolean triples of $L$ | 80 |
| $\mathsf{M}_3[L, a]$ | interval of $\mathsf{M}_3[L]$ | 85 |
| $\mathsf{M}_3[L, a, b]$ | interval of $\mathsf{M}_3[L]$ | 88 |
| $\mathsf{M}_3[a, b]$ | Boolean triples of $[a, b]$ | 80 |
| $\mathrm{M}(D)$ | ordered set of meet-irreducible elements of $D$ | 17 |
| $\mathsf{N}_5$ | five-element nonmodular lattice | 28 |
| $\mathsf{N}_{5,5}$ | seven-element nonmodular lattice | 118 |
| $\mathsf{N}_6 = N(p, q)$ | six-element nonmodular lattice | 102 |
| $N(A, B)$ | a lattice construction | 166 |
| $\mathrm{Part}\, A$ | partition lattice of $A$ | 7, 9 |
| $\mathrm{Pow}\, X$ | power set lattice of $X$ | 5 |
| $\mathrm{Princ}\, L$ | ordered set of principal congruences $L$ | 37 |
| rc(L) | right corner of a rectangular lattice | 47 |
| $\mathrm{Prime}(L)$ | set of prime intervals of $L$ | 38 |
| $\mathrm{re}\colon \mathrm{Con}\, L \to \mathrm{Con}\, K$ | restriction map: $\boldsymbol{\alpha} \mapsto \boldsymbol{\alpha}\rceil K$ | 40 |
| $(\mathrm{SD}_{\vee})$, $(\mathrm{SD}_{\wedge})$ | semidistributive laws | 50 |
| **SecComp** | class of sectionally complemented lattices | 18 |
| **SemiMod** | class of semimodular lattices | 141 |
| $\mathrm{Simp}\, K$ | simple extension of $K$ | 93 |
| $(\mathrm{SP}_{\vee})$, $(\mathrm{SP}_{\wedge})$ | substitution properties | 13 |
| $\mathrm{Split}(P, a)$ | splitting an element | 318 |
| $S(p, q)$, $W(p, q)$ | small lattices used in lattice constructions | 300 |
| $\mathrm{sub}(H)$ | sublattice generated by $H$ | 12 |
| Swing, $\hookrightarrow$, $\overset{\mathrm{in}}{\hookrightarrow}$, $\overset{\mathrm{ex}}{\hookrightarrow}$ | $\mathfrak{p} \hookrightarrow \mathfrak{q}$, $\mathfrak{p}$ swings to $\mathfrak{q}$ | 367 |
| $\mathrm{Traj}\, L$ | set of all trajectories of $L$ | 48 |

| Symbol | Explanation | Page |
|---|---|---|
| **Relations and Congruences** | | |
| $A^2$ | binary relations | 3 |
| $\boldsymbol{\alpha}, \boldsymbol{\beta}, \ldots$ | congruences | 13 |
| $\mathbf{0}$ | zero of $\operatorname{Part} A$ and $\operatorname{Con} L$ | 8 |
| $\mathbf{1}$ | unit of $\operatorname{Part} A$ and $\operatorname{Con} L$ | 8 |
| $a \equiv b \pmod{\pi}$ | $a$ and $b$ in the same block of $\pi$ | 7 |
| $a \varrho b$ | $a$ and $b$ in relation $\varrho$ | 3 |
| $a \equiv b \pmod{\boldsymbol{\alpha}}$ | $a$ and $b$ in relation $\boldsymbol{\alpha}$ | 3 |
| $a/\pi$ | block containing $a$ | 7, 13 |
| $H/\pi$ | blocks represented by $H$ | 7 |
| $\alpha \circ \beta$ | product of $\alpha$ and $\beta$ | 18 |
| $\alpha \overset{\mathrm{r}}{\circ} \beta$ | reflexive product of $\alpha$ and $\beta$ | 26 |
| $\boldsymbol{\alpha}\rceil_K$ | restriction of $\boldsymbol{\alpha}$ to the sublattice $K$ | 13 |
| $L/\boldsymbol{\alpha}$ | quotient lattice | 15 |
| $\boldsymbol{\beta}/\boldsymbol{\alpha}$ | quotient congruence | 16 |
| $\pi_i$ | projection map: $L_1 \times \cdots \times L_n \to L_i$ | 18 |
| $\boldsymbol{\alpha} \times \boldsymbol{\beta}$ | direct product of congruences | 18 |
| **Ordered sets** | | |
| $\leq, <$ | ordering | 3 |
| $\geq, >$ | ordering, inverse notation | 3 |
| $K \leq L$ | $K$ a sublattice of $L$ | 12 |
| $\leq_Q$ | ordering of $P$ restricted to a subset $Q$ | 4 |
| $a \parallel b$ | $a$ incomparable with $b$ | 3 |
| $a \prec b$ | $a$ is covered by $b$ | 5 |
| $b \succ a$ | $b$ covers $a$ | 5 |
| $0$ | zero, least element of an ordered set | 4 |
| $1$ | unit, largest element of an ordered set | 4 |
| $a \vee b$ | join operation | 8 |
| $\bigvee H$ | least upper bound of $H$ | 3 |
| $a \wedge b$ | meet operation | 8 |
| $\bigwedge H$ | greatest lower bound of $H$ | 4 |
| $P^\delta$ | dual of the ordered set (lattice) $P$ | 4, 10 |
| $[a, b]$ | interval | 12 |
| $\downarrow H$ | down set generated by $H$ | 5 |
| $\downarrow a$ | down set generated by $\{a\}$ | 5 |
| $P \cong Q$ | ordered set (lattice) $P$ isomorphic to $Q$ | 4, 11 |

| Symbol | Explanation | Page |
|---|---|---|
| **Constructions** | | |
| $P \times Q$ | direct product of $P$ and $Q$ | 6, 18 |
| $P + Q$ | sum of $P$ and $Q$ | 6 |
| $P \dotplus Q$ | glued sum of $P$ and $Q$ | 16 |
| $A[B]$ | tensor extension of $A$ by $B$ | 280 |
| $A \otimes B$ | tensor product of $A$ and $B$ | 277 |
| $U \circledast V$ | modular lattice construction | 153 |
| **Prime intervals** | | |
| $\mathfrak{p}, \mathfrak{q}, \dots$ | prime intervals | |
| $\mathrm{con}(\mathfrak{p})$ | principal congruence generated by $\mathfrak{p}$ | 38 |
| $\mathrm{Prime}(L)$ | set of prime intervals of $L$ | 38 |
| $\mathrm{Princ}\, L$ | the ordered set of principal congruences | 37 |
| **Perspectivities** | | |
| $[a,b] \sim [c,d]$ | $[a,b]$ perspective to $[c,d]$ | 32 |
| $[a,b] \overset{\text{up}}{\sim} [c,d]$ | $[a,b]$ up-perspective to $[c,d]$ | 32 |
| $[a,b] \overset{\text{dn}}{\sim} [c,d]$ | $[a,b]$ down-perspective to $[c,d]$ | 32 |
| $[a,b] \approx [c,d]$ | $[a,b]$ projective to $[c,d]$ | 32 |
| $[a,b] \twoheadrightarrow [c,d]$ | $[a,b]$ congruence-perspective onto $[c,d]$ | 36 |
| $[a,b] \overset{\text{up}}{\twoheadrightarrow} [c,d]$ | $[a,b]$ up congruence-perspective onto $[c,d]$ | 35 |
| $[a,b] \overset{\text{dn}}{\twoheadrightarrow} [c,d]$ | $[a,b]$ down congruence-perspective onto $[c,d]$ | 35 |
| $[a,b] \Rightarrow [c,d]$ | $[a,b]$ congruence-projective onto $[c,d]$ | 35 |
| $\mathfrak{p} \xrightarrow{\text{p}} \mathfrak{q}$ | $\mathfrak{p}$ prime-perspective to $\mathfrak{q}$ | 363 |
| $\mathfrak{p} \xrightarrow{\text{p-up}} \mathfrak{q}$ | $\mathfrak{p}$ prime-perspective up to $\mathfrak{q}$ | 363 |
| $\mathfrak{p} \xrightarrow{\text{p-dn}} \mathfrak{q}$ | $\mathfrak{p}$ prime-perspective down to $\mathfrak{q}$ | 363 |
| $\mathfrak{p} \xRightarrow{\text{p}} \mathfrak{q}$ | $\mathfrak{p}$ prime-projective to $\mathfrak{q}$ | 363 |
| $\mathfrak{p} \curvearrowright \mathfrak{q}$ | $\mathfrak{p}$ swings to $\mathfrak{q}$ | 367 |
| $\mathfrak{p} \overset{\text{in}}{\curvearrowright} \mathfrak{q}$ | internal swing | 367 |
| $\mathfrak{p} \overset{\text{ex}}{\curvearrowright} \mathfrak{q}$ | external swing | 367 |
| **Miscellaneous** | | |
| $\overline{x}$ | closure of $x$ | 10 |
| $\varnothing$ | empty set | 4 |
| $\rceil$ | restriction | 4 |
| **Acronyms** | | |
| SPS lattice | slim, planar, semimodular lattice | 47 |
| SR lattice | slim rectangular lattice | 47 |

# *Picture Gallery*

$C_2$

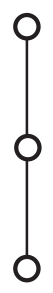

$C_3$

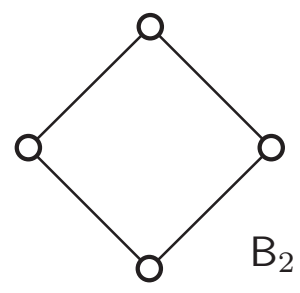

$B_2$

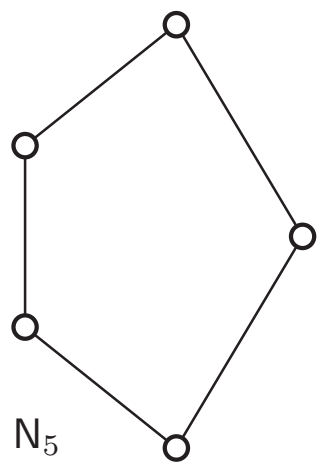

$N_5$

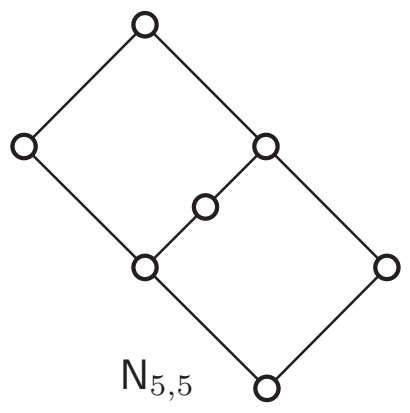

$N_{5,5}$

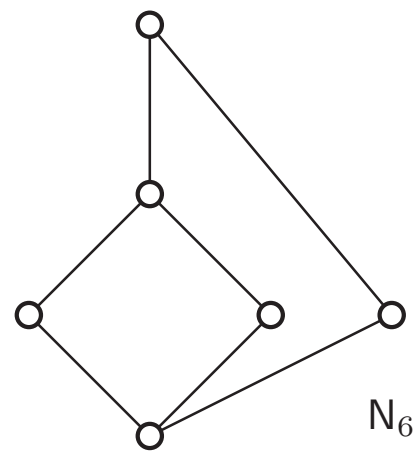

$N_6$

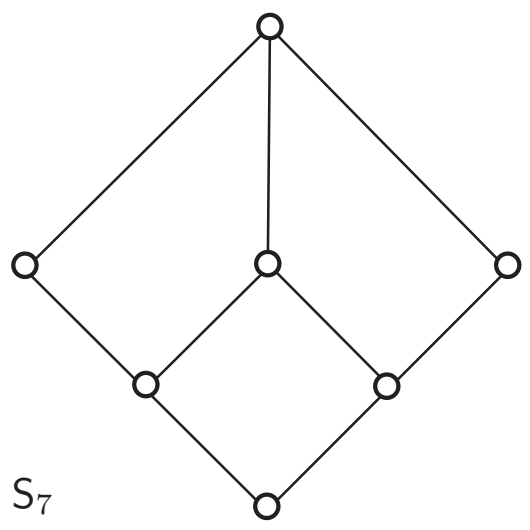

$S_7$

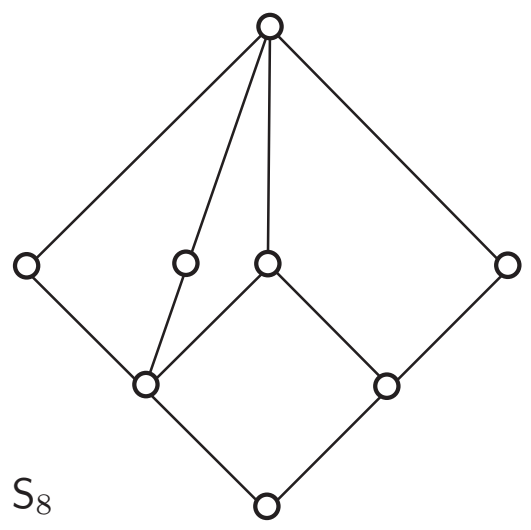

$S_8$

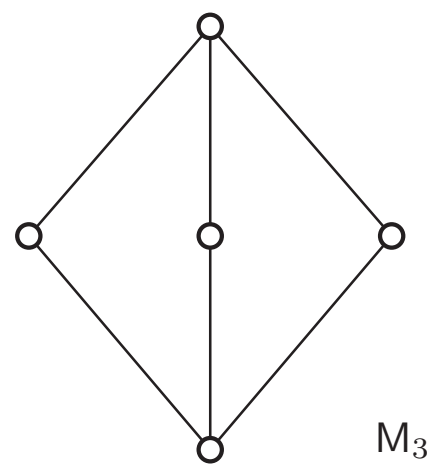

$M_3$

# *Preface*

A book such as this is largely autobiographical; it references about a third of my (mathematical) publications from 1956 to 2022.

Compared with the first edition (in 2006), this book grew from 281 to 432 pages, from 265 to 360 statements, and from 123 to 262 references.

The manuscript was read for the publisher by Friedrich Wehrung; he offered many corrections.

Gábor Czédli, Shriram K. Nimbhorkar, Sylvia Pulmannová, A. Tepavčević sent detailed reports, offering new evidence that Fred's Law is true (a manuscript of $n$ characters has $2^n$ typos; this book has about 900,000 characters).

And after all this, I received a 30-page report from Gregory L. Cherlin, full of corrections and suggestions. I'm not sure how I can thank him for his contribution.

I updated the Introduction and the freely available Part I:

`arXiv:2104.06539`
`https://www.researchgate.net/publication/360184868`

Toronto, Ontario
January 8, 2023

George Grätzer

# *Introduction*

## The topics

This book is an introduction to congruences of finite lattices, which naturally splits into four fields of research:

A. Congruence lattices of finite lattices.

B. The ordered set of principal congruences of finite lattices.

C. The congruence structure of finite lattices.

D. Congruence properties of slim, planar, semimodular (SPS) lattices.

These topics cover about 80 years of research and 250 papers.

Some of the results extend naturally, for instance, to lattices of finite length and to some classes of universal algebras. To keep this book reasonably short, we do not include these results.

### Topic A. Congruence lattices of finite lattices

The congruences of a finite lattice $L$ form a lattice, called the *congruence lattice of* $L$, denoted by $\operatorname{Con} L$. According to a 1942 result of N. Funayama and T. Nakayama [95], the lattice $\operatorname{Con} L$ is a finite *distributive* lattice. The converse is a result of R. P. Dilworth from around 1944 (see the book [22]).

**Basic Representation Theorem.** *Every finite distributive lattice $D$ can be represented as the congruence lattice of a finite lattice $L$.*

We will refer to this result as the *Basic RT* (Representation Theorem). The Basic RT was not published until 1962 in my joint paper [172] with E. T. Schmidt. A large number of papers have strengthened and generalized the Basic RT. These papers form two distinct subfields:

**A1.** RTs of finite distributive lattices as congruence lattices of finite lattices with special properties.

**A2.** The Congruence Lattice Problem (CLP): Can congruence lattices of general lattices be characterized as distributive algebraic lattices?

**A1.** A finite distributive lattice $D$ of more than one element is determined by the ordered set $\mathrm{J}(D)$ of join-irreducible elements. So a representation of $D$ as the congruence lattice of a finite lattice $L$ is the same as a representation of a finite ordered set $P$ $(= \mathrm{J}(D))$, as the ordered set, $\mathrm{J}(\mathrm{Con}\, L)$, of join-irreducible congruences of a finite lattice $L$. A join-irreducible congruence of a finite lattice of more than one element $L$ is exactly the same as a congruence of the form $\mathrm{con}(a, b)$, where $a \prec b$ in $L$; that is, the smallest congruence collapsing a prime interval $[a, b]$. Therefore, it is enough to concentrate on such congruences, and make sure that they are ordered as required by $P$.

**A2.** The infinite case is much different. There are really only two general positive results.

1. The ideal lattice of a distributive lattice with zero is the congruence lattice of a lattice (see E. T. Schmidt [241] and also P. Pudlák [230]).

2. Any distributive algebraic lattice with at most $\aleph_1$ compact elements is the congruence lattice of a lattice (A. P. Huhn [214] and [215], see also H. Dobbertin [80]).

The big breakthrough for negative results came in 2007 in F. Wehrung [258] (based on his paper [255]). Wehrung proved that there is a distributive algebraic lattice with $\aleph_{\omega+1}$ compact elements that is not representable as the congruence lattice of a lattice. P. Růžička [235] improved this result: there is a distributive algebraic lattice with $\aleph_2$ compact elements that is not representable as the congruence lattice of a lattice. This is the sharp bound.

This book deals almost exclusively with the finite case. F. Wehrung [259]–[261] (Chapters 7–9 of the book LTS1-[206]) provide a detailed review of the infinite case.

### The two types of RTs

RTs for the finite case are all of the same general type. We represent a finite distributive lattice $D$ as the congruence lattice of a "nice" finite lattice $L$. For instance, in my joint 1962 paper with E. T. Schmidt [172], we proved that the finite lattice $L$ for the Basic RT can be constructed as a *sectionally complemented lattice*.

To understand the second, more sophisticated, type of RT, we need the concept of a congruence-preserving extension.

Let $L$ be a lattice, and let $K$ be a sublattice of $L$. In general, there is not much connection between the congruence lattices of $L$ and $K$. If they happen to be naturally isomorphic, we call $L$ a *congruence-preserving extension* of $K$. (More formally, see Section 3.3.)

For sectionally complemented lattices, the congruence-preserving extension theorem, ET, was published in my 1999 paper with E. T. Schmidt [184], *Every* *finite lattice $K$ has a finite, sectionally complemented, congruence-preserving extension $L$.* In particular this reduces the result of [172] to the Basic RT.

Reading this statement for the first time, it is difficult to appreciate how much stronger this theorem is than the Basic RT. For a finite distributive lattice $D$, the 1962 theorem provides a finite sectionally complemented lattice $L$ whose congruence lattice is isomorphic to $D$; the 1999 theorem starts with an *arbitrary* finite lattice $K$, and builds a sectionally complemented lattice $L$ extending it with the "same" congruence structure.

#### Topic B. The ordered set of principal congruences of finite lattices

A large part of this book investigates the congruence lattice, $\operatorname{Con} L$, of a finite lattice $L$. But $\operatorname{Con} L$ is not the only interesting congruence construct we can associate with a finite lattice $L$. A newer one, from a decade ago, is $\operatorname{Princ} L$, the ordered set of principal congruences of $L$. We discuss this topic in Part VI.

#### Topic C. The congruence structure of finite lattices

The spreading of a congruence from a prime interval to another prime interval involves intervals of arbitrary size. Can we describe such a spreading with prime intervals only?

We can indeed, by introducing the concept of prime-projectivity (see Chapter 26), and obtaining the Prime-projectivity Lemma (see my paper [112]).

Then in Chapter 27, we develop much sharper forms of this result for SPS (slim, planar, semimodular) lattices. The main result is the Swing Lemma (see my paper [114]), from which we derive many of the known results of G. Czédli and myself concerning congruences of SPS lattices.

#### Topic D. Congruence properties of slim, planar, semimodular (SPS) lattices.

A finite ordered set $P$ satisfies the *Two-Cover Condition*, if any element of $P$ has at most two covers. The *Two-Cover Theorem* (Theorem 27.7) states that the ordered set of join-irreducible congruences of an SPS lattice $L$ has the Two-Cover Condition (see my paper [117]).

This theorem is the start of a new field covered in Part VIII.

## *Proof-by-Picture*

In 1960, trying with E. T. Schmidt to prove the Basic RT (unpublished at the time), we came up with the construction—more or less—as presented in Section 8.2. In 1960, we did not anticipate my 1968 result with H. Lakser [139], establishing that the construction of a chopped lattice solves the problem. So we translated the chopped lattice construction to a closure space, as in Section 8.4, proved that the closed sets form a sectionally complemented lattice $L$, and based on that, we verified that the congruence lattice of $L$ represents the given finite distributive lattice.

When we submitted the paper [172] for publication, it had a three-page section explaining the chopped lattice construction and its translation to a closure space. The referee was strict:

"You cannot have a three-page explanation for a two-page proof."

I believe that in the 50 plus years since the publication of that article, few readers have developed an understanding of the idea behind the published proof.

The referee's dictum is quite in keeping with mathematical tradition and practice. When mathematicians discuss new results, they explain the constructions and the ideas with examples; when these same results are published, the motivation and the examples are largely gone. We publish definitions, constructions, and formal proofs (and conjectures, Paul Erdős would have added).

After Gauss proved one of his famous results, he was not yet ready to publicize it because the proof gave away too much as to how the theorem was discovered. "I have had my results for a long time: but I do not yet know how I am to arrive at them," Gauss is quoted in A. Arber [13].

I try to break with this tradition in this book. In many chapters, after stating the main result, I include a section: *Proof-by-Picture*. This is a misnomer. A *Proof-by-Picture* is not a proof. The Pythagorean Theorem has many well-known *Proofs-by-Picture*—sometimes called "Visual Proofs;" these are really proofs. My *Proof-by-Picture* is an attempt *to convey the idea of the proof.* I trust that if the idea is properly understood, the reader should be able to skip the formal proof, or should at least have less trouble reading it. Think of a *Proof-by-Picture* as a lecture to an informed audience, concluding with "its formal details now you can provide."

## Outline and notation

### Part I. A Brief Introduction to Lattices

In the last paragraph, I call an audience "informed" if they are familiar with the basic concepts and techniques of lattice theory. Part I provides this. It can be downloaded at

`arXiv:submit/4276288`
`https://www.researchgate.net/publication/360184868`

I am quite selective as to what to include. There are few proofs in this part (with a few exceptions) they are easy enough for the readers to work them out on their own. For proofs, lots of exercises, and a more detailed exposition, I refer the reader to my book LTF-[105].

This part has only two sections with novel results.

### Part II. Some Special Techniques

Most of the research in this book deals with RTs; lattices with certain properties are constructed with prescribed congruence structures. The constructions are *ad hoc*. Nevertheless, there are three basic techniques to verify them.

- Chopped lattices, used in almost every chapter in Parts III–V.
- Boolean triples, used in Chapters 12, 14, and 17, and generalized in Chapter 21. also used in some papers that did not make it in this book, for instance, my joint paper with E. T. Schmidt [186].
- Cubic extensions, used in most chapters of Part IV.

These are presented in Part II with proofs.

There are two more basic techniques. *Multi-coloring* is used in several papers; however, it appears in the book only in Chapter 18, so we introduce it there. *Pruning* is utilized in Chapters 13 and 16—it would seem to qualify for Part II; however, there are only concrete uses of pruning, there is no general theory to discuss.

### Part III. RTs

This major part contains the RTs of congruence lattices of finite lattices, requiring only chopped lattices from Part II. I cover the following topics.

- The Basic RT and the RT for *sectionally complemented* lattices in Chapter 8 (my joint paper with E. T. Schmidt [172], P. Crawley and R. P. Dilworth [24]; see also the book [22]). The closure relation of Section 8.4 is generalized in Section 8.5 using the N-relation (called $\underline{D}$-relation or J-relation in the literature) due to J. B. Nation et al.
- *Minimal representations* in Chapter 9; that is, for a given $|\mathrm{J}(D)|$, we minimize the size of $L$ representing the finite distributive lattice $D$ (G. Grätzer, H. Lakser, and E. T. Schmidt [156], G. Grätzer, Rival, and N. Zaguia [168]).
- The *semimodular* RT in Chapter 10 (G. Grätzer, H. Lakser, and E. T. Schmidt [159]).
- The *rectangular* RT (my joint paper with E. Knapp [137]) in Chapter 11.
- The RT for *modular* lattices in Chapter 12 (E. T. Schmidt [238] and my joint paper with E. T. Schmidt [189]); we are forced to represent with a countable lattice $L$, since the congruence lattice of a finite modular lattice is always Boolean.
- The RT for *uniform* lattices (that is, lattices in which any two congruence classes of a congruence are of the same size) in Chapter 13 (G. Grätzer, E. T. Schmidt, and K. Thomsen [195]).

## Part IV. ETs

I present the ETs for the following classes of lattices.

- *Sectionally complemented* lattices in Chapter 14 (my joint paper with E. T. Schmidt [184]).
- *Semimodular* lattices in Chapter 15 (my joint paper with E. T. Schmidt [187]).
- *Isoform* lattices (that is, lattices in which any two congruence classes of a congruence are isomorphic) in Chapter 16 (G. Grätzer, R. W. Quackenbush, and E. T. Schmidt [167]).

These three constructions are based on cubic extensions, introduced in Part II.

Finally, in Chapter 17, I discuss two congruence "destroying" extensions, which we call "magic wands" (my joint paper with E. T. Schmidt [190], G. Grätzer, M. Greenberg, and E. T. Schmidt [132]).

## Part V. Congruence Lattices of Two Related Lattices

What happens if we consider the congruence lattices of two related lattices, such as a lattice and a sublattice? I take up some variants of this question in this part.

Let $L$ be a finite lattice, and let $K$ be a sublattice of $L$. As we discuss in Section 3.3, there is an extension map $\operatorname{ext}\colon \operatorname{Con} K \to \operatorname{Con} L$: for a congruence $\boldsymbol{\alpha}$ of $K$, let the image $\operatorname{ext} \boldsymbol{\alpha}$ be the congruence $\operatorname{con}_L(\boldsymbol{\alpha})$ of $L$ generated by $\boldsymbol{\alpha}$. The map ext is a $\{0\}$-separating join-homomorphism.

Chapter 18 proves the converse, a 1974 result of A. P. Huhn [213]. It is presented in a stronger form due to G. Grätzer, H. Lakser, and E. T. Schmidt [157].

I deal with ideals in Chapter 19. Let $K$ be an ideal of a lattice $L$. Then the restriction map $\operatorname{re}\colon \operatorname{Con} L \to \operatorname{Con} K$ (which assigns to a congruence $\boldsymbol{\alpha}$ of $L$, the restriction $\boldsymbol{\alpha}\rceil_K$ of $\boldsymbol{\alpha}$ to $K$) is a bounded homomorphism (that is, $\{0,1\}$-homomorphism). We prove the corresponding representation theorem for finite lattices, based on my joint paper with H. Lakser [140] (see Theorem 19.1).

*Let $D$ and $E$ be finite distributive lattices. Let $\varphi$ be a bounded homomorphism of $D$ into $E$. Then there exists a finite lattice $L$ and an ideal $I$ of $L$ such that $D \cong \operatorname{Con} L$, $E \cong \operatorname{Con} I$, and $\varphi$ is represented by* re*, the restriction map.*

This is an

abstract/abstract

result. The congruence lattices are given as $D$ and $E$, they are abstract finite distributive lattices, whereas the finite lattices $L$ and $G$ are constructed.

G. Czédli [28] improved on the abstract/abstract representation (see Theorem 20.2), obtaining a

concrete/abstract

result, where $E$ is given as $\operatorname{Con} L$. Earlier, E. T. Schmidt [246] proved this result for injective homomorphisms.

In Chapter 19, we also prove two variants. The first is in my joint paper [147] with H. Lakser, stating that this result also holds for *sectionally complemented* lattices. The second is in another joint paper with H. Lakser [145] stating that this result also holds for *planar* lattices.

Chapter 20 contains sharper forms of Theorem 19.1.

Now let $L$ be a lattice and let $F$ and $G$ be convex sublattices of $L$. How should we map $\operatorname{Con} F$ into $\operatorname{Con} G$? We could try to map a congruence $\boldsymbol{\alpha}$ of the convex sublattice $F$ to $\overline{\boldsymbol{\alpha}}$, the minimal extension of $\boldsymbol{\alpha}$ to $L$, and restrict it to the convex sublattice $G$ as follows,

$$\sigma \colon \boldsymbol{\alpha} \to \overline{\boldsymbol{\alpha}}\rceil G, \qquad \boldsymbol{\alpha} \in \operatorname{Con} F.$$

Minimal extensions do not preserve meets, so $\sigma$ is not a bounded homomorphism of $\operatorname{Con} F$ to $\operatorname{Con} G$, in general. This problem does not arise if we assume that $L$ is a congruence-preserving extension of $F$. This leads us to Theorem 20.3 (see my joint paper with H. Lakser [152]), which is a

concrete/concrete

result.

The final Chapter 21 is a first contribution to the following class of problems. Let $\circledast$ be a construction for finite lattices (that is, if $D$ and $E$ are finite lattices, then so is $D \circledast E$). Find a construction $\circledcirc$ of finite distributive lattices (that is, if $K$ and $L$ are finite distributive lattices, then so is $K \circledcirc L$) satisfying $\operatorname{Con}(K \circledast L) \cong \operatorname{Con} K \circledcirc \operatorname{Con} L$.

If the lattice construction is the direct product, the answer is obvious since $\operatorname{Con}(K \times L) \cong \operatorname{Con} K \times \operatorname{Con} L$.

In Chapter 21, we take up the construction defined as the distributive lattice of all isotone maps from $\operatorname{J}(E)$ to $D$.

In my joint paper [128] with M. Greenberg, we introduced another construction: the *tensor extension*, $A[B]$, for nontrivial finite lattices $A$ and $B$. In Chapter 21, we prove that $\operatorname{Con}(A[B]) \cong (\operatorname{Con} A)[\operatorname{Con} B]$.

## Part VI. The Ordered Set of Principal Congruences

In 2013, I raised the question whether one can associate with a finite lattice $L$ a structure of some of its congruences? We could take the ordered set of the principal congruences generated by prime intervals, but this is just $\operatorname{J}(\operatorname{Con} L)$, which is "equivalent" to $\operatorname{Con} L$ (see Section 2.5.2 for an explanation). I proposed to consider the ordered set $\operatorname{Princ} K$ of principal congruences of a lattice $K$. In Part VI, we state and prove the two major results of this new field.

Chapter 22 contains the first RT for Principal Congruences (Theorem 22.1), characterizing $\operatorname{Princ} K$ of a finite lattice $K$ as a finite, bounded, ordered set (my paper [107]). This chapter also contains the Independence Theorem: for a finite lattice $L$, the two related structures $\operatorname{Princ} L$ and $\operatorname{Aut} L$ are independent (G. Czédli [35], see also my paper [119]).

The second major result is the Minimal RT in Chapter 23, based on my joint paper with H. Lakser [150]. A finite lattice $L$ has a *minimal set of principal congruences*, if all proper principal congruences are join-irreducible. Similarly, for a finite distributive lattice $D$, we call the finite lattice $L$ a minimal representation of $D$, if $D$ and $\operatorname{Con} L$ are isomorphic and $L$ has a minimal set of principal congruences. The Minimal RT states that a finite distributive lattice $D$ has a minimal representation $L$ iff $D$ has at most two dual atoms.

Chapter 24 deals with a related result. Let $D$ be a finite distributive lattice and let $Q \subseteq D$. We call the subset $Q$ of $D$ *principal congruence representable*, if there is a finite lattice $L$ such that $\operatorname{Con} L$ is isomorphic to $D$ and $\operatorname{Princ} L$ corresponds to $Q$ under this isomorphism.

We introduce a simple combinatorial condition, called *chain representability* (see Section 24.1), for a subset $Q$ of a finite distributive lattice $D$ and we prove two results.

- The Necessity Theorem. Let $D$ be a finite distributive lattice and let $Q \subseteq D$. If $Q$ is representable, then it is chain representable.

- The Sufficiency Theorem. Let $D$ be a finite distributive lattice with a join-irreducible unit element. Then $Q \subseteq D$ is representable iff it is chain representable.

### Part VII. Congruence Extensions and Prime Interval

We discuss in Section 4.4, how an SPS lattice can be constructed by inserting *forks* and removing *corners*.

In Chapter 28, we examine the extendibility of congruences to a fork extension. Since by the Structure Theorem for SR Lattices (Theorem 4.11), any SR lattice can be obtained from a grid by inserting forks, this gives us an insight into the congruence lattice of SR lattices. The results are mostly technical, except for the *Two-cover Theorem* we discussed in Topic D.

The spreading of a congruence from a prime interval to another prime interval involves intervals of arbitrary size (as illustrated by Figure 3.2). We would like to describe such a spreading with prime intervals only. We do this in Chapter 26 with the Prime-projectivity Lemma (see my paper [112]).

Chapter 27 sharpens the Prime-projectivity Lemma to the Swing Lemma, a very strong form of the Prime-projectivity Lemma, for slim, planar, and semimodular lattices (see my paper [114]). Almost all results for congruences of slim, planar, and semimodular lattices can be derived from the Swing Lemma.

### Part VIII. Six Congruence Properties of SPS lattices

We introduced this field in Topic D. We state and prove Czédli's four properties (G. Czédli [46] and my paper [125]). The sixth major property is from my joint paper with G. Czédli [59], the proof is from my paper [127].

### Notation

Lattice-theoretic terminology and notation evolved from the three editions of G. Birkhoff's Lattice Theory, [21], by way of my books, UA-[97], GLT-[99], UA2-[100], GLT2-[102], CFL-[103], LTF-[105], CFL2-[115], and R. N. McKenzie, G. F. McNulty, and W. F. Taylor [224], changing quite a bit in the process.

Birkhoff's notation for the congruence lattice and ideal lattice of a lattice changed from $\Theta(L)$ and $I(L)$ to $\operatorname{Con} L$ and $\operatorname{Id} L$ , respectively. The advent of LaTeX promoted the use of operators for lattice constructions. I try to be consistent: I use an operator when a new structure is constructed; so I use $\operatorname{Con} L$, $\operatorname{Id} L$, and $\operatorname{Aut} L$, and so on, without parentheses, unless required for readability, for instance, $\mathrm{J}(D)$ and $\operatorname{Con}(\operatorname{Id} L)$. I use functional notation when sets are constructed, as in $\operatorname{Atom}(L)$ and $\mathrm{J}(a)$. "Generated by" uses the same letters as the corresponding lattice construction, but starting with a lower case letter: $\operatorname{Con} L$ is the congruence lattice of $L$ and $\operatorname{con}(H)$ is the congruence generated by $H$, whereas $\operatorname{Id} L$ is the ideal lattice of $L$ and $\operatorname{id}(H)$ is the ideal generated by $H$.

New concepts introduced in more recent research papers exhibit the usual richness in notation and terminology. I use this opportunity, with the wisdom of hindsight, to make their use more consistent. The reader will often find different notation and terminology when reading the original papers. The detailed Table of Notation and Index may help.

In combinatorial results, I use Landau's big $O$ notation: for the functions $f$ and $g$, we write $f = O(g)$ to mean that $|f| \leq C|g|$ for a suitable constant $C$. Natural numbers start at 1.

In Section 4.2, we introduce acronyms for two classes of lattices: SPS and SR.

We also use acronyms for my books; for instance, LTF for "Lattice Theory: Foundation;" we use them in the form LTF-[105].

Toronto, Ontario George Grätzer
Winter, 2022
Homepage: `http://www.maths.umanitoba.ca/homepages/gratzer.html/`

# Part I

# A Brief Introduction to Lattices

# Chapter 1

# *Basic Concepts*

In this chapter we introduce the most basic order theoretic concepts: ordered sets, lattices, diagrams, and the most basic algebraic concepts: sublattices, congruences, products.

## 1.1. Ordering

### 1.1.1 Ordered sets

A *binary relation* $\varrho$ on a nonempty set $A$ is a subset of $A^2$, that is, a set of ordered pairs $(a, b)$, with $a, b \in A$. For $(a, b) \in \varrho$, we will write $a \mathrel{\varrho} b$ or $a \equiv b \pmod{\varrho}$.

A binary relation $\leq$ on a set $P$ is called an *ordering* if it is *reflexive* ($a \leq a$ for all $a \in P$), *antisymmetric* ($a \leq b$ and $b \leq a$ imply that $a = b$ for all $a, b \in P$), and *transitive* ($a \leq b$ and $b \leq c$ imply that $a \leq c$ for all $a, b, c \in P$). An ordered set $(P, \leq)$ consists of a nonempty set $P$ and an ordering $\leq$.

$a < b$ means that $a \leq b$ and $a \neq b$. We also use the "inverse" relations, $a \geq b$ defined as $b \leq a$ and $a > b$ for $b < a$. If more than one ordering is being considered, we write $\leq_P$ for the ordering of $(P, \leq)$; on the other hand if the ordering is understood, we will say that $P$ (rather than $(P, \leq)$) is an ordered set. An ordered set $P$ is *trivial* if $P$ has only one element.

The elements $a$ and $b$ of the ordered set $P$ are *comparable* if $a \leq b$ or $b \leq a$. Otherwise, $a$ and $b$ are *incomparable*, in notation, $a \parallel b$.

Let $H \subseteq P$ and $a \in P$. Then $a$ is an upper bound of $H$ iff $h \leq a$ for all $h \in H$. An upper bound $a$ of $H$ is the *least upper bound* of $H$ iff $a \leq b$ for any upper bound $b$ of $H$; in this case, we will write $a = \bigvee H$. If $a = \bigvee H$ exists, then it is unique. By definition, $\bigvee \varnothing$ exist ($\varnothing$ is the empty set) iff $P$ has a smallest element, *zero*, denoted by 0. The concepts of *lower bound* and *greatest lower bound* are similarly defined; the latter is denoted by $\bigwedge H$. Note that $\bigwedge \varnothing$ exists iff $P$ has a largest element, *unit*, denoted by 1. A *bounded ordered set* has both 0 and 1. We often denote the 0 and 1 of $P$ by $0_P$ and $1_P$. The notation $\bigvee H$ and $\bigwedge H$ will also be used for families of elements.

The adverb "similarly" (in "similarly defined") in the previous paragraph can be given concrete meaning. Let $(P, \leq)$ be an ordered set. Then $(P, \geq)$ is also an ordered set, called the *dual* of $(P, \leq)$. The dual of the ordered set $P$ will be denoted by $P^\delta$. Now if $\Phi$ is a "statement" about ordered sets, and if we replace all occurrences of $\leq$ by $\geq$ in $\Phi$, then we get the *dual* of $\Phi$.

**Duality Principle for Ordered Sets.** *If a statement $\Phi$ is true for all ordered sets, then its dual is also true for all ordered sets.*

For $a, b \in P$, if $a$ is an upper bound of $\{b\}$, then $a$ is an *upper bound* of $b$. If for all $a, b \in P$, the set $\{a, b\}$ has an upper bound, then the ordered set $P$ is *directed*.

A chain (linearly ordered set, *totally ordered set*) is an ordered set with no incomparable elements. An *antichain* is one in which $a \parallel b$ for all $a \neq b$.

Let $(P, \leq)$ be an ordered set and let $Q$ be a nonempty subset of $P$. Then there is a natural ordering $\leq_Q$ on $Q$ induced by $\leq$: for $a, b \in Q$, let $a \leq_Q b$ iff $a \leq b$; we call $(Q, \leq_Q)$ (or simply, $(Q, \leq)$, or even simpler, $Q$) an ordered subset (or *suborder*) of $(P, \leq)$. We denote this ordered set $(P, \leq)\rceil Q$, the *restriction* of $(P, \leq)$ to $Q$.

A *chain $C$ in an ordered set $P$* is a nonempty subset, which, as a suborder, is a chain. An antichain *$C$ in an ordered set $P$* is a nonempty subset which, as a suborder, is an antichain.

The length of a finite chain $C$, $\operatorname{length} C$, is $|C| - 1$. An ordered set $P$ is said to be *of length $n$* (in symbols, $\operatorname{length} P = n$), where $n$ is a natural number iff there is a chain in $P$ of length $n$ and all chains in $P$ are of length $\leq n$.

The ordered sets $P$ and $Q$ *are isomorphic* (written as $P \cong Q$) and the map $\psi\colon P \to Q$ is an *isomorphism* iff $\psi$ is one-to-one and onto and

$$a \leq b \text{ in } P \quad \text{iff} \quad \psi a \leq \psi b \text{ in } Q.$$

Let $\mathsf{C}_n$ denote the set $\{0, \ldots, n-1\}$ ordered by

$$0 < 1 < 2 < \cdots < n-1.$$

Then $\mathsf{C}_n$ is an $n$-element chain. Observe that $\operatorname{length} \mathsf{C}_n = n - 1$. If $C = \{x_0, \ldots, x_{n-1}\}$ is an $n$-element chain and $x_0 < x_1 < \cdots < x_{n-1}$, then $\psi\colon i \mapsto x_i$ is an isomorphism between $\mathsf{C}_n$ and $C$. Therefore, the $n$-element chain is unique up to isomorphism.

Let $\mathsf{B}_n$ denote the set of all subsets of the set $\{0, \ldots, n-1\}$ ordered by containment. Observe that the ordered set $\mathsf{B}_n$ has $2^n$ elements and $\operatorname{length} \mathsf{B}_n = n$.

In general, for a set $X$, we denote by $\operatorname{Pow} X$ the *power set* of $X$, that is, the set of all subsets of $X$ ordered by set inclusion.

A *quasiordered set* is a nonempty set $Q$ with a binary relation $\leq$ that is reflexive and transitive. Let us define the binary relation $a \approx b$ on $Q$ as $a \leq b$ and $b \leq a$. Then $\approx$ is an equivalence relation. Define the set $P$ as $Q/{\approx}$, and on $P$ define the binary relation $\leq$:

$$a/{\approx} \leq b/{\approx} \quad \text{iff} \quad a \leq b \text{ in } Q.$$

It is easy to see that the definition of $\leq$ on $P$ is well defined and that $P$ is an ordered set. We will call $P$ the *ordered set associated with the quasiordered set* $Q$.

Starting with a binary relation $\prec$ on the set $Q$, we can define the *reflexive-transitive closure* $\leq$ of $\prec$ by the formula: for $a, b \in Q$, let $a \leq b$ iff $a = b$ or if $a = x_0 \prec x_1 \prec \cdots \prec x_n = b$ for elements $x_1, \ldots, x_{n-1} \in Q$. Then $\leq$ is a quasiordering on $Q$. A *cycle* on $Q$ is a sequence $x_1, \ldots, x_n \in Q$ satisfying $x_1 \prec x_2 \prec \cdots \prec x_n \prec x_1$ $(n > 1)$. The quasiordering $\leq$ is an ordering iff there are no cycles.

For an ordered set $P$, call $A \subseteq P$ a *down set* iff $x \in A$ and $y \leq x$ in $P$, imply that $y \in A$. For $H \subseteq P$, there is a smallest down set containing $H$, namely, $\{\, x \mid x \leq h, \text{ for some } h \in H \,\}$; we use the notation $\downarrow H$ for this set. If $H = \{a\}$, we write $\downarrow a$ for $\downarrow \{a\}$. Let $\operatorname{Dn} P$ denote the set of all down sets ordered by set inclusion. If $P$ is an antichain, then $\operatorname{Dn} P \cong \mathsf{B}_n$, where $n = |P|$.

The map $\psi \colon P \to Q$ is an isotone map (resp., antitone map) of the ordered set $P$ into the ordered set $Q$ iff $a \leq b$ in $P$ implies that $\psi a \leq \psi b$ (resp., $\psi a \geq \psi b$) in $Q$. Then $\psi P$ is a suborder of $Q$. Even if $\psi$ is one-to-one, the ordered sets $P$ and $\psi P$ need not be isomorphic. If both $P$ and $Q$ are bounded, then the map $\psi \colon P \to Q$ is *bounded* or a *bounded map*, if it preserves the bounds, that is, $\psi 0_P = 0_Q$ and $\psi 1_P = 1_Q$. Most often, we talk about bounded isotone maps (and bounded homomorphisms, see Section 1.3.1).

### 1.1.2 Diagrams

In the ordered set $P$, the element $a$ *is covered by* $b$ or $b$ *covers* $a$ (written as $a \prec b$ or $b \succ a$) iff $a < b$ and $a < x < b$ for no $x \in P$. The binary relation $\prec$ is called the *covering relation*. The covering relation determines the ordering.

Let $P$ be a finite ordered set. Then $a \leq b$ iff $a = b$ or if there exists a finite sequence of elements $x_1, x_2, \ldots, x_n$ such that

$$a = x_1 \prec x_2 \prec \cdots \prec x_n = b.$$

A diagram of an ordered set $P$ represents the elements with small circles $\circ$; the circles representing two elements $x, y$ are connected by a line segment iff one covers the other; if $x$ is covered by $y$, then the circle representing $x$ is placed lower than the circle representing $y$.

The diagram of a finite ordered set determines the order up to isomorphism.

In a diagram the intersection of two line segments does not indicate an element. A diagram is planar if no two line segments intersect. An ordered set $P$ is *planar* if it has a diagram that is planar. Figure 1.1 shows three diagrams of the same ordered set $P$. Since the third diagram is planar, $P$ is a planar ordered set.

### 1.1.3 Constructions of ordered sets

Given the ordered sets $P$ and $Q$, we can form the *direct product* $P \times Q$, consisting of all ordered pairs $(x_1, x_2)$, with $x_1 \in P$ and $x_2 \in Q$, ordered componentwise, that is, $(x_1, x_2) \leq (y_1, y_2)$ iff $x_1 \leq y_1$ and $x_2 \leq y_2$. Therefore, $(x_1, x_2) \prec (y_1, y_2)$ in $P \times Q$ iff $x_1 \prec y_1$, $x_2 = y_2$ or $x_1 = y_1$, $x_2 \prec y_2$. If $P = Q$, then we write $P^2$ for $P \times Q$. Similarly, we use the notation $P^n$ for $P^{n-1} \times P$ for $n > 2$. Figure 1.2 shows a diagram of $\mathsf{C}_2 \times P$, where $P$ is the ordered set with diagrams in Figure 1.1.

Another often used construction is the (ordinal) *sum* $P + Q$ of $P$ and $Q$, defined on the (disjoint) union $P \cup Q$ and ordered as follows:

$$x \leq y \quad \text{iff} \quad \begin{cases} x \leq_P y & \text{for } x, y \in P; \\ x \leq_Q y & \text{for } x, y \in Q; \\ x \in P,\ y \in Q. \end{cases}$$

Figure 1.3 shows diagrams of $\mathsf{C}_2 + P$ and $P + \mathsf{C}_2$, where $P$ is the ordered set with diagrams in Figure 1.1. In both diagrams, the elements of $\mathsf{C}_2$ are black-filled. Figure 1.3 also shows the diagram of $P \dotplus \mathsf{C}_2$.

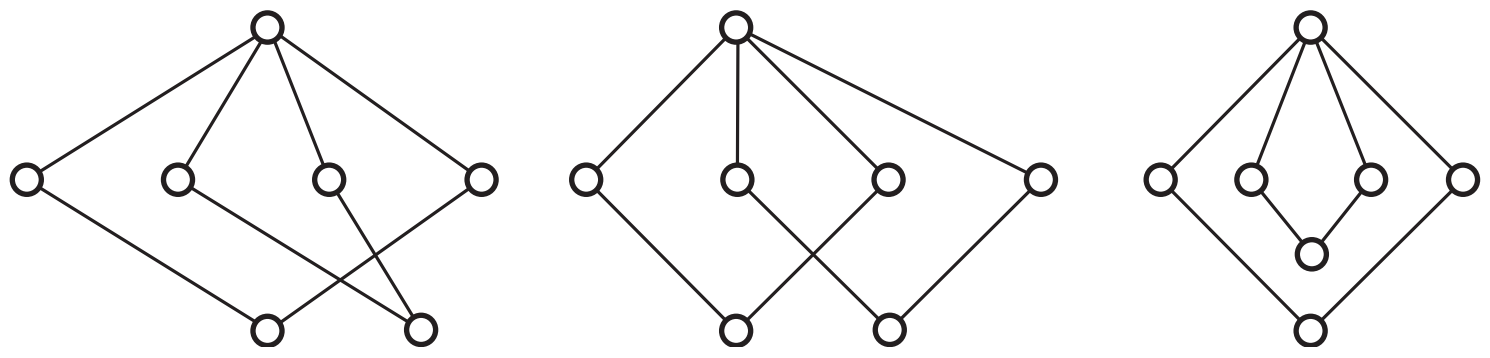

Figure 1.1: Three diagrams of the ordered set $P$

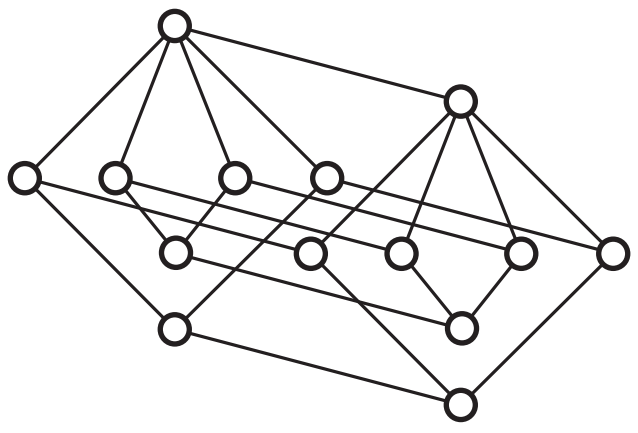

Figure 1.2: A diagram of $\mathsf{C}_2 \times P$

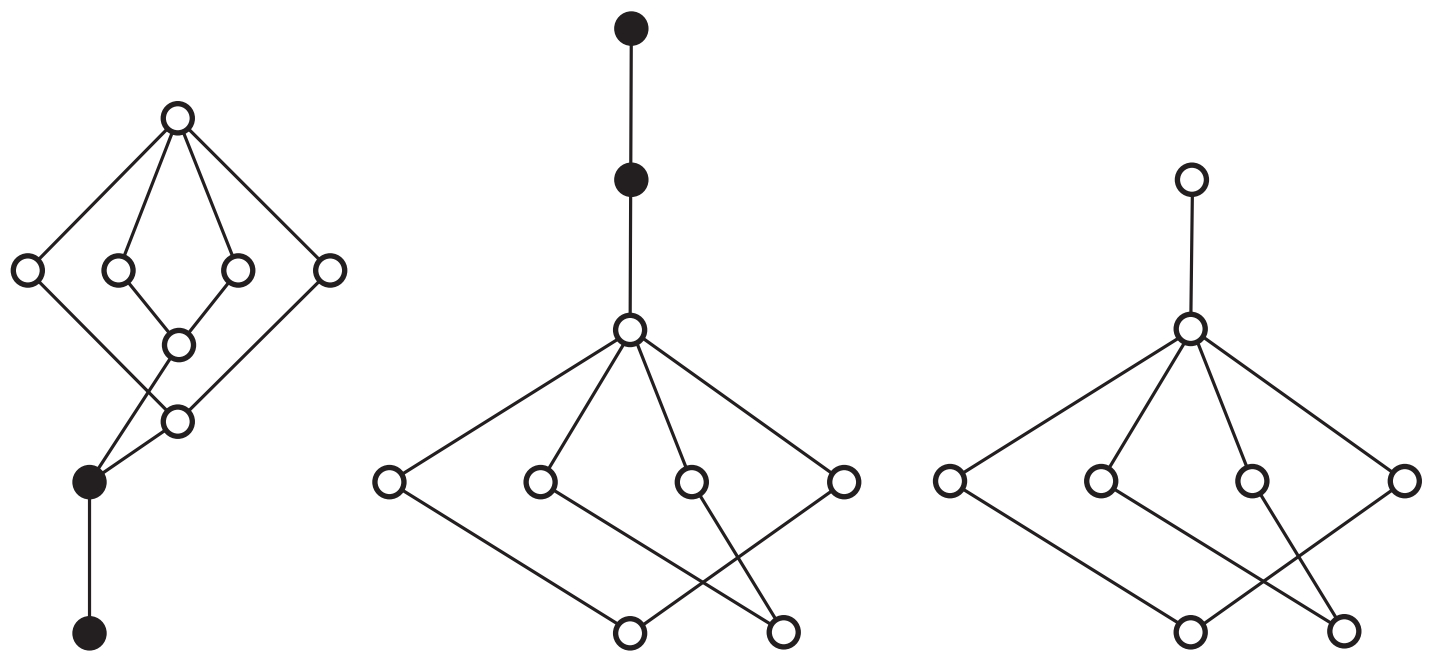

Figure 1.3: Diagrams of $\mathsf{C}_2 + P$, $P + \mathsf{C}_2$, and $P \dotplus \mathsf{C}_2$

A variant construction is the *glued sum*, $P \dotplus Q$, applied to an ordered set $P$ with largest element $1_P$ and an ordered set $Q$ with smallest element $0_Q$; then $P \dotplus Q$ is $P + Q$ in which $1_P$ and $0_Q$ are identified (that is, $1_P = 0_Q$ in $P \dotplus Q$).

### 1.1.4 Partitions

We now give a nontrivial example of an ordered set. A partition of a nonempty set $A$ is a set $\pi$ of nonempty pairwise disjoint subsets of $A$ whose union is $A$. The members of $\pi$ are called the blocks of $\pi$. The block containing $a \in A$ will be denoted by $a/\pi$. A singleton as a block is called trivial. If the elements $a$ and $b$ of $A$ belong to the same block, we write $a \equiv b \pmod{\pi}$ or $a \mathrel{\pi} b$ or $a/\pi = b/\pi$. In general, for $H \subseteq A$,

$$H/\pi = \{\, a/\pi \mid a \in H \,\},$$

a collection of blocks.

An equivalence relation $\varepsilon$ on the set $A$ is a reflexive, symmetric ($a \mathrel{\varepsilon} b$ implies that $b \mathrel{\varepsilon} a$, for all $a, b \in A$), and transitive binary relation. Given a partition $\pi$, we can define an equivalence relation $\varepsilon$ by $(x, y) \in \varepsilon$ iff $x/\pi = y/\pi$. Conversely, if $\varepsilon$ is an equivalence relation, then $\pi = \{\, a/\varepsilon \mid a \in A \,\}$ is a partition of $A$. There is a one-to-one correspondence between partitions and equivalence relations; we will use the two terms interchangeably.

$\operatorname{Part} A$ will denote the set of all partitions of $A$ ordered by

$$\pi_1 \leq \pi_2 \quad \text{iff} \quad x \equiv y \pmod{\pi_1} \text{ implies that } x \equiv y \pmod{\pi_2}.$$

We draw a picture of a partition by drawing the boundary lines of the (nontrivial) blocks. Then $\pi_1 \leq \pi_2$ iff the boundary lines of $\pi_2$ are also boundary lines of $\pi_1$ (but $\pi_1$ may have some more boundary lines). Equivalently, the blocks of $\pi_2$ are unions of blocks of $\pi_1$ (see Figure 1.4).

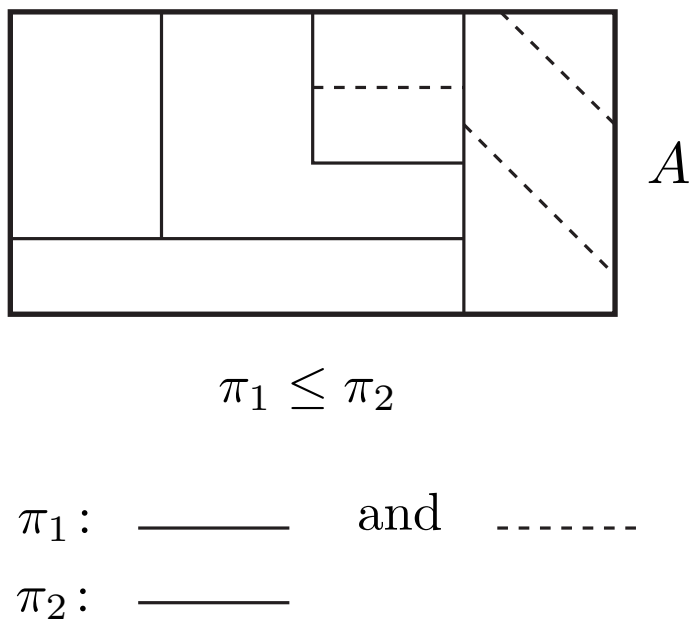

Figure 1.4: Drawing a partition

Part $A$ has a zero and a unit, denoted by **0** and **1**, respectively, defined by

$$x \equiv y \pmod{\mathbf{0}} \qquad \text{iff } x = y;$$
$$x \equiv y \pmod{\mathbf{1}} \qquad \text{for all } x, y \in A.$$

Figure 1.5 shows the diagrams of Part $A$ for $|A| \leq 4$. The partitions are labeled by listing the nontrivial blocks.

## 1.2. Lattices and semilattices

### 1.2.1 Lattices

We need two basic concepts from Universal Algebra. An ($n$-ary) *operation* on a nonempty set $A$ is a map from $A^n$ to $A$. For $n = 2$, we call the operation *binary*. An *algebra* is a nonempty set $A$ with operations defined on $A$.

An ordered set $(L, \leq)$ is a lattice if $\bigvee\{a, b\}$ and $\bigwedge\{a, b\}$ exist for all $a, b \in L$. A lattice $L$ is *trivial* if it has only one element; otherwise, it is *nontrivial*.

We will use the notations

$$a \vee b = \bigvee\{a, b\},$$
$$a \wedge b = \bigwedge\{a, b\},$$

and call $\vee$ the *join*, and $\wedge$ the *meet*. They are both binary operations that are *idempotent* ($a \vee a = a$ and $a \wedge a = a$), *commutative* ($a \vee b = b \vee a$ and $a \wedge b = b \wedge a$), *associative* ($(a \vee b) \vee c = a \vee (b \vee c)$ and $(a \wedge b) \wedge c = a \wedge (b \wedge c)$), and *absorptive* ($a \vee (a \wedge b) = a$ and $a \wedge (a \vee b) = a$). These properties of the operations are also called the *idempotent identities, commutative identities*, *associative identities*, and *absorption identities*, respectively. (Identities, in general, are introduced in Section 2.3.) As always in algebra, associativity makes it possible to write $a_1 \vee a_2 \vee \cdots \vee a_n$ without using parentheses (and the same for $\wedge$).

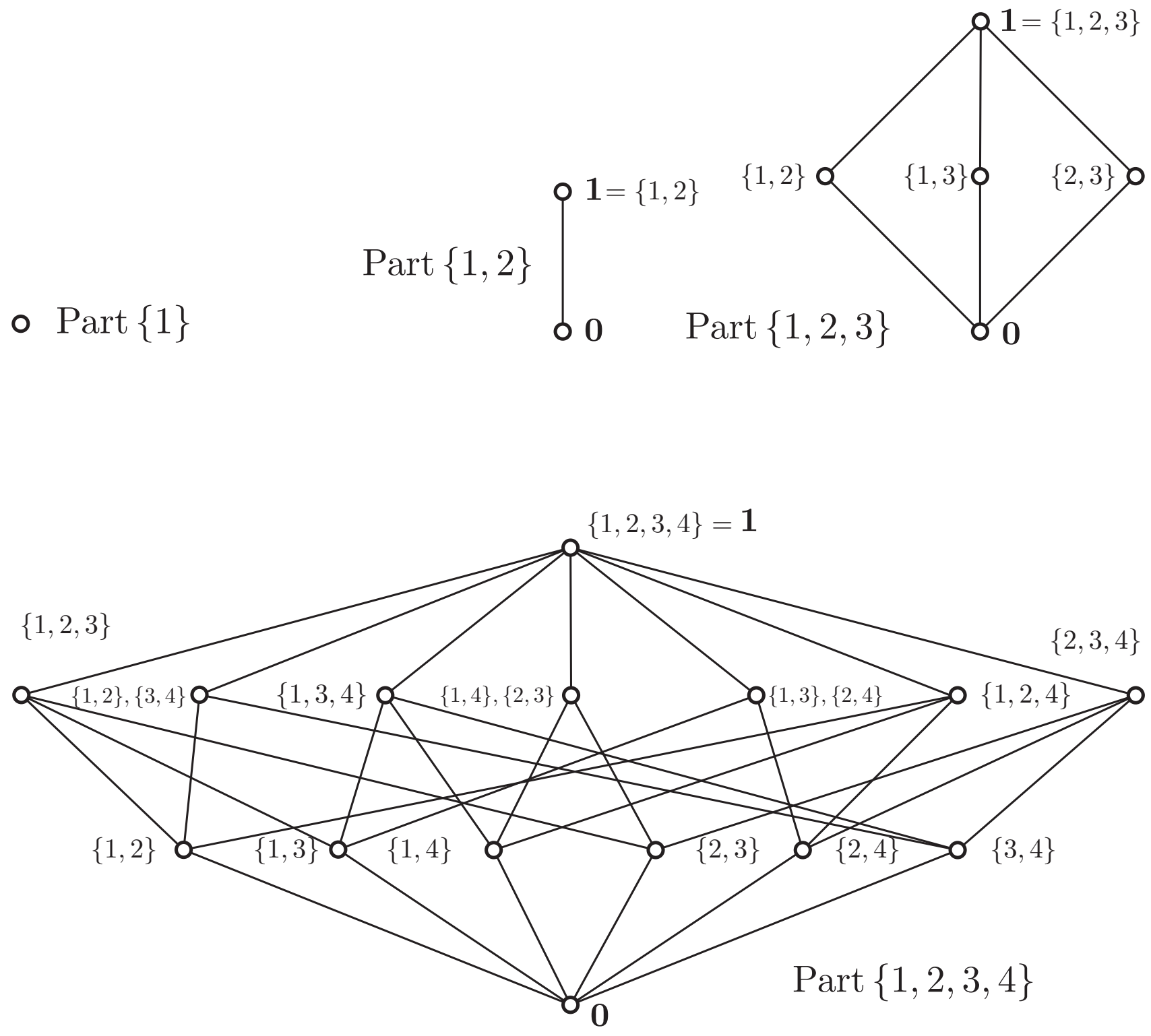

Figure 1.5: Part $A$ for $|A| \leq 4$

For instance, for $A, B \in \operatorname{Pow} X$, we have $A \vee B = A \cup B$ and $A \wedge B = A \cap B$. So $\operatorname{Pow} X$ is a lattice.

For $\boldsymbol{\alpha}, \boldsymbol{\beta} \in \operatorname{Part} A$, if we regard $\boldsymbol{\alpha}$ and $\boldsymbol{\beta}$ as equivalence relations, then the meet formula is trivial: $\boldsymbol{\alpha} \wedge \boldsymbol{\beta} = \boldsymbol{\alpha} \cap \boldsymbol{\beta}$, but the formula for joins is a bit more complicated:

$x \equiv y \pmod{\boldsymbol{\alpha} \vee \boldsymbol{\beta}}$ iff there is a sequence $x = z_0, z_1, \ldots, z_n = y$ of elements of $A$ such that $z_i \equiv z_{i+1} \pmod{\boldsymbol{\alpha}}$ or $z_i \equiv z_{i+1} \pmod{\boldsymbol{\beta}}$ for all $0 \leq i < n$.

So $\operatorname{Part} A$ is a lattice; it is called the *partition lattice* on $A$.

For an ordered set $P$, the order $\operatorname{Dn} P$ is a lattice: $A \vee B = A \cup B$ and $A \wedge B = A \cap B$ for $A, B \in \operatorname{Dn} P$.

To treat lattices as algebras, define an algebra $(L, \vee, \wedge)$ a lattice iff $L$ is a nonempty set, $\vee$ and $\wedge$ are binary operations on $L$, both $\vee$ and $\wedge$ are idempotent, commutative, and associative, and they satisfy the two absorption identities. A lattice as an algebra and a lattice as an ordered set are "equivalent" concepts: Let the order $\mathfrak{L} = (L, \leq)$ be a lattice. Then the algebra $\mathfrak{L}^a = (L, \vee, \wedge)$ is a lattice. Conversely, let the algebra $\mathfrak{L} = (L, \vee, \wedge)$ be a lattice. Define $a \leq b$ iff $a \vee b = b$. Then $\mathfrak{L}^p = (L, \leq)$ is an ordered set, and the ordered set $\mathfrak{L}^p$ is a lattice. For an ordered set $\mathfrak{L}$ that is a lattice, we have $\mathfrak{L}^{ap} = \mathfrak{L}$; for an algebra $\mathfrak{L}$ that is a lattice, we have $\mathfrak{L}^{pa} = \mathfrak{L}$.

Note that for lattices as algebras, the Duality Principle takes on the following very simple form.

**Duality Principle for Lattices.** *Let $\Phi$ be a statement about lattices expressed in terms of $\vee$ and $\wedge$. The dual of $\Phi$ is the statement we get from $\Phi$ by interchanging $\vee$ and $\wedge$. If $\Phi$ is true for all lattices, then the dual of $\Phi$ is also true for all lattices.*

If the operations are understood, we will say that $L$ (rather than $(L, \vee, \wedge)$) is a lattice. The dual of the lattice $L$ will be denoted by $L^\delta$; the ordered set $L^\delta$ is also a lattice. In this book, we deal almost exclusively with finite lattices. Some concepts, however, are more natural to introduce in a more general context. An ordered set $(L, \leq)$ is a *complete lattice* if $\bigvee X$ and $\bigwedge X$ exist for all $X \subseteq L$. All finite lattices are complete, of course.

### 1.2.2 Semilattices and closure systems

A *semilattice* $(S, \circ)$ is an algebra: a nonempty set $S$ with an idempotent, commutative, and associative binary operation $\circ$. A join-semilattice $(S, \vee, \leq)$ is a structure, where $(S, \vee)$ is a semilattice, $(S, \leq)$ is an ordered set, and $a \leq b$ iff $a \vee b = b$. In the ordered set $(S, \leq)$, we have $\bigvee\{a, b\} = a \vee b$. As conventional, we write $(S, \vee)$ for $(S, \vee, \leq)$ or just $S$ if the operation is understood.

Similarly, a *meet-semilattice* $(S, \wedge, \leq)$ is a structure, where $(S, \wedge)$ is a semilattice, $(S, \leq)$ is an ordered set, and $a \leq b$ iff $a \wedge b = a$. In the ordered set $(S, \leq)$, we have $\bigwedge\{a, b\} = a \wedge b$. As conventional, we write $(S, \wedge)$ for $(S, \wedge, \leq)$ or just $S$ if the operation is understood.

If $(L, \vee, \wedge)$ is a lattice, then $(L, \vee)$ is a join-semilattice and $(L, \wedge)$ is a meet-semilattice; moreover, the orderings agree. The converse also holds.

Let $L$ be a lattice and let $C$ be a nonempty subset of $L$ with the property that for every $x \in L$, there is a smallest element $\overline{x}$ of $C$ with $x \leq \overline{x}$. We call $C$ a *closure system* in $L$, and $\overline{x}$ the *closure* of $x$ in $C$.

Obviously, $C$, as an ordered subset of $L$, is a lattice: For $x, y \in C$, the meet in $C$ is the same as the meet in $L$, and the join is

$$x \vee_C y = \overline{x \vee_L y}.$$

Let $L$ be a complete lattice and let $C$ be $\bigwedge$-closed subset of $L$, that is, if $X \subseteq C$, then $\bigwedge X \in C$. (Since $\bigwedge \varnothing = 1$, such a subset is nonempty and contains the 1 of $L$.) Then $C$ is a closure system in $L$, and for every $x \in L$,

$$\overline{x} = \bigwedge(\, y \in C \mid x \leq y \,).$$

# 1.3. Some algebraic concepts

## 1.3.1 Homomorphisms

The lattices $\mathfrak{L}_1 = (L_1, \vee, \wedge)$ and $\mathfrak{L}_2 = (L_2, \vee, \wedge)$ are isomorphic as algebras (in symbols, $\mathfrak{L}_1 \cong \mathfrak{L}_2$), and the map $\varphi\colon L_1 \to L_2$ is an *isomorphism* iff $\varphi$ is one-to-one and onto and

$$\varphi(a \vee b) = \varphi a \vee \varphi b, \tag{1}$$

$$\varphi(a \wedge b) = \varphi a \wedge \varphi b \tag{2}$$

for $a, b \in L_1$.

A map, in general, and a homomorphism, in particular, is called *injective* if it is one-to-one, *surjective* if it is onto, and *bijective* if it is one-to-one and onto.

An isomorphism of a lattice with itself is called an automorphism. The automorphisms of a lattice $L$ form a group $\operatorname{Aut} L$ under composition. A lattice $L$ is *rigid* if the identity map is the only automorphism of $L$, that is, if $\operatorname{Aut} L$ is the one-element group.

It is easy to see that two lattices are isomorphic as ordered sets iff they are isomorphic as algebras.

Let us define a *homomorphism* of the join-semilattice $(S_1, \vee)$ into the join-semilattice $(S_2, \vee)$ as a map $\varphi\colon S_1 \to S_2$ satisfying (1); similarly, for meet-semilattices, we require (2). A *lattice homomorphism* (or simply, *homomorphism*) $\varphi$ of the lattice $L_1$ into the lattice $L_2$ is a map of $L_1$ into $L_2$ satisfying both (1) and (2). A homomorphism of a lattice into itself is called an *endomorphism.* A one-to-one homomorphism is also called an *embedding.*

Note that meet-homomorphisms, join-homomorphisms, and (lattice) homomorphisms are all isotone.

Figure 1.6 shows three maps of the four-element lattice $\mathsf{B}_2$ into the three-element chain $\mathsf{C}_3$. The first map is isotone but it is neither a meet- nor a join-homomorphism. The second map is a join-homomorphism but is not a meet-homomorphism, thus not a homomorphism. The third map is a (lattice) homomorphism.

Various versions of homomorphisms and embeddings will be used. For instance, for lattices and join-semilattices, there are also $\{\vee, 0\}$-homomorphism, and so on, with obvious meanings. An onto homomorphism $\varphi$ is also called *surjective*, whereas a one-to-one homomorphism is called *injective*; it is the same as an *embedding*. For bounded lattices, we often use *bounded homomorphisms* and *bounded embeddings*, that is, $\{0,1\}$-homomorphisms and $\{0,1\}$-embeddings. (In the literature, bounded homomorphisms sometimes have a different definition; this is unlikely to cause any confusion.)

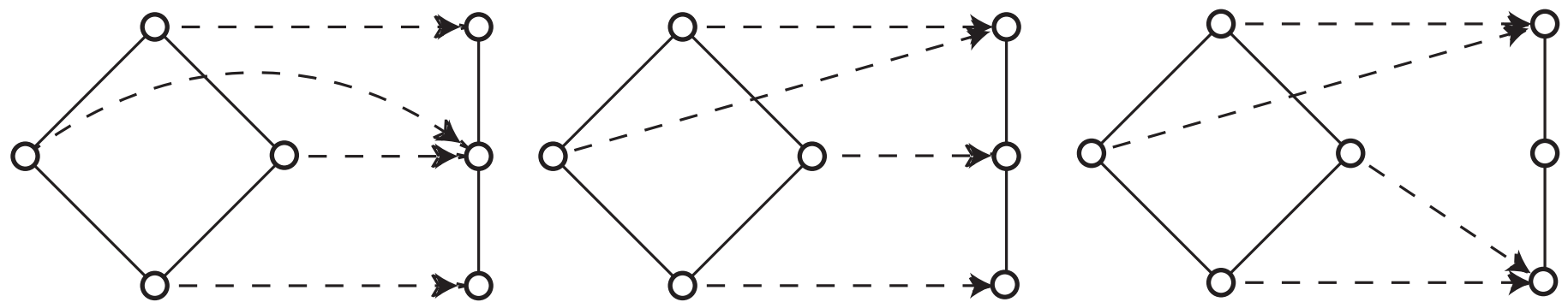

Figure 1.6: Morphism

It should always be clear from the context what kind of homomorphism we are considering. If we say, "let $\varphi$ be a homomorphism of $K$ into $L$," where $K$ and $L$ are lattices, then $\varphi$ is a lattice homomorphism, unless otherwise stated.

### 1.3.2 Sublattices

A *sublattice* $(K, \vee, \wedge)$ of the lattice $(L, \vee, \wedge)$ is defined on a nonempty subset $K$ of $L$ with the property that $a$, $b \in K$ implies that $a \vee b, a \wedge b \in K$ (where the operations $\vee, \wedge$ are formed in $(L, \vee, \wedge)$), and the $\vee$ and the $\wedge$ of $(K, \vee, \wedge)$ are restrictions to $K$ of the $\vee$ and the $\wedge$ of $(L, \vee, \wedge)$, respectively. Instead of "$(K, \vee, \wedge)$ is a sublattice of $(L, \vee, \wedge)$," we will simply say that "$K$ is a sublattice of $L$"—in symbols, $K \leq L$. Of course, a sublattice of a lattice is again a lattice. If $K$ is a sublattice of $L$, then we call $L$ an *extension* of $K$—in symbols, $L \geq K$.

For every $H \subseteq L$, $H \neq \varnothing$, there is a smallest sublattice $\operatorname{sub}(H) \subseteq L$ containing $H$ called the *sublattice of* $L$ *generated by* $H$. We say that $H$ is a *generating set* of $\operatorname{sub}(H)$.

For a bounded lattice $L$, the sublattice $K$ is *bounded* (also called a bounded sublattice) if the 0 and 1 of $L$ are in $K$. Similarly, we can define a $\{0\}$-sublattice, bounded extension, and so on.

The subset $K$ of the lattice $L$ is called *convex* iff $a, b \in K$, $c \in L$, and $a \leq c \leq b$ imply that $c \in K$. We can add the adjective "convex" to sublattices, extensions, and embeddings. A sublattice $K$ of the lattice $L$ is *convex* if it is a convex subset of $L$. Let $L$ be an extension of $K$; then $L$ is a *convex extension* if $K$ is a convex sublattice. An embedding is *convex* if the image is a convex sublattice.

For $a, b \in L$, $a \leq b$, the *interval*

$$I = [a, b] = \{\, x \mid a \leq x \leq b \,\}$$

is an important example of a convex sublattice. We will use the notation $1_I$ for the largest element of $I$, that is, $b$, and $0_I$ for the smallest element of $I$, that is, $a$.

An interval $[a, b]$ is *trivial* if $a = b$. The smallest nontrivial intervals are called *prime*; that is, $[a, b]$ is *prime* iff $a \prec b$. For planar lattices, the term prime interval is used interchangeably with *edge*. In Chapters 27 and 29, we use only edges. Another important example of a convex sublattice is an ideal. A nonempty subset $I$ of $L$ is an *ideal* iff it is a down set with the property:

(Id) $a, b \in I$ implies that $a \vee b \in I$.

An ideal $I$ of $L$ is *proper* if $I \neq L$. Since the intersection of any number of ideals is an ideal, unless empty, we can define $\mathrm{id}(H)$, the *ideal generated by a subset* $H$ of the lattice $L$, provided that $H \neq \varnothing$. If $H = \{a\}$, we write $\mathrm{id}(a)$ for $\mathrm{id}(\{a\})$, and call it a *principal ideal*. Obviously, $\mathrm{id}(a) = \{\, x \mid x \leq a \,\} = \downarrow a$.

The set $\mathrm{Id}\, L$ of all ideals of $L$ is an ordered set under set inclusion, and as an ordered set it is a lattice. In fact, for $I, J \in \mathrm{Id}\, L$, the lattice operations in $\mathrm{Id}\, L$ are $I \vee J = \mathrm{id}(I \cup J)$ and $I \wedge J = I \cap J$. So we obtain the formula for the ideal join:

$x \in I \vee J$ iff $x \leq i \vee j$ for some $i \in I$, $j \in J$.

We call $\mathrm{Id}\, L$ the *ideal lattice* of $L$. Now observe the formulas: $\mathrm{id}(a) \vee \mathrm{id}(b) = \mathrm{id}(a \vee b)$, $\mathrm{id}(a) \wedge \mathrm{id}(b) = \mathrm{id}(a \wedge b)$. Since $a \neq b$ implies that $\mathrm{id}(a) \neq \mathrm{id}(b)$, these yield:

The map $a \mapsto \mathrm{id}(a)$ embeds $L$ into $\mathrm{Id}\, L$.

Since the definition of an ideal uses only $\vee$ and $\leq$, it applies to any join-semilattice $S$. The ordered set $\mathrm{Id}\, S$ is a join-semilattice and the same join formula holds as the one for lattices. Since the intersection of two ideals could be empty, $\mathrm{Id}\, S$ is not a lattice, in general. However, for a $\{\vee, 0\}$-semilattice (a join-semilattice with zero), $\mathrm{Id}\, S$ is a lattice.

For lattices (join-semilattices) $S$ and $T$, let $\varepsilon\colon S \to T$ be an embedding. We call $\varepsilon$ an *ideal-embedding* if $\varepsilon S$ is an ideal of $T$. Then, of course, for any ideal $I$ of $S$, we have that $\varepsilon I$ is an ideal of $T$. Ideal-embeddings play a major role in Chapter 19.

By dualizing, we get the concepts of a *filter*, $\mathrm{fil}(H)$, the *filter generated by a subset* $H$ of the lattice $L$, provided that $H \neq \varnothing$, *principal filter* $\mathrm{fil}(a)$, and so on.

### 1.3.3 Congruences

An equivalence relation $\boldsymbol{\alpha}$ on a lattice $L$ is called a *congruence relation*, or *congruence*, of $L$ iff $a \equiv b \pmod{\boldsymbol{\alpha}}$ and $c \equiv d \pmod{\boldsymbol{\alpha}}$ imply that

$$(\mathrm{SP}_\wedge) \qquad a \wedge c \equiv b \wedge d \pmod{\boldsymbol{\alpha}},$$
$$(\mathrm{SP}_\vee) \qquad a \vee c \equiv b \vee d \pmod{\boldsymbol{\alpha}}$$

(*Substitution Properties*). Trivial examples are the relations $\mathbf{0}$ and $\mathbf{1}$ (introduced in Section 1.1.4). As in Section 1.1.4, for $a \in L$, we write $a/\boldsymbol{\alpha}$ for the *congruence class* (*congruence block*) containing $a$; observe that $a/\boldsymbol{\alpha}$ is a convex sublattice.

If $L$ is a lattice, $K \leq L$, and $\boldsymbol{\alpha}$ a congruence on $L$, then $\boldsymbol{\alpha}\rceil K$, the *restriction of* $\boldsymbol{\alpha}$ *to* $K$, is a congruence of $K$. Formally, for $x, y \in K$,

$$x \equiv y \pmod{\boldsymbol{\alpha}\rceil K} \quad \text{iff} \quad x \equiv y \pmod{\boldsymbol{\alpha}} \text{ in } L.$$

We call $\boldsymbol{\alpha}$ *discrete* on $K$ if $\boldsymbol{\alpha}\rceil K = \mathbf{0}$.

Sometimes it is tedious to compute that a binary relation is a congruence relation. Lemma 1.1, referred to as the Technical Lemma in the literature, often facilitates such computations (see my joint paper [170] with E. T. Schmidt and the paper F. Maeda [222]).

**Lemma 1.1** (Technical Lemma)**.** *A reflexive binary relation $\boldsymbol{\alpha}$ on a lattice $L$ is a congruence relation iff the following three properties are satisfied for $x, y, z, t \in L$:*

(i) $x \equiv y \pmod{\boldsymbol{\alpha}}$ *iff* $x \wedge y \equiv x \vee y \pmod{\boldsymbol{\alpha}}$.

(ii) $x \leq y \leq z$, $x \equiv y \pmod{\boldsymbol{\alpha}}$, *and* $y \equiv z \pmod{\boldsymbol{\alpha}}$ *imply that* $x \equiv z \pmod{\boldsymbol{\alpha}}$.

(iii) $x \leq y$ *and* $x \equiv y \pmod{\boldsymbol{\alpha}}$ *imply that* $x \wedge t \equiv y \wedge t \pmod{\boldsymbol{\alpha}}$ *and* $x \vee t \equiv y \vee t \pmod{\boldsymbol{\alpha}}$.

For finite lattices there is a stronger form (see my paper [111]).

**Lemma 1.2** (Technical Lemma for Finite Lattices)**.** *Let $L$ be a finite lattice. Let $\boldsymbol{\alpha}$ be an equivalence relation on $L$ with intervals as equivalence classes. Then $\boldsymbol{\alpha}$ is a congruence relation iff the following condition and its dual hold for $L$.*

($\mathrm{C}_\vee$) *If* $x \prec y, z \in L$ *and* $x \equiv y \pmod{\boldsymbol{\alpha}}$, *then* $x \equiv y \vee z \pmod{\boldsymbol{\alpha}}$.

Let $\operatorname{Con} L$ denote the set of all congruence relations on $L$ ordered by set inclusion (remember that we can view $\boldsymbol{\alpha} \in \operatorname{Con} L$ as a subset of $L^2$).

We use the Technical Lemma to prove the following result.

**Theorem 1.3.** $\operatorname{Con} L$ *is a lattice. For* $\boldsymbol{\alpha}, \boldsymbol{\beta} \in \operatorname{Con} L$,

$$\boldsymbol{\alpha} \wedge \boldsymbol{\beta} = \boldsymbol{\alpha} \cap \boldsymbol{\beta}.$$

*The join,* $\boldsymbol{\alpha} \vee \boldsymbol{\beta}$, *can be described as follows.*

$x \equiv y \pmod{\boldsymbol{\alpha} \vee \boldsymbol{\beta}}$ *iff there is a sequence*

$$x \wedge y = z_0 \leq z_1 \leq \cdots \leq z_n = x \vee y$$

*of elements of $L$ such that* $z_i \equiv z_{i+1} \pmod{\boldsymbol{\alpha}}$ *or* $z_i \equiv z_{i+1} \pmod{\boldsymbol{\beta}}$ *for every $i$ with* $0 \leq i < n$.

*Remark.* For the binary relations $\gamma$ and $\delta$ on a set $A$, we define the binary relation $\gamma \circ \delta$, the *product of $\gamma$ and $\delta$*, as follows: for $a, b \in A$, the relation $a\,(\gamma \circ \delta)\,b$ holds iff $a \mathrel{\gamma} x$ and $x \mathrel{\delta} b$ for some $x \in A$. The relation $\boldsymbol{\alpha} \vee \boldsymbol{\beta}$ is formed by repeated products. Theorem 1.3 strengthens this statement.

The integer $n$ in Theorem 1.3 can be restricted for some congruence joins. We call the congruences $\boldsymbol{\alpha}$ and $\boldsymbol{\beta}$ *permutable* if $\boldsymbol{\alpha} \vee \boldsymbol{\beta} = \boldsymbol{\alpha} \circ \boldsymbol{\beta}$. A lattice $L$ is *congruence permutable* if any pair of congruences of $L$ are permutable. The chain $\mathsf{C}_n$ is congruence permutable iff $n \leq 2$.

$\operatorname{Con} L$ is called the *congruence lattice* of $L$. Observe that $\operatorname{Con} L$ is a sublattice of $\operatorname{Part} L$; that is, the join and meet of congruence relations as congruence relations and as equivalence relations (partitions) coincide.

If $L$ is nontrivial, then $\operatorname{Con} L$ contains the two-element sublattice $\{\mathbf{0}, \mathbf{1}\}$. If $\operatorname{Con} L = \{\mathbf{0}, \mathbf{1}\}$, we call the lattice $L$ *simple*. All the nontrivial lattices of Figure 1.5 are simple. Of the many lattices of Figure 4.1, only $\mathsf{M}_3$ is simple.

Given $a, b \in L$, there is a smallest congruence $\operatorname{con}(a, b)$—called a *principal congruence*—under which $a \equiv b$. The formula

$$\boldsymbol{\alpha} = \bigvee (\, \operatorname{con}(a, b) \mid a \equiv b \pmod{\boldsymbol{\alpha}}\,) \tag{3}$$

is trivial but important. For $H \subseteq L$, we form the smallest congruence under which $H$ is in one class as $\operatorname{con}(H) = \bigvee(\, \operatorname{con}(a, b) \mid a,\ b \in H\,)$.

Homomorphisms and congruence relations express two sides of the same phenomenon. Let $L$ be a lattice and let $\boldsymbol{\alpha}$ be a congruence relation on $L$. Let $L/\boldsymbol{\alpha} = \{\, a/\boldsymbol{\alpha} \mid a \in L\,\}$. Define $\wedge$ and $\vee$ on $L/\boldsymbol{\alpha}$ by $a/\boldsymbol{\alpha} \wedge b/\boldsymbol{\alpha} = (a \wedge b)/\boldsymbol{\alpha}$ and $a/\boldsymbol{\alpha} \vee b/\boldsymbol{\alpha} = (a \vee b)/\boldsymbol{\alpha}$. The lattice axioms are easily verified. The lattice $L/\boldsymbol{\alpha}$ is the *quotient lattice* of $L$ modulo $\boldsymbol{\alpha}$.

**Lemma 1.4.** *The map*

$$\varphi_{\boldsymbol{\alpha}} \colon x \mapsto x/\boldsymbol{\alpha} \qquad \text{for } x \in L,$$

*is a homomorphism of $L$ onto $L/\boldsymbol{\alpha}$.*

The lattice $K$ is a *homomorphic image* of the lattice $L$ iff there is a homomorphism of $L$ *onto* $K$. Theorem 1.5 (illustrated in Figure 1.7) states that any quotient lattice is a homomorphic image. To state it, we need one more concept: Let $\varphi\colon L \to L_1$ be a homomorphism of the lattice $L$ into the lattice $L_1$, and define the binary relation $\boldsymbol{\alpha}$ on $L$ by $x \mathrel{\boldsymbol{\alpha}} y$ iff $\varphi x = \varphi y$; the relation $\boldsymbol{\alpha}$ is a congruence relation of $L$, called the *kernel* of $\varphi$, in notation, $\ker(\varphi) = \boldsymbol{\alpha}$.

**Theorem 1.5** (Homomorphism Theorem)**.** *Let $L$ be a lattice. Any homomorphic image of $L$ is isomorphic to a suitable quotient lattice of $L$. In fact, if $\varphi\colon L \to L_1$ is a homomorphism of $L$ onto $L_1$ and $\boldsymbol{\alpha}$ is the kernel of $\varphi$, then $L/\boldsymbol{\alpha} \cong L_1$; an isomorphism* (*see Figure 1.7*) *is given by $\boldsymbol{\beta}\colon x/\boldsymbol{\alpha} \mapsto \varphi x$ for $x \in L$.*

We also know the congruence lattice of a homomorphic image.

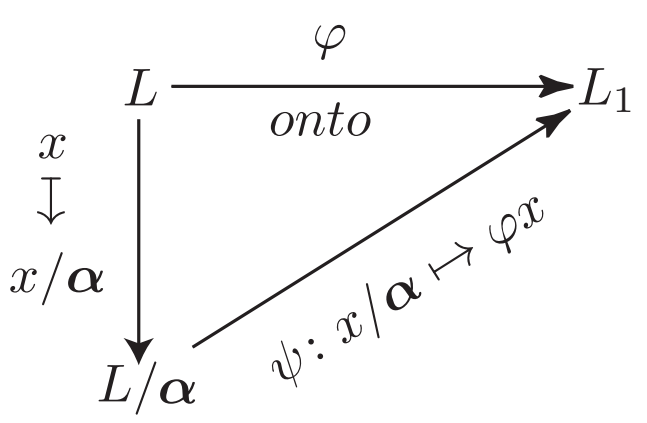

Figure 1.7: The Homomorphism Theorem

**Theorem 1.6** (Second Isomorphism Theorem)**.** *Let $L$ be a lattice and let $\boldsymbol{\alpha}$ be a congruence relation of $L$. For any congruence $\boldsymbol{\beta}$ of $L$ such that $\boldsymbol{\beta} \geq \boldsymbol{\alpha}$, define the relation $\boldsymbol{\beta}/\boldsymbol{\alpha}$ on $L/\boldsymbol{\alpha}$ by*

$$x/\boldsymbol{\alpha} \equiv y/\boldsymbol{\alpha} \pmod{\boldsymbol{\beta}/\boldsymbol{\alpha}} \qquad \textit{iff} \qquad x \equiv y \pmod{\boldsymbol{\beta}}.$$

*Then $\boldsymbol{\beta}/\boldsymbol{\alpha}$ is a congruence of $L/\boldsymbol{\alpha}$. Conversely, every congruence $\boldsymbol{\gamma}$ of $L/\boldsymbol{\alpha}$ can be* (*uniquely*) *represented in the form $\boldsymbol{\gamma} = \boldsymbol{\beta}/\boldsymbol{\alpha}$ for some congruence $\boldsymbol{\beta} \geq \boldsymbol{\alpha}$ of $L$. In particular, the congruence lattice of $L/\boldsymbol{\alpha}$ is isomorphic with the interval* $[\boldsymbol{\alpha}, \mathbf{1}]$ *of the congruence lattice of $L$.*

Let $L$ be a bounded lattice. A congruence $\boldsymbol{\alpha}$ of $L$ *separates* 0 if $0/\boldsymbol{\alpha} = \{0\}$, that is, $x \equiv 0 \pmod{\boldsymbol{\alpha}}$ implies that $x = 0$. Similarly, a congruence $\boldsymbol{\alpha}$ of $L$ *separates* 1 if $1/\boldsymbol{\alpha} = \{1\}$, that is, $x \equiv 1 \pmod{\boldsymbol{\alpha}}$ implies that $x = 1$. We call the lattice $L$ *nonseparating* if there is no congruence $\boldsymbol{\alpha} \neq \mathbf{0}$ separating both 0 and 1 (see Figure 1.8).

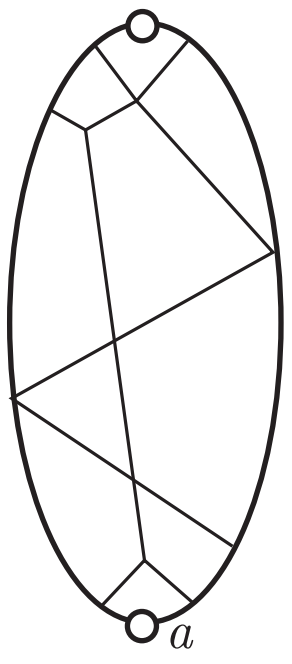

Figure 1.8: Illustrating a separating lattice

Similarly, a homomorphism $\varphi$ of the lattices $L_1$ and $L_2$ with zero is 0-*separating* if $\varphi 0 = 0$, but $\varphi x \neq 0$ for $x \neq 0$. We also use 1-*separating* and $\{0, 1\}$-*separating.*

**Chapter**

# 2

# *Special Concepts*

In this chapter, we introduce special elements, constructions, and classes of lattices that play an important role in the representation of finite distributive lattices as congruence lattices of finite lattices.

## 2.1. Elements and lattices

In a nontrivial finite lattice $L$, an element $a$ is *join-reducible* if $a = 0$ or if $a = b \vee c$ for some $b < a$ and $c < a$; otherwise, it is *join-irreducible*. If $a \in L$ is join-irreducible, then it has a unique lower cover, denoted by $a_*$. Let $\mathrm{J}(L)$ denote the set of all join-irreducible elements of $L$, regarded as an ordered set under the ordering of $L$. By definition, $0 \notin \mathrm{J}(L)$. For $a \in L$, set

$$\mathrm{J}(a) = \{\, x \mid x \leq a,\ x \in \mathrm{J}(L) \,\} = \mathrm{id}(a) \cap \mathrm{J}(L),$$

that is, $\mathrm{J}(a)$ is $\mathrm{id}(a)$ formed in $\mathrm{J}(L)$.

In a finite lattice, every element is a join of join-irreducible elements (indeed, $a = \bigvee \mathrm{J}(a)$), and similarly for meets.

Dually, we define *meet-reducible*, *meet-irreducible* element, its unique upper cover, denoted by $a^*$, and $\mathrm{M}(L)$ the set of all meet-irreducible elements of $L$, regarded as an ordered set under the ordering of $L$. By definition, $1 \notin \mathrm{M}(L)$.

An element $a$ is an *atom* if $0 \prec a$ and a *dual atom* if $a \prec 1$. Atoms are join-irreducible.

A lattice $L$ is *atomistic* if every element is a finite join of atoms.

In a bounded lattice $L$, the element $a$ is a *complement* of the element $b$ iff $a \wedge b = 0$ and $a \vee b = 1$. A *complemented lattice* is a bounded lattice in which every element has a complement. The lattices of Figure 1.5 are complemented and so are all but one of the lattices of Figure 4.1. The lattice $\mathsf{B}_n$ is also complemented.

Let $a \in [b, c]$; the element $x$ is a *relative complement of a in* $[b, c]$ iff $a \wedge x = b$, and $a \vee x = c$. A *relatively complemented lattice* is a lattice in which every element has a relative complement in any interval containing it. The lattice $\mathsf{N}_5$ of the Picture Gallery is complemented but not relatively complemented.

In a lattice $L$ with zero, let $a \leq b$. A complement of $a$ in $[0, b]$ is called a *sectional complement* of $a$ in $b$. A lattice $L$ with zero is called *sectionally complemented* if $a$ has a sectional complement in $b$ for all $a \leq b$ in $L$. The lattice $\mathsf{N}_6$ of the Picture Gallery is sectionally complemented but not relatively complemented.

## 2.2. Direct and subdirect products

Let $L$ and $K$ be lattices and form the *direct product* $L \times K$ as in Section 1.1.3. Then $L \times K$ is a lattice and $\vee$ and $\wedge$ are computed "componentwise:;"

$$(a, b) \vee (c, d) = (a \vee c, b \vee d),$$
$$(a, b) \wedge (c, d) = (a \wedge c, b \wedge d).$$

Figure 2.1 shows the example $\mathsf{C}_2 \times \mathsf{N}_5$.

Obviously, $\mathsf{B}_n$ is isomorphic to a direct product of $n$ copies of $\mathsf{B}_1$.

There are two *projection maps* (homomorphisms) associated with this construction:

$$\pi_L \colon L \times K \to L \quad \text{and} \quad \pi_K \colon L \times K \to K,$$

defined by $\pi_L \colon (x, y) \mapsto x$ and by $\pi_K \colon (x, y) \mapsto y$, respectively.

Similarly, we can form the *direct product* $L_1 \times \cdots \times L_n$ with elements $(x_1, \ldots, x_n)$, where $x_i \in L_i$ for $i \leq n$; we denote the projection map

$$(x_1, \ldots, x_i, \ldots, x_n) \mapsto x_i$$

by $\pi_i$. If $L_i = L$ for all $i \leq n$, we get the *direct power* $L^n$. By identifying $x \in L_i$ with $(0, \ldots, 0, x, 0, \ldots, 0)$ ($x$ is the $i$-th coordinate), we regard $L_i$ as an ideal of $L_1 \times \cdots \times L_n$ for $i \leq n$. The black-filled elements in Figure 2.1 show how we consider $\mathsf{C}_2$ and $\mathsf{N}_5$ ideals of $\mathsf{C}_2 \times \mathsf{N}_5$.

A very important property of direct products is:

**Theorem 2.1.** *Let $L$ and $K$ be lattices, let $\boldsymbol{\alpha}$ be a congruence relation of $L$, and let $\boldsymbol{\beta}$ be a congruence relation of $K$. Define the relation $\boldsymbol{\alpha} \times \boldsymbol{\beta}$ on $L \times K$ by*

$$(a_1, b_1) \equiv (a_2, b_2)\,(\boldsymbol{\alpha} \times \boldsymbol{\beta}) \quad \textit{iff} \quad a_1 \equiv a_2 \pmod{\boldsymbol{\alpha}} \textit{ and } b_1 \equiv b_2 \pmod{\boldsymbol{\beta}}.$$

*Then $\boldsymbol{\alpha} \times \boldsymbol{\beta}$ is a congruence relation on $L \times K$. Conversely, every congruence relation of $L \times K$ is of this form.*

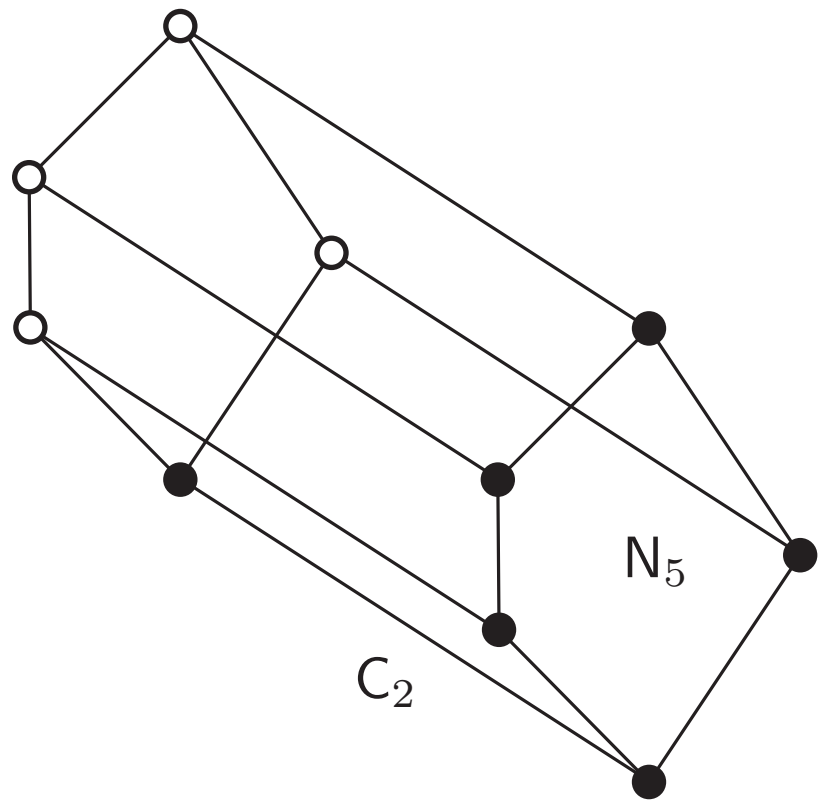

Figure 2.1: $\mathsf{C}_2 \times \mathsf{N}_5$, a direct product of two lattices

A more general construction is subdirect products. If $L \leq K_1 \times \cdots \times K_n$ and the projection maps $\pi_i$ are onto maps for all $i \leq n$, then we call $L$ a *subdirect product* of $K_1, \ldots, K_n$.

Trivial examples: $L$ is a subdirect product of $L$ and $L$ if we identify $x \in L$ with $(x, x) \in L^2$ (diagonal embedding). In this example, the projection map is an isomorphism. To exclude such trivial cases, let us call a representation of $L$ as a subdirect product of $K_1, \ldots, K_n$ *trivial* if one of the projection maps $\pi_1, \ldots, \pi_n$ is an isomorphism.

A lattice $L$ is called *subdirectly irreducible* iff all representations of $L$ as a subdirect product are trivial. For a subdirect product $L$ of $K_1$ and $K_2$, we have $\ker(\pi_1) \wedge \ker(\pi_2) = \mathbf{0}$. This subdirect product is trivial iff $\ker(\pi_1) = \mathbf{0}$ or $\ker(\pi_2) = \mathbf{0}$. Conversely, if $\boldsymbol{\alpha}_1 \wedge \boldsymbol{\alpha}_2 = \mathbf{0}$ in $\operatorname{Con} K$, then $K$ is a subdirect product of $K/\boldsymbol{\alpha}_1$ and $K/\boldsymbol{\alpha}_2$, and this representation is trivial iff $\boldsymbol{\alpha}_1 = \mathbf{0}$ or $\boldsymbol{\alpha}_2 = \mathbf{0}$.

Every simple lattice is subdirectly irreducible. The lattice $\mathsf{N}_5$ is subdirectly irreducible but not simple.

There is a natural correspondence between subdirect representations of a lattice $L$ and sets of congruences $\{\boldsymbol{\gamma}_1, \ldots, \boldsymbol{\gamma}_n\}$ satisfying $\boldsymbol{\gamma}_1 \wedge \cdots \wedge \boldsymbol{\gamma}_n = \mathbf{0}$. This representation is nontrivial iff $\boldsymbol{\gamma}_i \neq \mathbf{0}$ for all $i \leq n$. In this subdirect representation, the factors (the lattices $L/\boldsymbol{\gamma}_i$) are subdirectly irreducible iff the congruences $\boldsymbol{\gamma}_i$ are meet-irreducible for all $i \leq n$, by the Second Isomorphism Theorem (Theorem 1.6).

For a finite lattice $L$, the lattice $\operatorname{Con} L$ is finite, so we can represent $\mathbf{0}$ as a meet of meet-irreducible congruences, and we obtain *Birkhoff Subdirect RT*.

**Theorem 2.2.** *Every finite lattice $L$ is a subdirect product of subdirectly irreducible lattices.*

This result is true for any algebra in any variety (a class of algebras defined by identities, such as the class of all lattices or the class of all groups).

Finite subdirectly irreducible lattices are easy to recognize. If $L$ is such a lattice, then the meet $\boldsymbol{\alpha}$ of all the $> \mathbf{0}$ elements is $> \mathbf{0}$. Obviously, $\boldsymbol{\alpha}$ is an atom, the unique atom of $\operatorname{Con} L$. Conversely, if $\operatorname{Con} L$ has a unique atom, then all $> \mathbf{0}$ congruences are $\geq \boldsymbol{\alpha}$, so their meet cannot be $\mathbf{0}$. We call $\boldsymbol{\alpha}$ the *base congruence* of $L$ (called *monolith* in many publications).

If $u \neq v$ and $u \equiv v \pmod{\boldsymbol{\alpha}}$, then $\boldsymbol{\alpha} = \operatorname{con}(u, v)$. So

$$\operatorname{Con} L = \{\mathbf{0}\} \cup [\operatorname{con}(u, v), \mathbf{1}],$$

as illustrated in Figure 2.2.

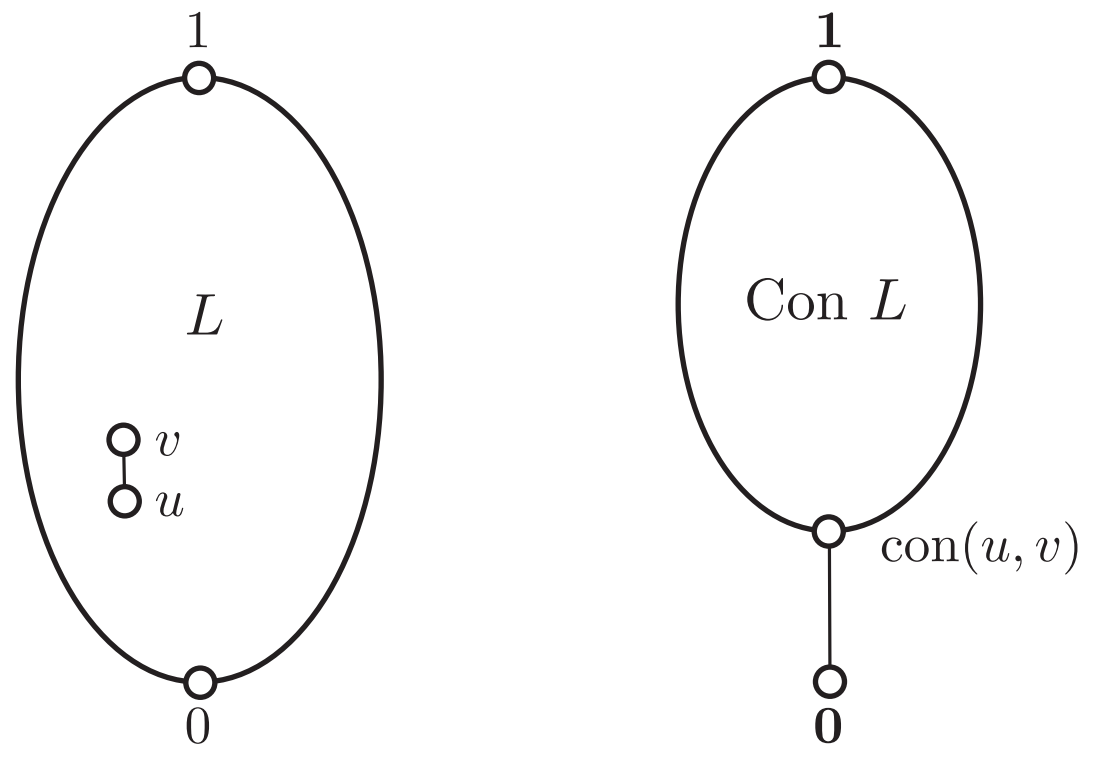

Figure 2.2: A subdirectly irreducible lattice and its congruence lattice

Let $L$ be a finite subdirectly irreducible lattice with $\operatorname{con}(u, v)$ the base congruence, where $u \prec v \in L$. By inserting two elements as shown in Figure 2.3, we embed $L$ into a simple lattice.

**Lemma 2.3.** *Every finite subdirectly irreducible lattice can be embedded into a simple lattice with at most two extra elements.*

Note that every finite lattice can be embedded into a finite simple lattice; in general, we need more than two elements. Lemma 14.3 provides a stronger statement.

## 2.3. Terms and identities

From the variables $x_1, \ldots, x_n$, we can form ($n$-ary) *terms* in the usual manner using $\vee, \wedge$, and parentheses. Examples of terms are: $x_1$, $x_3$, $x_1 \vee x_1$, $(x_1 \wedge x_2) \vee (x_3 \wedge x_1)$, $(x_3 \wedge x_1) \vee ((x_3 \vee x_2) \wedge (x_1 \vee x_2))$.

An $n$-ary term $p$ defines a function in $n$ variables (a *term function*, or simply, a *term*) on a lattice $L$. For example, if

$$p = (x_1 \wedge x_3) \vee (x_3 \vee x_2)$$

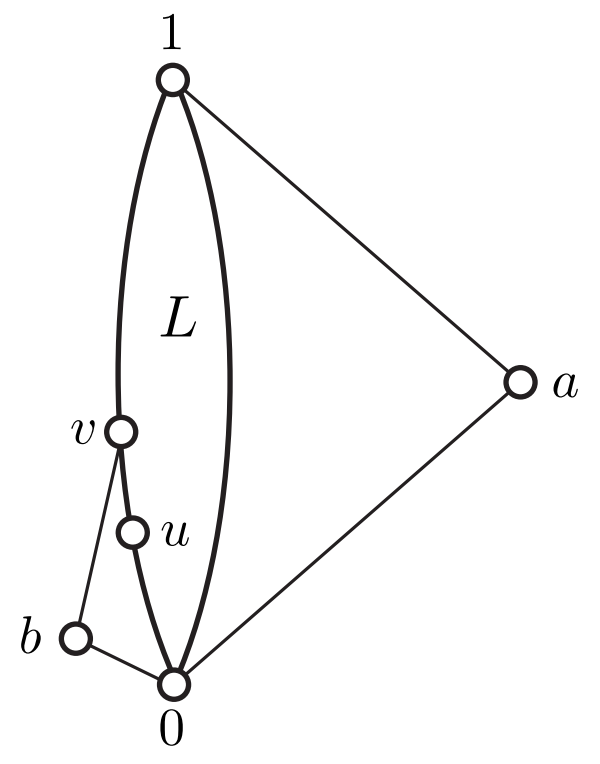

Figure 2.3: Embedding into a simple lattice

and $a, b, c \in L$, then

$$p(a, b, c) = (a \wedge c) \vee (c \vee b) = b \vee c.$$

If $p$ is a *unary* ($n = 1$) lattice term, then $p(a) = a$ for any $a \in L$. If $p$ is *binary*, then $p(a, b) = a$, or $p(a, b) = b$, or $p(a, b) = a \wedge b$, or $p(a, b) = a \vee b$ for all $a$, $b \in L$.

If $p = p(x_1, \dots, x_n)$ is an $n$-ary term and $L$ is a lattice, then by substituting some variables by elements of $L$, we get a function on $L$ of $n$-variables. Such functions are called *term functions*. Unary term functions of the form

$$p(x) = p(a_1, \dots, a_{i-1}, x, a_{i+1}, \dots, x_n),$$

where $a_1, \dots, a_{i-1}, a_{i+1}, \dots, a_n \in L$, play the most important role (see Section 3.1).

A term (function), in fact, any term function, $p$ is *isotone*; that is, if $a_1 \leq b_1$, …, $a_n \leq b_n$, then $p(a_1, \dots, a_n) \leq p(b_1, \dots, b_n)$. Furthermore,

$$a_1 \wedge \cdots \wedge a_n \leq p(a_1, \dots, a_n) \leq a_1 \vee \cdots \vee a_n.$$

Note that many earlier publications in Lattice Theory and Universal Algebra use *polynomials* and *polynomial functions* for terms and term functions.

Terms have many uses. We briefly discuss three.

### (i) The sublattice generated by a set

**Lemma 2.4.** *Let $L$ be a lattice and let $H$ be a nonempty subset of $L$. Then $a \in \operatorname{sub}(H)$ (the sublattice generated by $H$) iff $a = p(h_1, \dots, h_n)$ for some integer $n \geq 1$, for some $n$-ary term $p$, and for some $h_1$, …, $h_n \in H$.*

**(ii) Identities**

A *lattice identity* (resp., *lattice inequality*)—also called *equation*—is an expression of the form $p = q$ (or $p \leq q$), where $p$ and $q$ are terms. An *identity* $p = q$ (or $p \leq q$) *holds in the lattice* $L$ iff $p(a_1, \ldots, a_n) = q(a_1, \ldots, a_n)$ (or $p(a_1, \ldots, a_n) \leq q(a_1, \ldots, a_n)$) holds for all $a_1, \ldots, a_n \in L$. The identity $p = q$ is equivalent to the two inequalities $p \leq q$ and $q \leq p$; the inequality $p \leq q$ is equivalent to the identity $p \vee q = q$.

The most important properties of identities are given by

**Lemma 2.5.** *Identities are preserved under the formation of sublattices, homomorphic images, direct products, and ideal lattices.*

A lattice $L$ is called *distributive* if the identities

$$x \wedge (y \vee z) = (x \wedge y) \vee (x \wedge z),$$
$$x \vee (y \wedge z) = (x \vee y) \wedge (x \vee z)$$

hold in $L$. In fact, it is enough to assume one of these identities, because the two identities are equivalent.

As we have just noted, the identity $x \wedge (y \vee z) = (x \wedge y) \vee (x \wedge z)$ is equivalent to the two inequalities:

$$x \wedge (y \vee z) \leq (x \wedge y) \vee (x \wedge z),$$
$$(x \wedge y) \vee (x \wedge z) \leq x \wedge (y \vee z).$$

However, the second inequality holds in any lattice. So a lattice is distributive iff the inequality $x \wedge (y \vee z) \leq (x \wedge y) \vee (x \wedge z)$ holds. By duality, we get a similar statement about the second identity defining distributivity.

The class of all distributive lattices will be denoted by $\mathbf{D}$. A *Boolean lattice* is a distributive complemented lattice. A finite Boolean lattice is isomorphic to some $\mathsf{B}_n$.

A lattice is called *modular* if the identity

$$(x \wedge y) \vee (x \wedge z) = x \wedge (y \vee (x \wedge z))$$

holds. Note that this identity is equivalent to the following implication:

$$x \geq z \text{ implies that } (x \wedge y) \vee z = x \wedge (y \vee z).$$

The class of all modular lattices will be denoted by $\mathbf{M}$.

Every distributive lattice is modular. The lattice $\mathsf{M}_3$ is modular but not distributive. All the lattices of Figures 1.5 and 4.1 are modular except for $\operatorname{Part}\{1, 2, 3, 4\}$ and $\mathsf{N}_5$.

A class of lattices $\mathbf{V}$ is called a *variety* if it is defined by a set of identities. The classes $\mathbf{D}$ and $\mathbf{M}$ are examples of varieties, and so is the class $\mathbf{L}$, the variety of all lattices, and $\mathbf{T}$, the (trivial) variety of one-element lattices.

### (iii) Free lattices

Starting with a set $H$, we can form the set of all terms over $H$, collapsing two terms if their equality follows from the lattice axioms. We thus form the "free-est" lattice over $H$. For instance, if we start with $H = \{a, b\}$, then we obtain the four-element lattice, $\mathsf{F}(2)$, of Figure 2.4.

We obtain more interesting examples if we start with an ordered set $P$, and require that the ordering in $P$ be preserved. For instance, if we start with $P = \{a, b, c\}$ with $a < b$, then we get the corresponding nine-element "free" lattice, $\mathsf{F}(P)$, of Figure 2.4.[1]

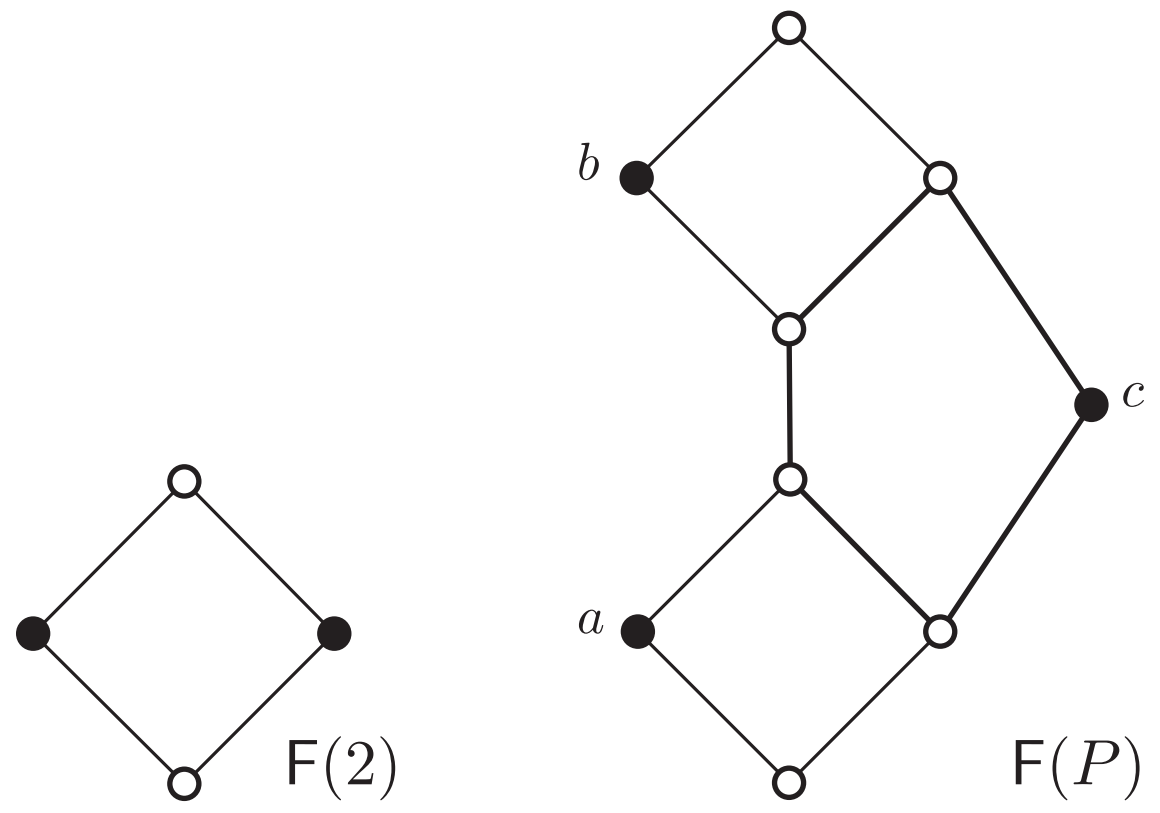

Figure 2.4: Two free lattices

Sometimes, we need free lattices with respect to some special conditions. The following result illustrates this.

**Lemma 2.6.** *Let $x, y$, and $z$ be elements of a lattice $L$ and let $x \vee y$, $y \vee z$, and $z \vee x$ be pairwise incomparable. Then* $\mathrm{sub}(\{x \vee y, y \vee z, z \vee x\}) \cong \mathsf{B}_3$.

Lemma 2.6 is illustrated by Figure 2.5.

We will also need "free distributive lattices," obtained by collapsing two terms if their equality follows from the lattice axioms and the distributive identities. Starting with a three-element set $H = \{x, y, z\}$, we then obtain the lattice, $\mathrm{Free}_{\mathbf{D}}(3)$, of Figure 2.6.

Similarly, we can define "free modular lattices." Starting with a three-element set $H = \{x, y, z\}$, we then obtain the lattice, $\mathrm{Free}_{\mathbf{M}}(3)$, of Figure 2.7.

An equivalent definition of freeness is the following.

Let $H$ be a set and let $\mathbf{K}$ be a variety of lattices. A lattice $\mathrm{Free}_{\mathbf{K}}(H)$ is called a *free lattice* over $\mathbf{K}$ *generated by* $H$ iff the following three conditions are satisfied:

[1]What we call "free", is called *finitely presented* in the literature.

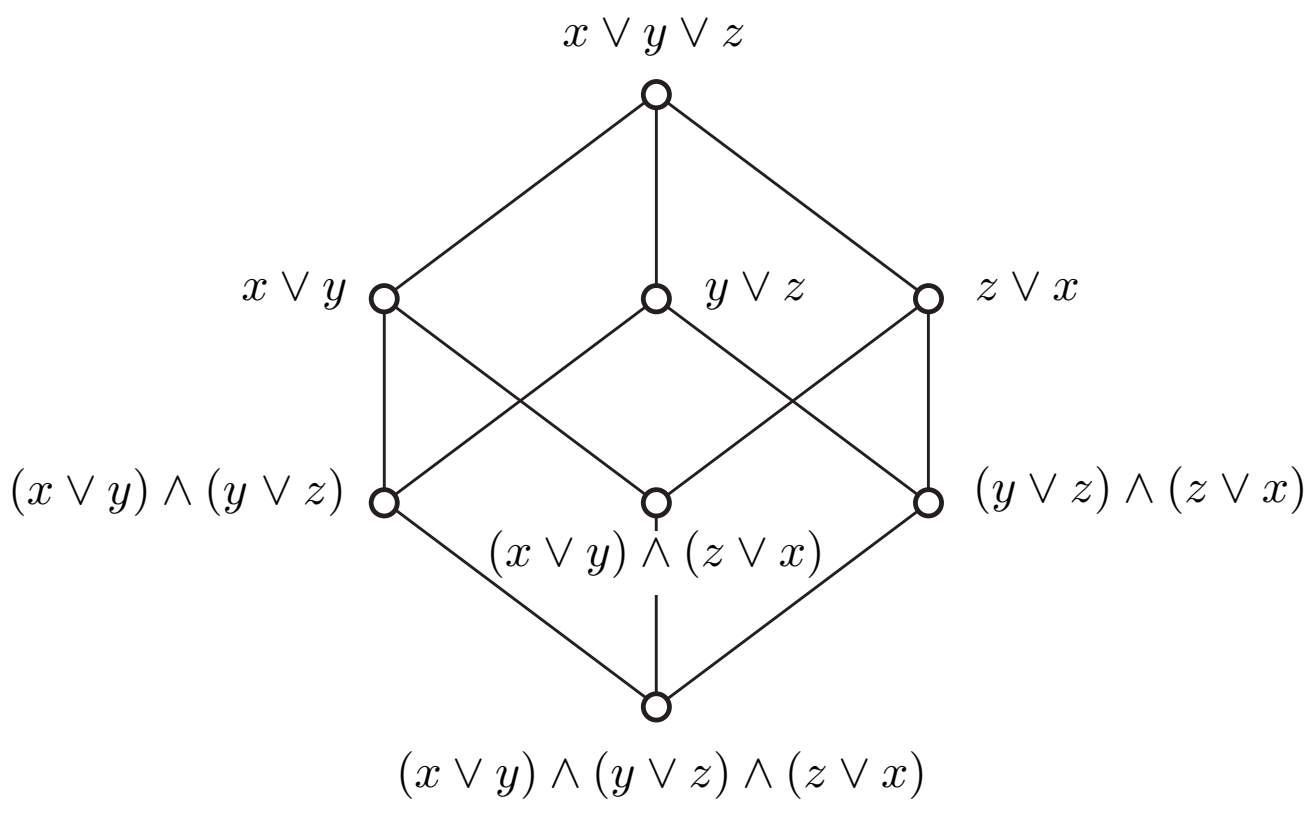

Figure 2.5: A free lattice with special relations

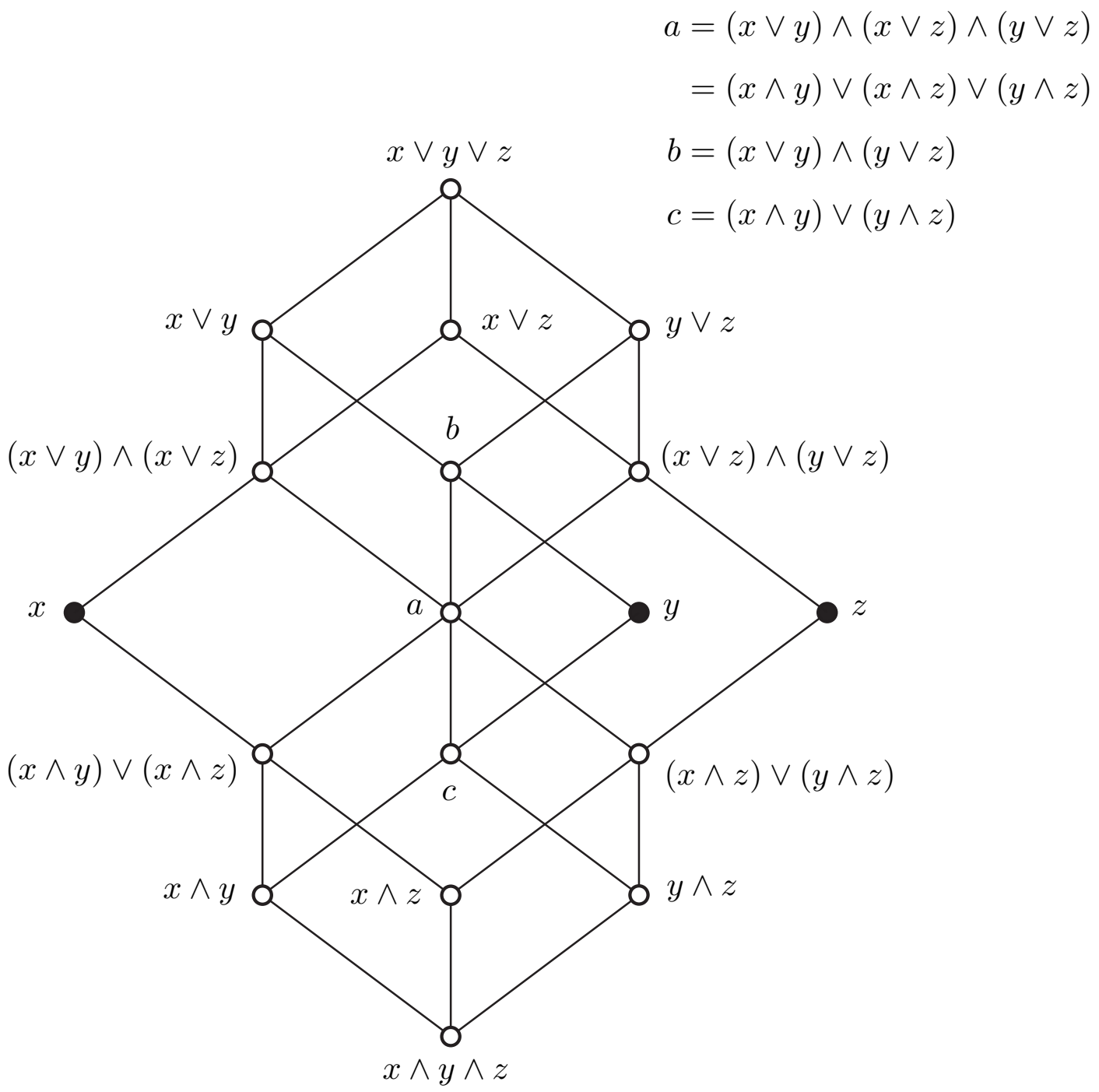

Figure 2.6: The free distributive lattice on three generators, $\mathrm{Free}_{\mathbf{D}}(3)$

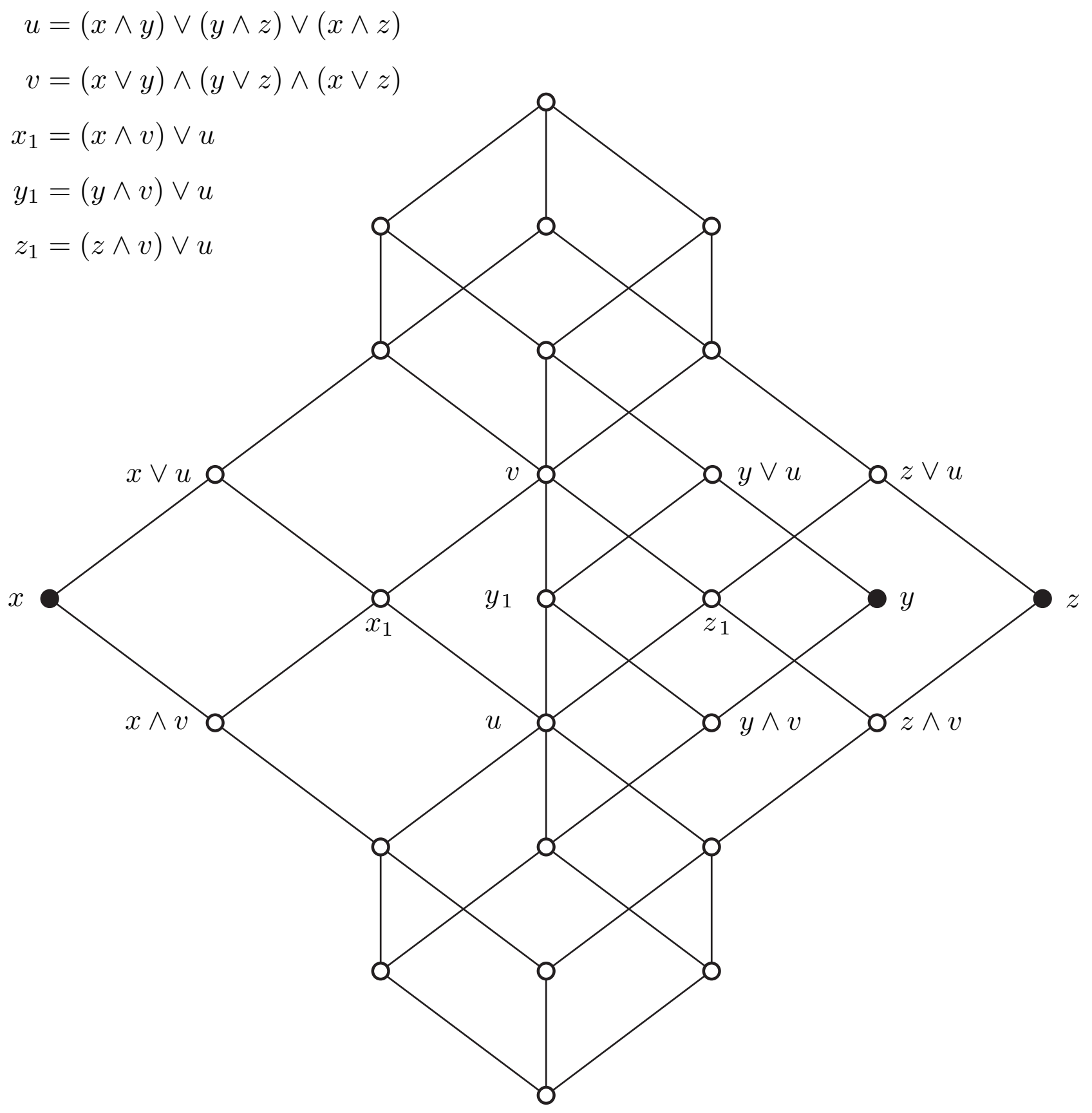

Figure 2.7: The free modular lattice on three generators, $\mathrm{Free}_{\mathbf{M}}(3)$

(i) $\mathrm{Free}_{\mathbf{K}}(H) \in \mathbf{K}$.

(ii) $\mathrm{Free}_{\mathbf{K}}(H)$ is generated by $H$.

(iii) Let $L \in \mathbf{K}$ and let $\psi\colon H \to L$ be a map; then there exists a (lattice) homomorphism $\varphi\colon \mathrm{Free}_{\mathbf{K}}(H) \to L$ extending $\psi$ (that is, satisfying $\varphi a = \psi a$ for all $a \in H$).

## 2.4. Gluing and generalizations

### 2.4.1 Gluing

In Section 1.1.3, for an ordered set $P$ with the largest element $1_P$ and an ordered set $Q$ with the smallest element $0_Q$, we introduced the glued sum $P \dotplus Q$. This applies to any two lattices $K$ with a unit and $L$ with a zero. A natural generalization of this construction is gluing.

Let $K$ and $L$ be lattices, let $F$ be a filter of $K$, and let $I$ be an ideal of $L$. If $F$ is isomorphic to $I$ (with $\varphi$ the isomorphism), then we can form the lattice $G$, the *gluing* of $K$ and $L$ over $F$ and $I$ (with respect to $\varphi$), defined as follows.

We form the disjoint union $K \cup L$, and identify $a \in F$ with $\varphi a \in I$, for all $a \in F$, to obtain the set $G$. We order $G$ as follows (see Figure 2.8):

$$a \leq b \quad \text{iff} \quad \begin{cases} a \leq_K b & \text{if } a, b \in K; \\ a \leq_L b & \text{if } a, b \in L; \\ a \leq_K x \text{ and } \varphi x \leq_L b & \text{if } a \in K \text{ and } b \in L \\ & \qquad \text{for some } x \in F. \end{cases}$$

**Lemma 2.7.** *$G$ is an ordered set, in fact, $G$ is a lattice. The join in $G$ is described by*

$$a \vee_G b = \begin{cases} a \vee_K b & \textit{if } a, b \in K; \\ a \vee_L b & \textit{if } a, b \in L; \\ \varphi(a \vee_K x) \vee_L b & \textit{if } a \in K \textit{ and } b \in L \textit{ for any } b \geq x \in F, \end{cases}$$

*and dually for the meet. If $L$ has a zero, $0_L$, then the last clause for the join may be rephrased:*

$$a \vee_G b = \varphi(a \vee_K 0_L) \vee_L b \quad \textit{if } a \in K \textit{ and } b \in L.$$

*$G$ contains $K$ and $L$ as sublattices; in fact, $K$ is an ideal and $L$ is a filter of $G$.*

An example of gluing is shown in Figure 2.9. There are more sophisticated examples in this book, for instance, in Chapters 12, 17, and Section 2.4 (see also Section 20.3).

**Lemma 2.8.** *Let $K$, $L$, $F$, $I$, and $G$ be given as above. Let $A$ be a lattice containing $K$ and $L$ as sublattices so that $K \cap L = I = F$. Then $K \cup L$ is a sublattice of $A$ and it is isomorphic to $G$.*

Now if $\boldsymbol{\alpha}_K$ is a binary relation on $K$ and $\boldsymbol{\alpha}_L$ is a binary relation on $L$, we define the *reflexive product* $\boldsymbol{\alpha}_K \overset{\mathrm{r}}{\circ} \boldsymbol{\alpha}_L$ as $\boldsymbol{\alpha}_K \cup \boldsymbol{\alpha}_L \cup (\boldsymbol{\alpha}_K \circ \boldsymbol{\alpha}_L) \cup (\boldsymbol{\alpha}_L \circ \boldsymbol{\alpha}_K)$.

We can easily describe the congruences of $G$.

**Lemma 2.9.** *A congruence $\boldsymbol{\alpha}$ of $G$ can be uniquely written in the form*

$$\boldsymbol{\alpha} = \boldsymbol{\alpha}_K \overset{r}{\circ} \boldsymbol{\alpha}_L,$$

*where* $\boldsymbol{\alpha}_K$ *is a congruence of* $K$ *and* $\boldsymbol{\alpha}_L$ *is a congruence of* $L$ *satisfying the condition that* $\boldsymbol{\alpha}_K$ *restricted to* $F$ *equals* $\boldsymbol{\alpha}_L$ *restricted to* $I$ *(under the identification of elements by* $\varphi$*).*

*Conversely, if* $\boldsymbol{\alpha}_K$ *is a congruence of* $K$ *and* $\boldsymbol{\alpha}_L$ *is a congruence of* $L$ *satisfying the condition that* $\boldsymbol{\alpha}_K$ *restricted to* $F$ *equals* $\boldsymbol{\alpha}_L$ *restricted to* $I$*, then* $\boldsymbol{\alpha} = \boldsymbol{\alpha}_K \overset{r}{\circ} \boldsymbol{\alpha}_L$ *is a congruence of* $G$*.*

Let $A$, $B$, and $C$ be lattices, $F_A$ a filter of $A$, $I_B$ an ideal of $B$, $F_B$ a filter of $B$, and $I_C$ an ideal of $C$. Let us assume that the lattices $F_A$, $I_B$, $F_B$, and $I_C$ are isomorphic.

We now define what it means to obtain $L$ by *double gluing* $A$, $B$, and $C$. Let $K$ be the gluing of $A$ and $B$ over $F_A$ and $I_B$ with the filter $F_B$ regarded as a filter $F_K$ of $K$ (see Figure 2.11). Then we glue $K$ and $C$ over $F_K$ and $I_C$ to obtain $L$, the double gluing of $A$ and $B$, and $C$.

Let $A$ and $B$ be lattices, $F_A$ a filter of $A$, $I_B$ an ideal of $B$, satisfying that $F_A$ and $I_B$ are isomorphic as lattices. For $k \geq 1$, we now define the *$k$-times gluing of $B$ to $A$*. For $k = 1$, let $L_1$ be the gluing of $A$ and $B$ over $F_A$ and $I_B$ with the filter $F_B$ regarded as a filter $F_{L_1}$ of $L_1$. Now if $L_{k-1}$ with the filter $F_{L_{k-1}}$ is the $(k-1)$-times gluing of $B$ to $A$, then we glue $L_{k-1}$ and $B$ over $F_{L_{k-1}}$ and $I_B$, to obtain $L$, the $k$-times gluing of $B$ to $A$ with the filter $F_B$ regarded as a filter $F_L$ of $L$.

Two references on gluings: my joint papers [90] and [91] with E. Fried.

### 2.4.2 Generalizations

Gluing is one of the most useful constructions for lattices, because it retains so many important properties: distributivity, modularity, semimodularity, planarity... So it is not surprising that there are so many generalizations. There are some very easy ones, such as the *triple gluing* of the section by the same title on page 55. E. Fried, G. Grätzer, and E. T. Schimidt [92] introduces *multipasting*, a common generalization of *$S$-glued sum* (where $S$ is a lattice of finite length) of C. Herrmann [209], and *pasting*, see V. Slavík [248] and my book GLT-[99], Exercise 12 of Section V.4. G. Czédli and E. T. Schmidt [74] introduced the *patchwork systems*, another generalization of gluing.

## 2.5. Modular and distributive lattices

### 2.5.1 The characterization theorems

The two typical examples of nondistributive lattices are $\mathsf{N}_5$ and $\mathsf{M}_3$, whose diagrams are given (again) in Figure 2.10. The following characterization theorem follows immediately by inspecting the diagrams of the free lattices in Figures 2.4 and 2.7.

**Theorem 2.10.**

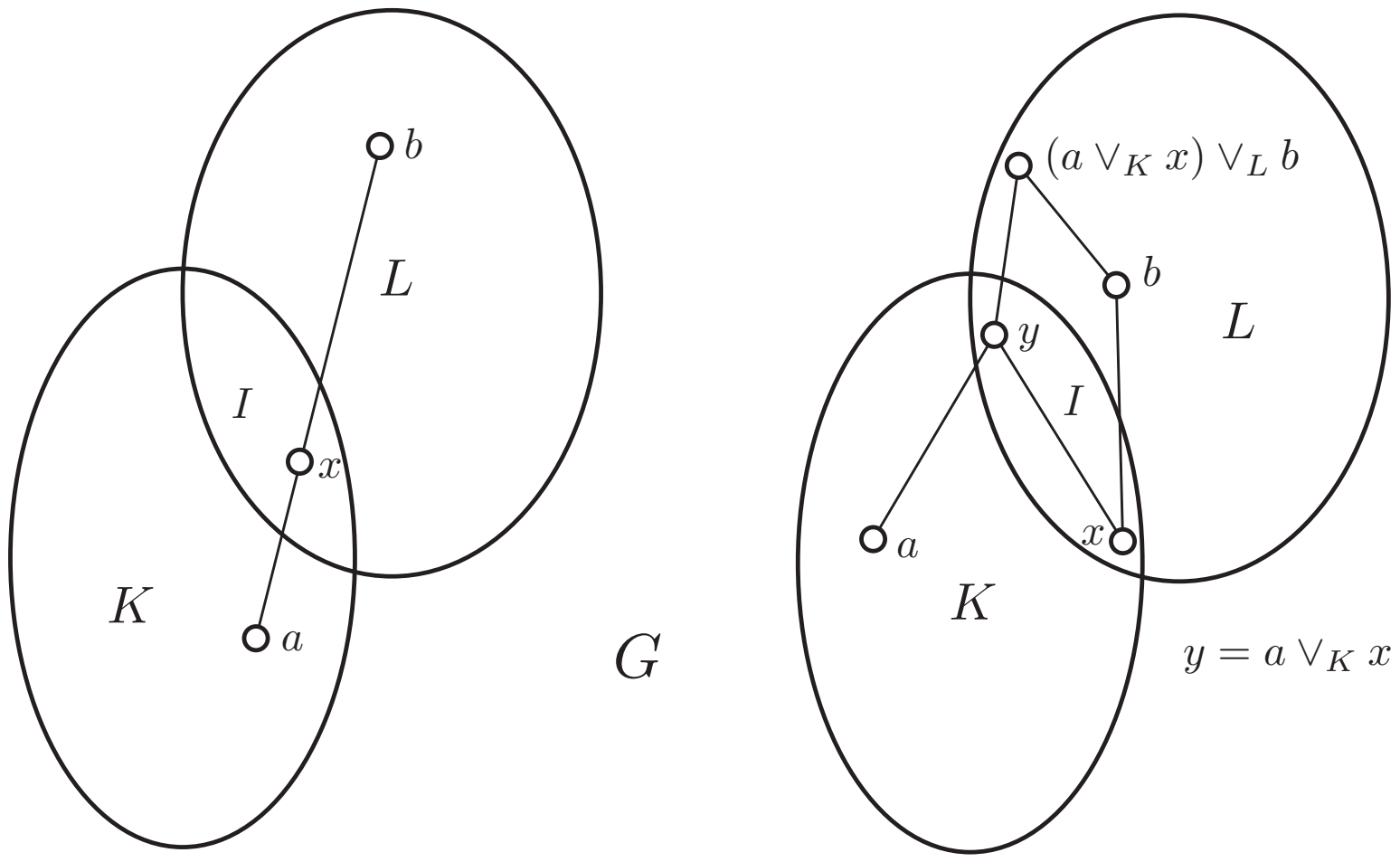

Figure 2.8: Defining gluing

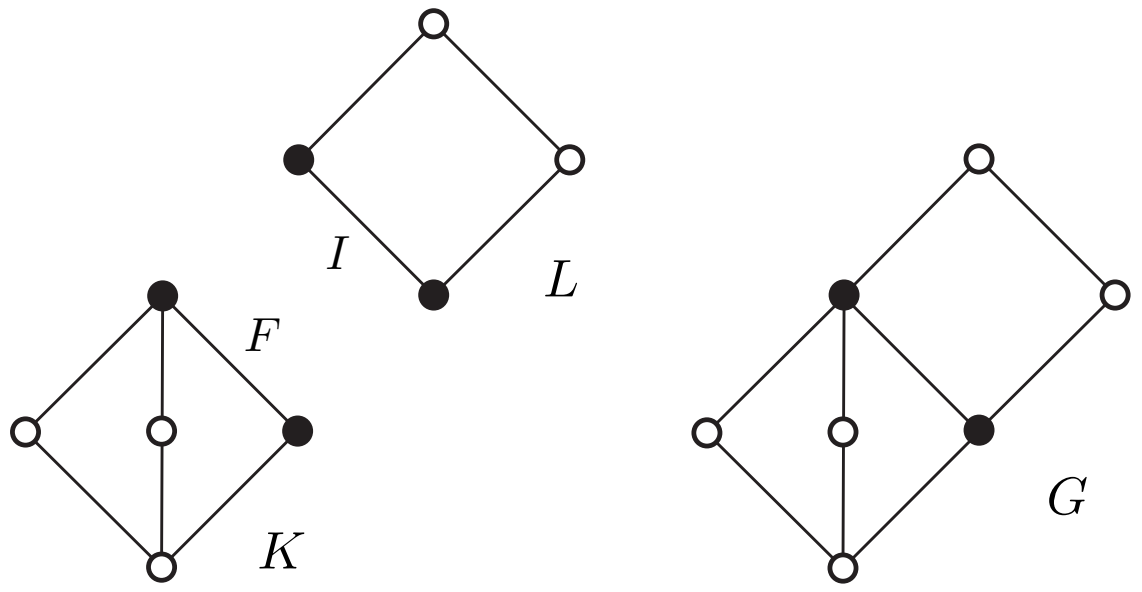

Figure 2.9: An easy gluing example

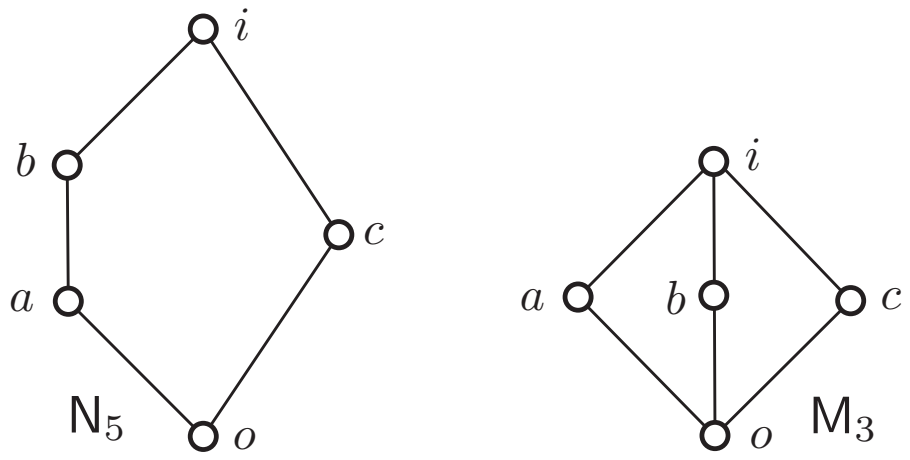

Figure 2.10: The two characteristic nondistributive lattices

(i) *A lattice $L$ is modular iff it does not contain* $\mathsf{N}_5$ *as a sublattice.*

(ii) *A modular lattice $L$ is distributive iff it does not contain* $\mathsf{M}_3$ *as a sublattice.*

(iii) *A lattice $L$ is distributive iff $L$ contains neither* $\mathsf{N}_5$ *nor* $\mathsf{M}_3$ *as a sublattice.*

**Theorem 2.11.** *Let $L$ be a modular lattice, let $a \in L$, and let $U$ and $V$ be sublattices with the property $u \wedge v = a$ for all $u \in U$ and $v \in V$. Then* $\operatorname{sub}(U \cup V)$ *is isomorphic to $U \times V$ under the isomorphism ($u \in U$ and $v \in V$)*

$$u \vee v \mapsto (u, v).$$

*Conversely, a lattice $L$ satisfying this property is modular.*

**Corollary 2.12.** *Let $L$ be a modular lattice and let $a, b \in L$. Then*

$$\operatorname{sub}([a \wedge b, a] \cup [a \wedge b, b]),$$

*that is, the sublattice of $L$ generated by $[a \wedge b, a] \cup [a \wedge b, b]$, is isomorphic to the direct product*

$$[a \wedge b, a] \times [a \wedge b, b].$$

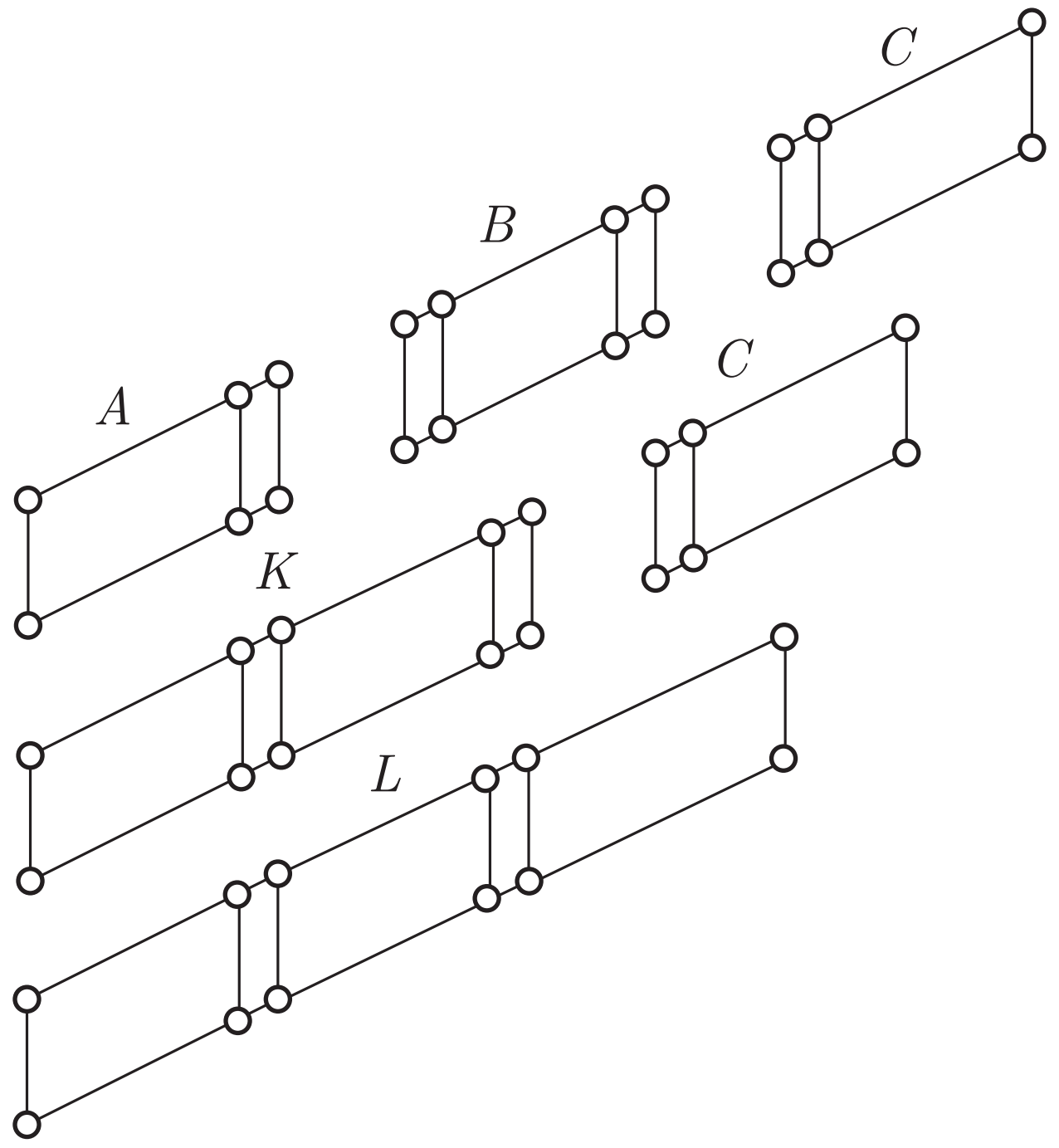

Figure 2.11: Double gluing

In the distributive case, the sublattice generated by $[a \wedge b, a] \cup [a \wedge b, b]$ is the interval $[a \wedge b, a \vee b]$; this does not hold for modular lattices, as exemplified by $\mathsf{M}_3$.

Let $G$ be the gluing of the lattices $K$ and $L$ over $F$ and $I$, as in Section 2.4.

**Lemma 2.13.** *If $K$ and $L$ are modular, so is the gluing $G$ of $K$ and $L$. If $K$ and $L$ are distributive, so is $G$.*

The distributive identity easily implies that every $n$-ary term equals one we get by joining meets of variables (the *disjunctive normal form*). So we get:

**Lemma 2.14.** *A finitely generated distributive lattice is finite.*

### 2.5.2 Finite distributive lattices

For a finite distributive lattice $L$, set

$$\operatorname{spec}(a) = \{\, x \in \mathrm{J}(L) \mid x \leq a \,\} = \operatorname{id}(a) \cap \mathrm{J}(L) = \downarrow a \cap \mathrm{J}(L),$$

the *spectrum* of $a$.

The structure of finite distributive lattices is described by the following result.

**Theorem 2.15.** *Let $L$ be a finite distributive lattice. Then the map*

$$\varphi\colon a \mapsto \operatorname{spec}(a)$$

*is an isomorphism between $L$ and* $\mathrm{Dn}(\mathrm{J}(L))$.

**Corollary 2.16.** *The correspondence $L \mapsto \mathrm{J}(L)$ makes the class of all finite distributive lattices with more than one element correspond to the class of all finite orders; isomorphic lattices correspond to isomorphic orders, and vice versa.*

*Proof.* This is obvious from $\mathrm{J}(\mathrm{Dn}\, P) \cong P$ and $\mathrm{Dn}(\mathrm{J}(L)) \cong L$. □

The RT relates bounded homomorphisms of distributive lattices and isotone maps of ordered sets.

**Theorem 2.17.** *Let $D$ and $E$ be finite distributive lattices and set $P = \mathrm{J}(D)$ and $Q = \mathrm{J}(E)$. Then*

(i) *With every bounded homomorphism $\varphi\colon D \to E$, we can associate an isotone map $\mathrm{J}(\varphi)\colon Q \to P$ defined by*

$$\mathrm{J}(\varphi)x = \bigwedge(\, e \in D \mid x \leq \varphi e \,)$$

*for $x \in Q$.*

(ii) *With every isotone map* $\psi\colon Q \to P$, *we can associate a bounded homomorphism* $\mathrm{Dn}(\psi)\colon D \to E$ *defined by*

$$\mathrm{Dn}(\psi)(e) = \bigvee \psi^{-1}(\mathrm{spec}(e))$$

*for* $e \in D$.

(iii) *The constructions of* (i) *and* (ii) *are inverses to one another, and so yield together a bijection between bounded homomorphisms* $\mathrm{J}(\varphi)\colon Q \to P$ *and isotone maps* $\mathrm{J}(\varphi)\colon D \to E$.

(iv) $\varphi$ *is one-to-one iff* $\mathrm{J}(\varphi)$ *is onto.*

(v) $\varphi$ *is onto iff* $\mathrm{J}(\varphi)$ *is an order-embedding.*

This result is illustrated by Figures 2.12–2.13.

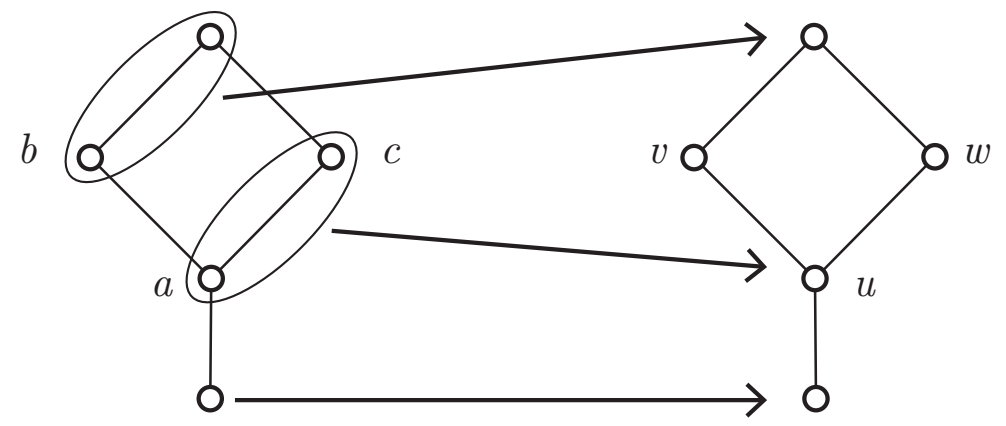

Figure 2.12: An example of the distributive lattices $D$ and $E$ and a bounded homomorphism $\varphi\colon D \to E$

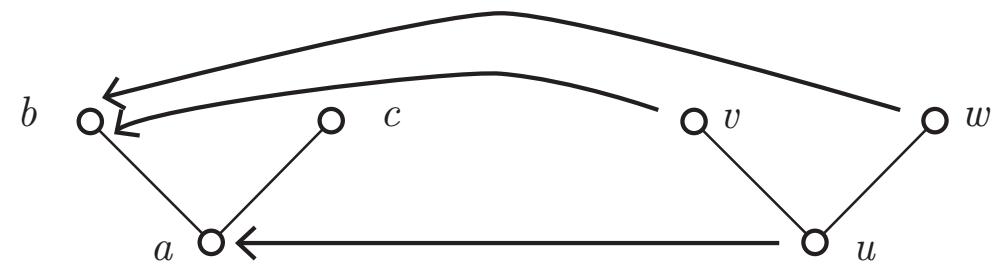

Figure 2.13: An example of the ordered sets $P$ and $Q$ and an isotone map $\psi\colon Q \to P$

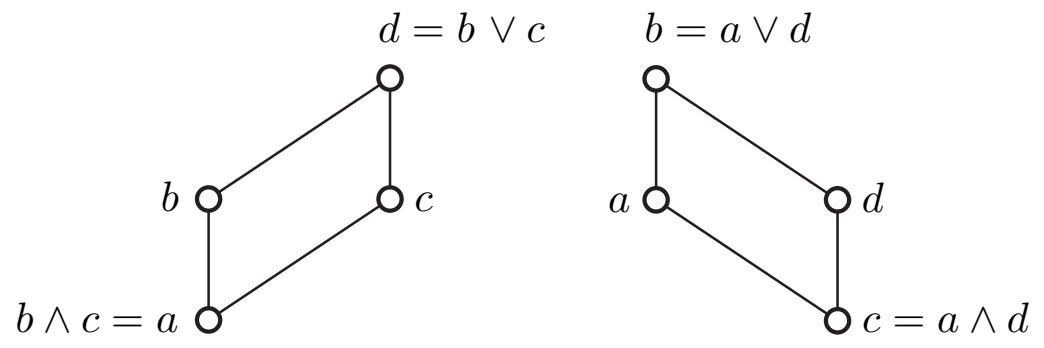

Figure 2.14: $[a,b] \overset{\mathrm{up}}{\sim} [c,d]$ and $[a,b] \overset{\mathrm{dn}}{\sim} [c,d]$

The two main formulas in this result are easy to visualize. To verbalize the formula in (i), for $x \in Q$, take the set $E$ of all elements of $D$ mapped by $\varphi$ onto $x$ or above. Since $\varphi$ is a homomorphism, it follows that the meet of all the elements (let us denote it by $x^\dagger$) is still mapped by $\varphi$ onto $x$ or above. The map $\mathrm{J}(\varphi)$ maps $x$ into $x^\dagger$.

With every isotone map $\psi\colon Q \to P$, we can associate a bounded homomorphism $\mathrm{Dn}(\psi)\colon D \to E$ defined by

$$\mathrm{Dn}(\psi)(e) = \bigvee \psi^{-1}(\mathrm{spec}(e))$$

for $e \in D$.

Now for the formula in (ii), take an element $e \in D$. We form the set of join-irreducible elements of $E$ mapped by $\psi$ to an element $\leq e$. The join of this set, denoted by $e^\ddagger$, is an element of $E$. The map $\mathrm{Dn}(\psi)$ maps $e$ into $e^\ddagger$.

The constructions of (i) and (ii) are inverses to one another, and so yield together a bijection between bounded homomorphisms $\varphi\colon D \to E$ and isotone maps $\psi\colon Q \to P$.

These results express the duality between finite distributive lattices, with bounded homomorphisms, and finite ordered sets, with isotone maps ("duality" means "equivalence between a category and the opposite — aka dual — of another one").

### 2.5.3 Finite modular lattices

Take a look at the two positions of the pair of intervals $[a, b]$ and $[c, d]$ in Figure 2.14. In either case, we will write $[a, b] \sim [c, d]$, and say that $[a, b]$ is *perspective* to $[c, d]$. If we want to show whether the perspectivity is "up" or "down," we will write $[a, b] \overset{\mathrm{up}}{\sim} [c, d]$ in the first case and $[a, b] \overset{\mathrm{dn}}{\sim} [c, d]$ in the second.

If for some natural number $n$ and intervals $[e_i, f_i]$, for $0 \leq i \leq n$,

$$[a, b] = [e_0, f_0] \sim [e_1, f_1] \sim \cdots \sim [e_n, f_n] = [c, d],$$

then we say that $[a, b]$ is *projective* to $[c, d]$ and write $[a, b] \approx [c, d]$.

One of the most important properties of a modular lattice is stated in the following result.

**Theorem 2.18** (Isomorphism Theorem for Modular Lattices)**.** *Let $L$ be a modular lattice and let $[a, b] \overset{up}{\sim} [c, d]$ in $L$. Then*

$$\varphi_c\colon x \mapsto x \vee c, \quad x \in [a, b],$$

*is an isomorphism of $[a, b]$ and $[c, d]$. The inverse isomorphism is*

$$\boldsymbol{\beta}_b\colon y \mapsto y \wedge b, \quad y \in [c, d].$$

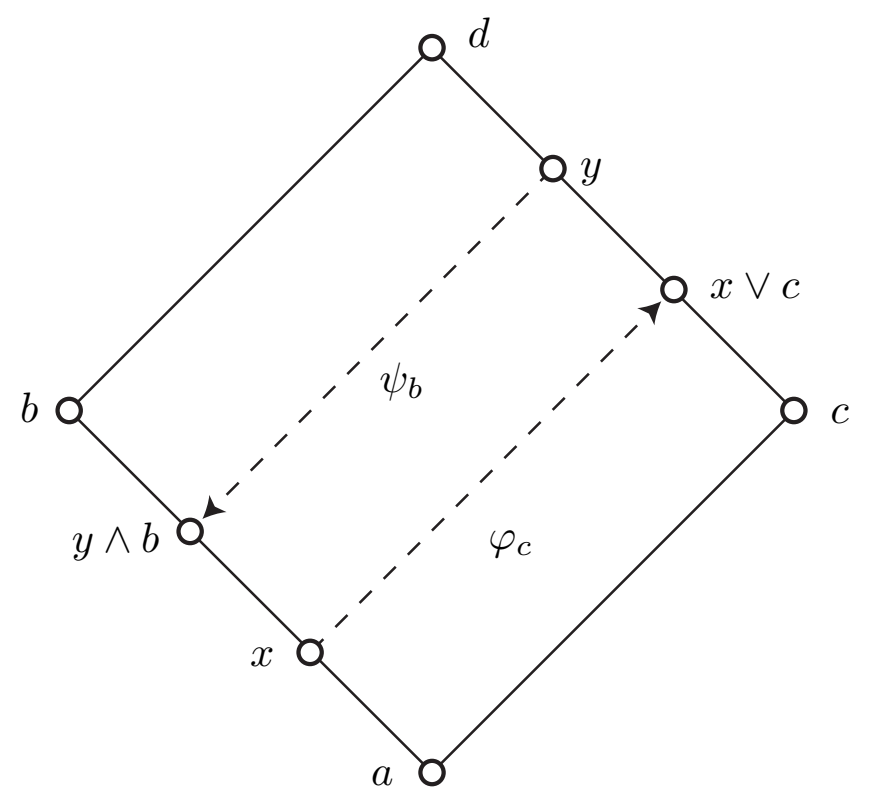

Figure 2.15: The isomorphisms $\varphi_c$ and $\boldsymbol{\beta}_b$

(See Figure 2.15.)

**Corollary 2.19.** *In a modular lattice, projective intervals are isomorphic.*

**Corollary 2.20.** *In a modular lattice if a prime interval* $\mathfrak{p}$ *is projective to an interval* $\mathfrak{q}$*, then* $\mathfrak{q}$ *is also prime.*

Let us call the finite lattice $L$ *semimodular* or *upper semimodular* if for $a, b, c \in L$, the covering $a \prec b$ implies that $a \vee c \prec b \vee c$ or $a \vee c = b \vee c$. The dual of an upper semimodular lattice is a *lower semimodular* lattice.

**Lemma 2.21.** *A modular lattice is both upper and lower semimodular. For a* finite *lattice, the converse also holds: a finite upper and lower semimodular lattice is modular, and conversely.*

The lattice $\mathsf{S}_8$ (see the Picture Gallery) is an example of a semimodular lattice that is not modular. Section 10.2 provides an interesting use of this lattice.

The following is even more trivial than Lemma 2.13.

**Lemma 2.22.** *If $K$ and $L$ are finite semimodular lattices, so is the gluing $G$ of $K$ and $L$.*

A large class of semimodular lattices is provided by

**Lemma 2.23.** *Let $A$ be a nonempty set. Then* $\operatorname{Part} A$ *is semimodular; it is not modular unless* $|A| \leq 3$.

**Chapter**

**3**

# *Congruences*

## 3.1. Congruence spreading

Let $a, b, c, d$ be elements of a lattice $L$. If

$$a \equiv b \pmod{\boldsymbol{\alpha}} \quad \text{implies that} \quad c \equiv d \pmod{\boldsymbol{\alpha}},$$

for any congruence relation $\boldsymbol{\alpha}$ of $L$, then we say that $a \equiv b$ *congruence-forces* $c \equiv d$.

In Section 1.3.3 we saw that $a \equiv b \pmod{\boldsymbol{\alpha}}$ iff $a \wedge b \equiv a \vee b \pmod{\boldsymbol{\alpha}}$; therefore, to investigate congruence-forcing, it is enough to deal with comparable pairs, $a \leq b$ and $c \leq d$. Instead of comparable pairs, we will deal with intervals $[a, b]$ and $[c, d]$.

Projectivity (see Section 2.5.3) is sufficient for the study of congruence-forcing (or congruence-spreading) in some classes of lattices (for instance, in the class of modular lattices). In general, however, we have to introduce somewhat more general concepts and notation.

As illustrated in Figure 3.1, we say that $[a, b]$ is *up congruence-perspective* onto $[c, d]$ and write $[a, b] \overset{\text{up}}{\twoheadrightarrow} [c, d]$ iff there is an $a_1 \in [a, b]$ with $[a_1, b] \overset{\text{up}}{\sim} [c, d]$. Similarly, $[a, b]$ is *down congruence-perspective* onto $[c, d]$ and we shall write $[a, b] \overset{\text{dn}}{\twoheadrightarrow} [c, d]$ iff there is a $b_1 \in [a, b]$ with $[a, b_1] \overset{\text{dn}}{\sim} [c, d]$. If $[a, b] \overset{\text{up}}{\twoheadrightarrow} [c, d]$ or $[a, b] \overset{\text{dn}}{\twoheadrightarrow} [c, d]$, then $[a, b]$ is *congruence-perspective onto* $[c, d]$ and we write $[a, b] \twoheadrightarrow [c, d]$. If for some natural number $n$ and and intervals $[e_i, f_i]$ for $0 \leq i \leq n$,

$$[a, b] = [e_0, f_0] \twoheadrightarrow [e_1, f_1] \twoheadrightarrow \cdots \twoheadrightarrow [e_n, f_n] = [c, d],$$

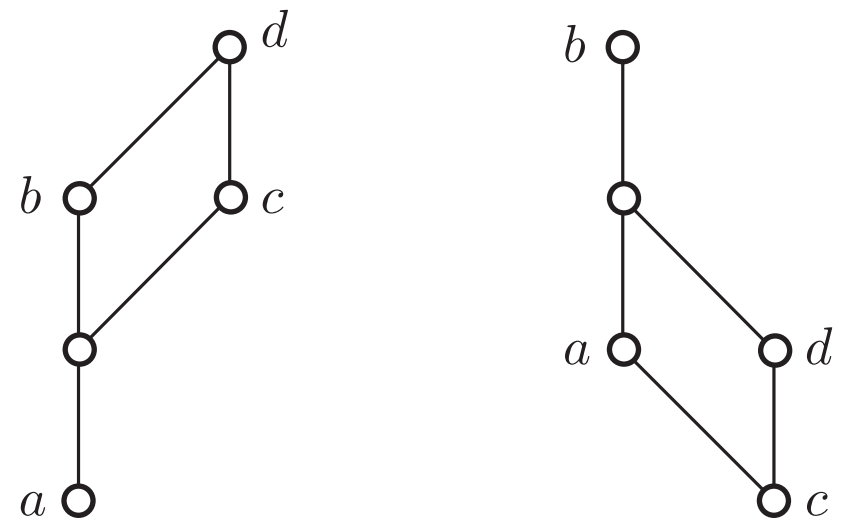

Figure 3.1: $[a,b] \overset{\text{up}}{\twoheadrightarrow} [c,d]$ and $[a,b] \overset{\text{dn}}{\twoheadrightarrow} [c,d]$

then we call $[a,b]$ *congruence-projective onto* $[c,d]$, and we write $[a,b] \Rightarrow [c,d]$.

If $[a,b] \Rightarrow [c,d]$ and $[c,d] \Rightarrow [a,b]$, then we write $[a,b] \Leftrightarrow [c,d]$.

Also if (see Section 2.5.3) $[a,b] \overset{\text{up}}{\sim} [c,d]$, then $[a,b] \overset{\text{up}}{\twoheadrightarrow} [c,d]$ and $[c,d] \overset{\text{dn}}{\twoheadrightarrow} [a,b]$; if $[a,b] \overset{\text{dn}}{\sim} [c,d]$, then $[a,b] \overset{\text{dn}}{\twoheadrightarrow} [c,d]$ and $[c,d] \overset{\text{up}}{\twoheadrightarrow} [a,b]$; if $[a,b] \sim [c,d]$, then $[a,b] \twoheadrightarrow [c,d]$ and $[c,d] \twoheadrightarrow [a,b]$; if $[a,b] \approx [c,d]$, then $[a,b] \Rightarrow [c,d]$ and $[c,d] \Rightarrow [a,b]$. But while $\sim$, $\approx$, and $\Leftrightarrow$ are symmetric, the relations $\twoheadrightarrow$ and $\Rightarrow$ are not. In particular, if $a \leq c \leq d \leq b$, then $[a,b] \Rightarrow [c,d]$.

If $[a,b] \Rightarrow [c,d]$, then there is a unary term function $p$ with $p(a) = c$ and $p(c) = d$. It is easy to see what special kinds of unary term functions are utilized in the relation $\Rightarrow$.

Intuitively, "$a \equiv b$ congruence-forces $c \equiv d$" iff $[c,d]$ is put together from pieces $[c', d']$, each of which satisfies $[a,b] \Rightarrow [c', d']$. To state this more precisely, we describe the smallest congruence relation under which $a \equiv b$, denoted by $\mathrm{con}(a,b)$, as introduced in Section 1.3.3 (see R. P. Dilworth [79]). We use the Technical Lemma to prove this result.

**Theorem 3.1.** *Let $L$ be a lattice, $a, b, c, d \in L$, $a \leq b$, $c \leq d$. Then $c \equiv d$ $(\mathrm{mod}\ \mathrm{con}(a,b))$ iff, for some sequence*

$$c = e_0 \leq e_1 \leq \cdots \leq e_m = d,$$

*we have*

$$[a,b] \Rightarrow [e_j, e_{j+1}] \qquad \textit{for} \quad j = 0, \ldots, m-1.$$

This result can be usefully augmented by the following lemma (for an application, see Lemma 25.2).

**Lemma 3.2.** *Let $L$ be a lattice, $a, b, c, d \in L$ with $a \leq b$ and $c \leq d$. Then $[a,b]$ is congruence-projective to $[c,d]$ iff the following condition is satisfied:*

*There is an integer $m$ and there are elements $e_0, \ldots, e_{m-1} \in L$ such that*

$$p_m(a, e_0, \ldots, e_{m-1}) = c,$$
$$p_m(b, e_0, \ldots, e_{m-1}) = d,$$

*where the term* $p_m$ *is defined by*

$$p_m(x, y_0, \ldots, y_{m-1}) = \cdots (((x \vee y_0) \wedge y_1) \vee y_2) \wedge \cdots .$$

Let $L$ be a lattice and $H \subseteq L^2$. To compute $\operatorname{con}(H)$, the smallest congruence relation $\boldsymbol{\alpha}$ under which $a \equiv b \pmod{\boldsymbol{\alpha}}$ for all $(a, b) \in H$, we use the formula

$$\operatorname{con}(H) = \bigvee(\, \operatorname{con}(a, b) \mid (a, b) \in H\,).$$

We also need a formula for joins.

**Lemma 3.3.** *Let* $L$ *be a lattice and let* $\boldsymbol{\alpha}_i$*,* $i \in I$*, be congruence relations of* $L$*. Then* $a \equiv b \pmod{\bigvee(\, \boldsymbol{\alpha}_i \mid i \in I\,)}$ *iff there is a sequence*

$$z_0 = a \wedge b \leq z_1 \leq \cdots \leq z_n = a \vee b$$

*such that for each* $j$ *with* $0 \leq j < n$*, there is an* $i_j \in I$ *satisfying* $z_j \equiv z_{j+1} \pmod{\boldsymbol{\alpha}_{i_j}}$*.*

This is an easy but profoundly important result. For instance, the result of N. Funayama and T. Nakayama [95]—which provides the foundation for this book—immediately follows. (Another typical application is Lemma 3.12.)

**Theorem 3.4.** *The lattice* $\operatorname{Con} L$ *is distributive for any lattice* $L$*.*

By combining Theorem 3.1 and Lemma 3.3, we get:

**Corollary 3.5.** *Let* $L$ *be a lattice, let* $H \subseteq L^2$*, and let* $a, b \in L$ *with* $a \leq b$*. Then* $a \equiv b \pmod{\operatorname{con}(H)}$ *iff, for some integer* $n$*, there exists a sequence*

$$a = c_0 \leq c_1 \leq \cdots \leq c_n = b$$

*such that for each* $i$ *with* $0 \leq i < n$*, there exists a* $(d_i, e_i) \in H$ *satisfying*

$$[d_i \wedge e_i, d_i \vee e_i] \Rightarrow [c_i, c_{i+1}].$$

There is another congruence structure we can associate with a lattice $L$, the ordered set of principal congruences, $\operatorname{Princ} L$, a subset of $\operatorname{Con} L$. We dedicate Part VI to the study of this structure.

**Lemma 3.6.** *For a lattice* $L$*, the ordered set* $\operatorname{Princ} L$ *is a directed ordered set with zero. If* $L$ *is bounded, so is* $\operatorname{Princ} L$*.*

Indeed, for $a, b, c, d \in L$, an upper bound of $\operatorname{con}(a, b)$ and $\operatorname{con}(c, d)$ is $\operatorname{con}(a \wedge b \wedge c \wedge d, a \vee b \vee c \vee d)$.

## 3.2. Finite lattices and prime intervals

Let $\mathrm{Prime}(L)$ denote the set of prime intervals of a finite lattice $L$. Let $\mathfrak{p} = [a, b]$ be a prime interval in $L$ (that is, $a \prec b$ in $L$) and let $\boldsymbol{\alpha}$ be a congruence relation of $L$.

In a finite lattice $L$, the formula

$$\boldsymbol{\alpha} = \bigvee(\,\mathrm{con}(\mathfrak{p}) \mid 0_{\mathfrak{p}} \equiv 1_{\mathfrak{p}} \pmod{\boldsymbol{\alpha}}\,)$$

immediately yields that the congruences in $\mathrm{J}(\mathrm{Con}\, L)$ (which we shall denote by $\mathrm{J}(\mathrm{Con}\, L)$) are the congruences of the form $\mathrm{con}(\mathfrak{p})$ for some $\mathfrak{p} \in \mathrm{Prime}(L)$. Of course, a join-irreducible congruence $\boldsymbol{\alpha}$ can be expressed, as a rule, in many ways in the form $\mathrm{con}(\mathfrak{p})$.

Since a prime interval cannot contain a three-element chain, the following two lemmas easily follow from the Technical Lemma and from Corollary 3.5.

**Lemma 3.7.** *Let $L$ be a finite lattice and let $\mathfrak{p}$ and $\mathfrak{q}$ be prime intervals in $L$. Then* $\mathrm{con}(\mathfrak{p}) \geq \mathrm{con}(\mathfrak{q})$ *iff* $\mathfrak{p} \Rightarrow \mathfrak{q}$.

This condition is easy to visualize using Figure 3.2; the sequence of congruence-perspectivities, as a rule, has to go through nonprime intervals (intervals of arbitrary size) to get from $\mathfrak{p}$ to $\mathfrak{q}$.

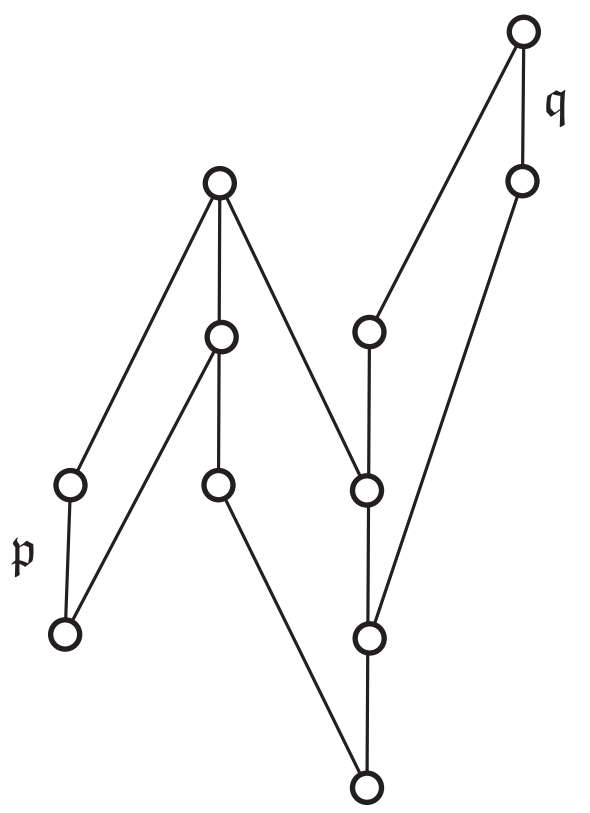

Figure 3.2: Congruence spreading from prime interval to prime interval

**Lemma 3.8.** *Let $L$ be a finite lattice, let $\mathfrak{p}$ be a prime interval of $L$, and let $[a, b]$ be an interval of $L$. If $\mathfrak{p}$ is collapsed under* $\mathrm{con}(a, b)$*, then there is a prime interval $\mathfrak{q}$ in $[a, b]$ satisfying* $\mathfrak{q} \Rightarrow \mathfrak{p}$.

In view of Theorem 2.15, we get the following.

**Theorem 3.9.** *Let $L$ be a finite lattice. The relation $\Rightarrow$ is a quasiordering on* $\mathrm{Prime}(L)$*. The equivalence classes under $\Leftrightarrow$ form an ordered set isomorphic to* $\mathrm{J}(\mathrm{Con}\, L)$*.*

If the finite lattice $L$ is atomistic, then the join-irreducible congruences are even simpler to find. Indeed, if $[a, b]$ is a prime interval, and $p$ is an atom with $p \leq b$ and $p \nleq a$, then $[a, b] \overset{\text{dn}}{\sim} [0, p]$, so $\mathrm{con}(a, b) = \mathrm{con}(0, p)$.

**Corollary 3.10.** *Let $L$ be a finite atomistic lattice. Then every join-irreducible congruence can be represented in the form* $\mathrm{con}(0, p)$*, where $p$ is an atom. The relation $p \Rightarrow q$, defined as $[0, p] \Rightarrow [0, q]$, introduces a quasiordering on the set of atoms of $L$. The equivalence classes under the quasiordering form an ordered set isomorphic to* $\mathrm{J}(\mathrm{Con}\, L)$*.*

We use Figure 3.3 to illustrate how we compute the congruence lattice of $\mathsf{N}_5$ using Theorem 3.9. $\mathsf{N}_5$ has five prime intervals: $[o, a]$, $[a, b]$, $[b, i]$, $[o, c]$, $[c, i]$. The equivalence classes are $\alpha = \{[a, b]\}$, $\beta = \{[o, c], [b, i]\}$, and $\gamma = \{[o, a], [c, i]\}$. The ordering $\alpha < \gamma$ holds because $[c, i] \overset{\text{dn}}{\twoheadrightarrow} [o, b] \overset{\text{up}}{\twoheadrightarrow} [a, b]$. Similarly, $\alpha < \beta$.

It is important to note that the computation of $\mathfrak{p} \Rightarrow \mathfrak{q}$ may involve non-prime intervals. For instance, let $\mathfrak{p} = [o, c]$ and $\mathfrak{q} = [a, b]$ in $\mathsf{N}_5$. Then $\mathfrak{p} \Rightarrow \mathfrak{q}$, because $\mathfrak{p} \overset{\text{up}}{\twoheadrightarrow} [a, i] \overset{\text{dn}}{\twoheadrightarrow} \mathfrak{q}$, but we cannot get $\mathfrak{p} \Rightarrow \mathfrak{q}$ involving only prime intervals.

As another example, we compute the congruence lattice of $\mathsf{S}_8$ (see Figure 3.4).

By Corollary 2.21, in a modular lattice if $\mathfrak{p}$ and $\mathfrak{q}$ are prime intervals, then $\mathfrak{p} \overset{\text{up}}{\twoheadrightarrow} \mathfrak{q}$ implies that $\mathfrak{p} \overset{\text{up}}{\sim} \mathfrak{q}$ and $\mathfrak{p} \overset{\text{dn}}{\twoheadrightarrow} \mathfrak{q}$ implies that $\mathfrak{p} \overset{\text{dn}}{\sim} \mathfrak{q}$; thus $\mathfrak{p} \Rightarrow \mathfrak{q}$ implies that $\mathfrak{p} \approx \mathfrak{q}$. Therefore, Theorem 3.9 tells us that $\mathrm{J}(\mathrm{Con}\, L)$ is an antichain in a finite modular lattice $L$, so $\mathrm{Con}\, L$ is Boolean.

**Corollary 3.11.** *The congruence lattice of a finite modular lattice is Boolean.*

By a *colored lattice* we will mean a finite lattice (some) of whose prime intervals are labeled so that if the prime intervals $\mathfrak{p}$ and $\mathfrak{q}$ are of the same color, then $\operatorname{con}(\mathfrak{p}) = \operatorname{con}(\mathfrak{q})$. These labels represent (a subset of) the equivalence classes of prime intervals, as stated in Theorem 3.9. In Figure 3.4, every prime interval of $\mathsf{S}_8$ is labeled; in Figure 9.2, only some are labeled. The coloring helps in the intuitive understanding of some constructions.

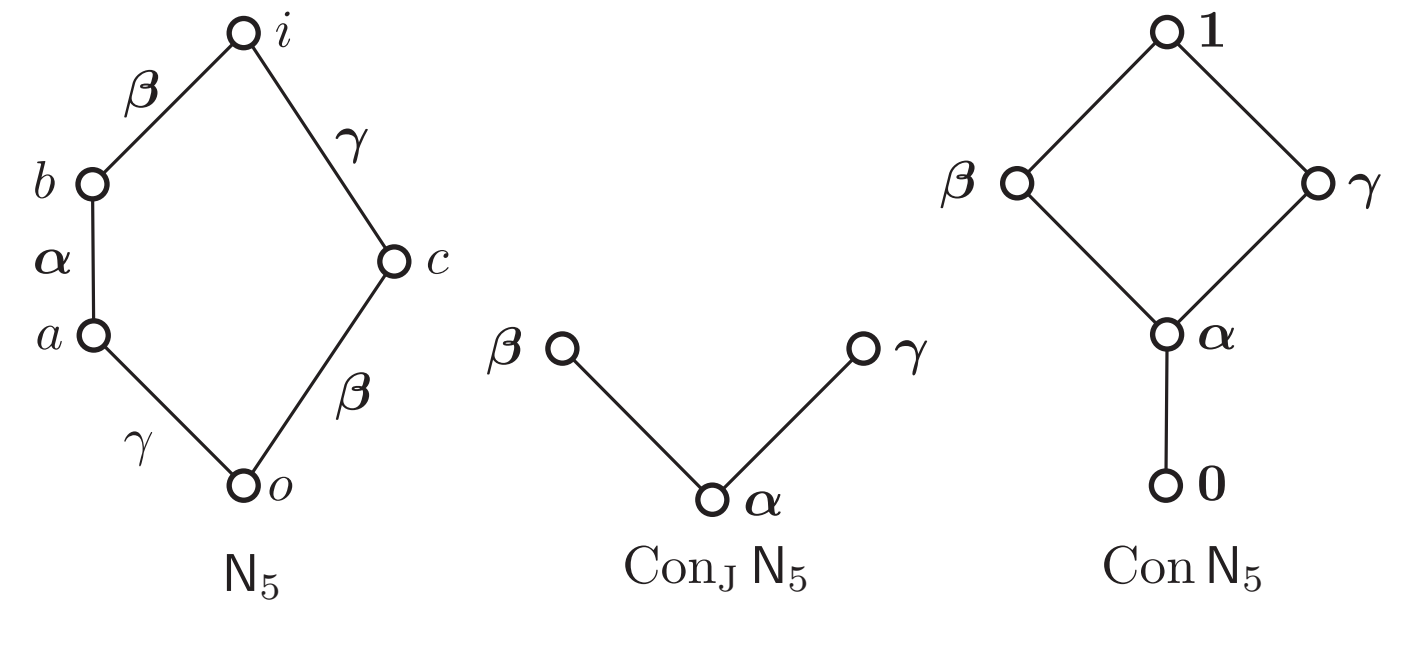

Figure 3.3: Computing the congruence lattice of $\mathsf{N}_5$

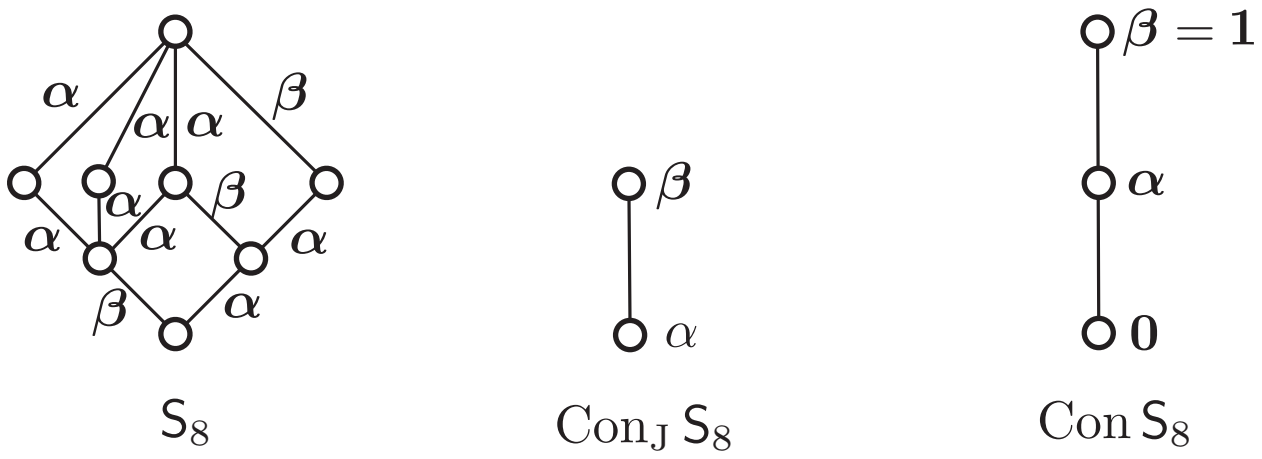

Figure 3.4: Computing the congruence lattice of $\mathsf{S}_8$

## 3.3. Congruence-preserving extensions and variants

Let $L$ be a lattice and $K \leq L$. How do the congruences of $L$ relate to the congruences of $K$?

Every congruence $\boldsymbol{\alpha}$ *restricts* to $K$: the relation $\boldsymbol{\alpha} \cap K^2 = \boldsymbol{\alpha}\rceil K$ on $K$ is a congruence of $K$. So we get the *restriction map*:

$$\operatorname{re}\colon \operatorname{Con} L \to \operatorname{Con} K,$$

which maps a congruence $\boldsymbol{\alpha}$ of $L$ to $\boldsymbol{\alpha}\rceil K$.

**Lemma 3.12.** *Let $K \leq L$ be lattices. Then* $\operatorname{re}\colon \operatorname{Con} L \to \operatorname{Con} K$ *is a* $\{\wedge, 0, 1\}$-*homomorphism.*

For instance, if $K = \{o, a, i\}$ and $L = \mathsf{M}_3$ (see Figure 2.10), then $\operatorname{Con} K$ is isomorphic to $\mathsf{B}_2$, but only $\mathbf{0}$ and $\mathbf{1}$ are restrictions of congruences in $L$. As another example, take the lattice $L$ of Figure 3.5 and its sublattice $K$, the black-filled elements; in this case, $\operatorname{Con} L \cong \operatorname{Con} K \cong \mathsf{B}_2$, but again only $\mathbf{0}$ and $\mathbf{1}$ are restrictions. There is no natural relationship between the congruences of $K$ and $L$.

If $K$ is an ideal in $L$ (or any convex sublattice), we can say a lot more.

**Lemma 3.13.** *Let $K \leq L$ be lattices. If $K$ is an ideal of $L$, then* $\operatorname{re}\colon \operatorname{Con} L \to \operatorname{Con} K$ *is a bounded homomorphism.*

Let $K \leq L$ be lattices, and let $\boldsymbol{\alpha}$ be a congruence of $K$. The congruence $\operatorname{con}_L(\boldsymbol{\alpha})$ (the congruence $\operatorname{con}(\boldsymbol{\alpha})$ formed in $L$) is the smallest congruence $\boldsymbol{\gamma}$ of $L$ such that $\boldsymbol{\alpha} \leq \boldsymbol{\gamma}\rceil K$. Unfortunately, $\operatorname{con}_L(\boldsymbol{\alpha})\rceil K$ may be different from $\boldsymbol{\alpha}$, as in the example of Figure 3.5. We say that the congruence $\boldsymbol{\alpha}$ of $K$ *extends to $L$,* iff $\boldsymbol{\alpha}$ is the restriction of $\operatorname{con}_L(\boldsymbol{\alpha})$. Figure 3.6 illustrates this in part.

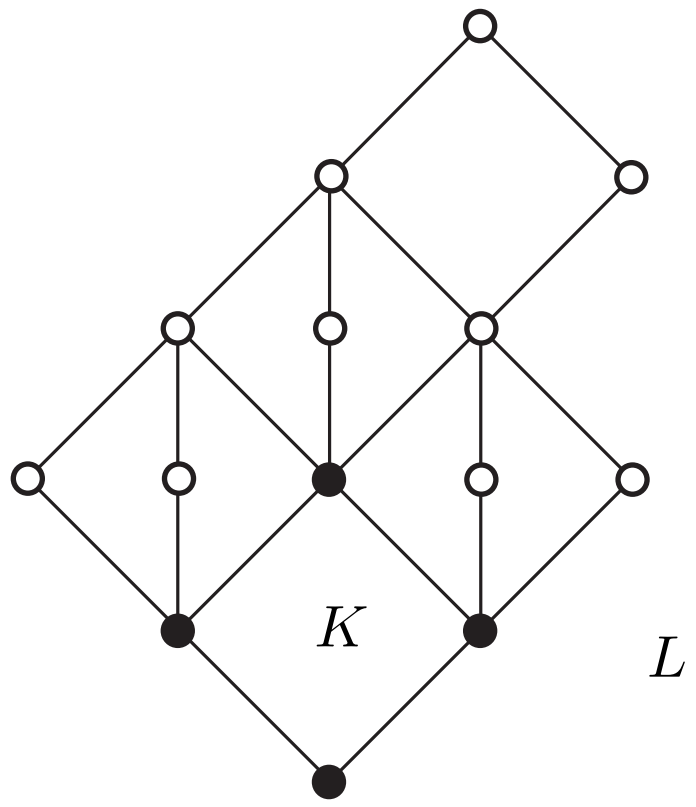

Figure 3.5: A sublattice for the discussion of the map re

If a congruence $\boldsymbol{\alpha}$ extends, then the congruence classes of $\boldsymbol{\alpha}$ in $K$ extend to congruence classes in $L$, but there may be congruence classes in $L$ that are not such extensions.

The *extension map*:

$$\mathrm{ext}\colon \operatorname{Con} K \to \operatorname{Con} L$$

maps a congruence $\boldsymbol{\alpha}$ of $K$ to the congruence $\mathrm{con}_L(\boldsymbol{\alpha})$ of $L$. The map ext is a $\{\vee, 0\}$-homomorphism of $\operatorname{Con} K$ into $\operatorname{Con} L$. In addition, ext preserves the zero, that is, ext is $\{0\}$*-separating.* To summarize,

**Lemma 3.14.** *Let $K \leq L$ be lattices. Then* $\mathrm{ext}\colon \operatorname{Con} K \to \operatorname{Con} L$ *is a $\{0\}$-separating join-homomorphism.*

The extension $L$ of $K$ is a *congruence-reflecting extension* (and the sublattice $K$ of $L$ is a *congruence-reflecting sublattice*) if every congruence of $K$ extends to $L$.

Utilizing the results of Sections 3.1 and 3.2, we can find many equivalent formulations of the congruence-reflecting extension property for finite lattices.

**Lemma 3.15.** *Let $L$ be a finite lattice, and $K \leq L$. Then the following conditions are equivalent:*

(i) *$K$ is a congruence-reflecting sublattice of $L$.*

(ii) *$L$ is a congruence-reflecting extension of $K$.*

(iii) *Let $\mathfrak{p}$ and $\mathfrak{q}$ be prime intervals in $K$; if $\mathfrak{p} \Rightarrow \mathfrak{q}$ in $L$, then $\mathfrak{p} \Rightarrow \mathfrak{q}$ in $K$.*

As an example, the reader may want to verify that any sublattice of a distributive lattice is congruence-reflecting.

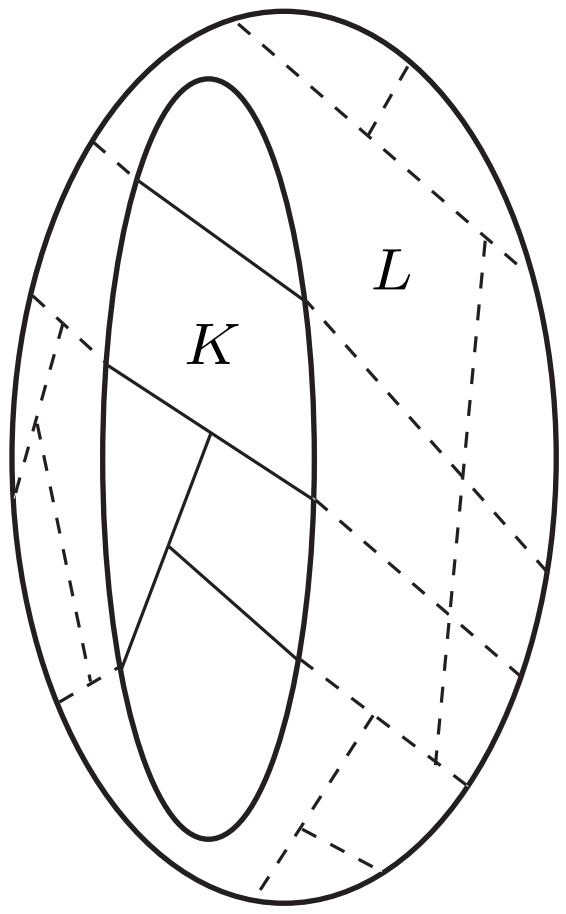

Figure 3.6: A congruence extends

A much stronger concept—central to this book—is the following. Let $K$ be a lattice. A lattice $L$ is a *congruence-preserving extension* of $K$ (or $K$ is a *congruence-preserving sublattice* of $L$), if $L$ is an extension and every congruence $\boldsymbol{\alpha}$ of $K$ has *exactly one* extension $\overline{\boldsymbol{\alpha}}$ to $L$ satisfying $\overline{\boldsymbol{\alpha}}\rceil_K = \boldsymbol{\alpha}$. Of course, $\overline{\boldsymbol{\alpha}} = \mathrm{con}_L(\boldsymbol{\alpha})$. It follows that $\boldsymbol{\alpha} \mapsto \mathrm{con}_L(\boldsymbol{\alpha})$ is an isomorphism between $\mathrm{Con}\, K$ and $\mathrm{Con}\, L$.

Two congruence-preserving extensions are shown in Figure 3.7, whereas Figure 3.8 shows two other extensions that are not congruence-preserving.

We can also obtain congruence-preserving extensions using gluing, based on the following result.

**Lemma 3.16.** *Let $K$ and $L$ be lattices, let $F$ be a filter of $K$, and let $I$ be an ideal of $L$. Let $\varphi$ be an isomorphism between $F$ and $I$. Let $G$ be the gluing of $K$ and $L$ over $F$ and $I$ with respect to $\varphi$. If $L$ is a congruence-preserving extension of $I$, then $G$ is a congruence-preserving extension of $K$.*

If $I$ and $L$ are simple, then $L$ is a congruence-preserving extension of $I$. So we obtain the following result.

**Corollary 3.17.** *Let $K$, $L$, $F$, $I$, and $\varphi$ be given as above. If $I$ and $L$ are simple lattices, then $G$ is a congruence-preserving extension of $K$.*

**Lemma 3.18.** *Let the lattice $L$ be an extension of the lattice $K$. Then $L$ is a congruence-preserving extension of $K$ iff the following two conditions hold*:

(i) $\mathrm{re}(\mathrm{ext}\,\boldsymbol{\alpha}) = \boldsymbol{\alpha}$ *for any congruence* $\boldsymbol{\alpha}$ *of* $K$;

(ii) $\mathrm{ext}(\mathrm{re}\,\boldsymbol{\alpha}) = \boldsymbol{\alpha}$ *for any congruence* $\boldsymbol{\alpha}$ *of* $L$.

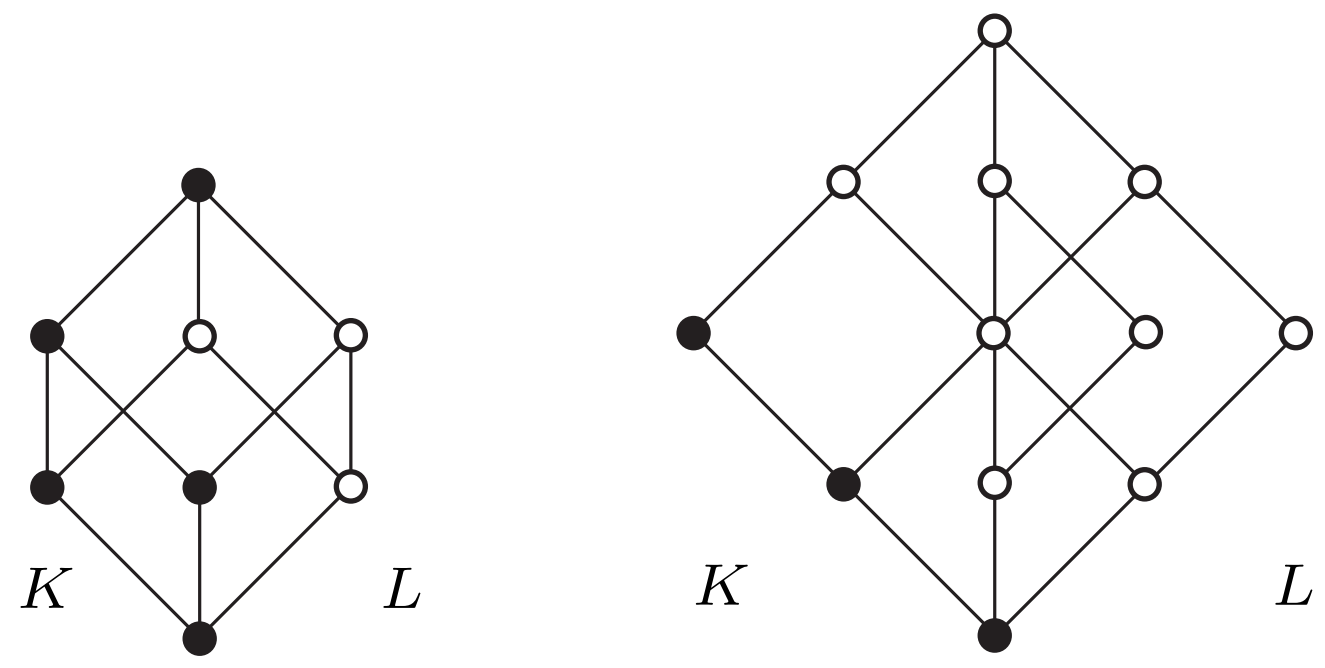

Figure 3.7: Examples of congruence-preserving extensions

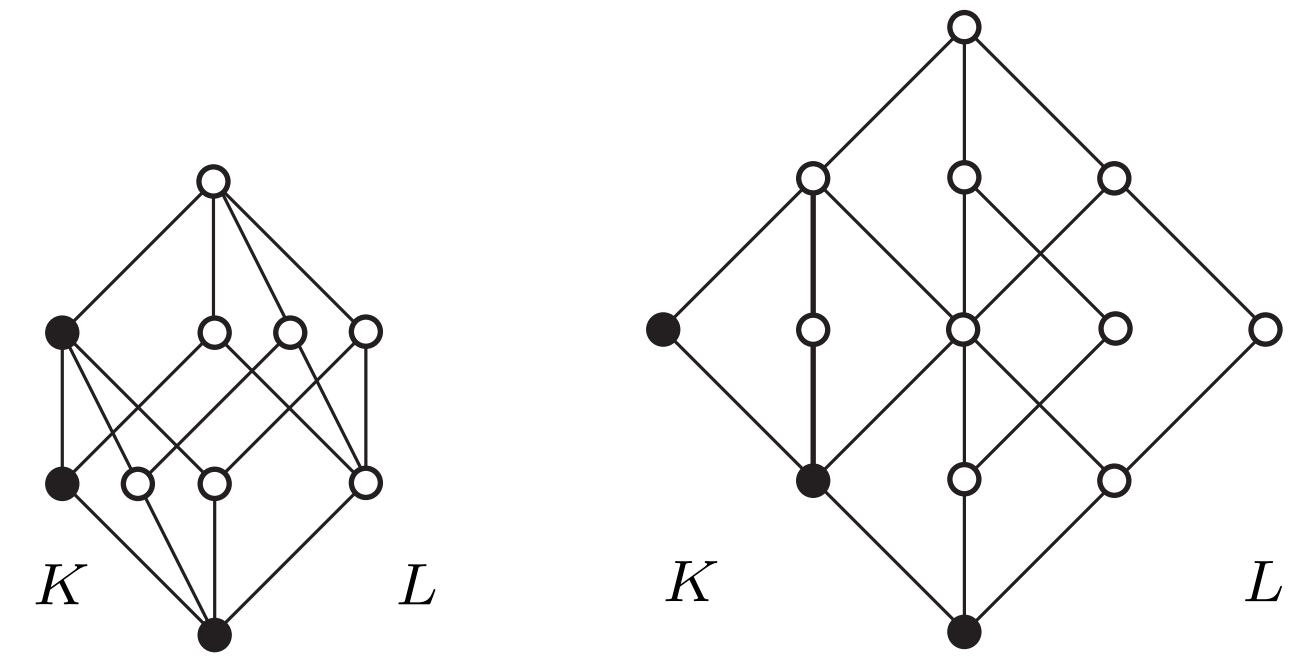

Figure 3.8: Extensions that are not congruence preserving

We can say a lot more for finite lattices.

**Lemma 3.19.** *Let $L$ be a finite lattice, and $K \leq L$. Then $L$ is a congruence-preserving extension of $K$ iff the following two conditions hold:*

(a) *Let $\mathfrak{p}$ and $\mathfrak{q}$ be prime intervals in $K$; if $\mathfrak{p} \Rightarrow \mathfrak{q}$ in $L$, then $\mathfrak{p} \Rightarrow \mathfrak{q}$ in $K$.*

(b) *Let $\mathfrak{p}$ be a prime interval of $L$. Then there exists a prime interval $\mathfrak{q}$ in $K$ such that $\mathfrak{p} \Leftrightarrow \mathfrak{q}$ in $L$.*

Lemma 3.18.(ii) is very interesting by itself. It says that every congruence $\boldsymbol{\alpha}$ of $L$ is determined by its restriction to $K$. In other words,

$$\boldsymbol{\alpha} = \operatorname{con}_L(\boldsymbol{\alpha}\rceil K).$$

We will call such a sublattice $K$ a *congruence-determining sublattice.* We can easily modify Lemma 3.19 to characterize congruence-determining sublattices in finite lattices.

**Lemma 3.20.** *Let $L$ be a finite lattice $L$, and $K \leq L$. Then $K$ is a congruence-determining sublattice of $L$ iff for any prime interval $\mathfrak{p}$ in $L$, there is a prime interval $\mathfrak{q}$ in $K$ satisfying $\mathfrak{p} \Leftrightarrow \mathfrak{q}$ in $L$.*

Of course, a congruence-preserving sublattice is always congruence-determining. In fact, a sublattice is congruence-preserving iff it is congruence-reflecting and congruence-determining.

**Chapter**

# 4

# *Planar Semimodular Lattices*

## 4.1. Planar lattices

A finite lattice $L$ is *planar* if it is planar as an ordered set (see Section 1.1.2), that is, it has a planar diagram as an ordered set. For planar lattices, the term prime interval is used interchangeably with edge. In Chapters 27 and 29, we use edges.

We have quite a bit of flexibility to construct a planar diagram for an ordered set, but we are much more constrained for a lattice because it has a zero, which must be the lowest element and a unit, which must be the highest element—contrast this with Figure 1.1. All lattices with five or fewer elements are planar; all but the five chains are shown in the first two rows of Figure 4.1.

The third row of Figure 4.1 provides an example of "good" and "bad" lattice diagrams; the two diagrams represent the same lattice, $\mathsf{C}_3^2$. Planar diagrams are the best. Diagrams in which meets and joins are hard to figure out are not of much value. For me, even Figure 1.5 is of limited value.

In the last row of Figure 4.1 there are two more diagrams. The one on the left is not planar; nevertheless, joins and meets are easy to see (the notation $\mathsf{M}_3[\mathsf{C}_3]$ will be explained in Section 6.1). The one on the right is not a lattice: the two black-filled elements have no join.

An *X-configuration* in a planar ordered set $P$ is formed by two edges $E$ and $F$ of $P$ satisfying the following properties:

(i) $0_E$ is to the left of $0_F$;

(ii) $1_E$ is to the right of $1_F$.

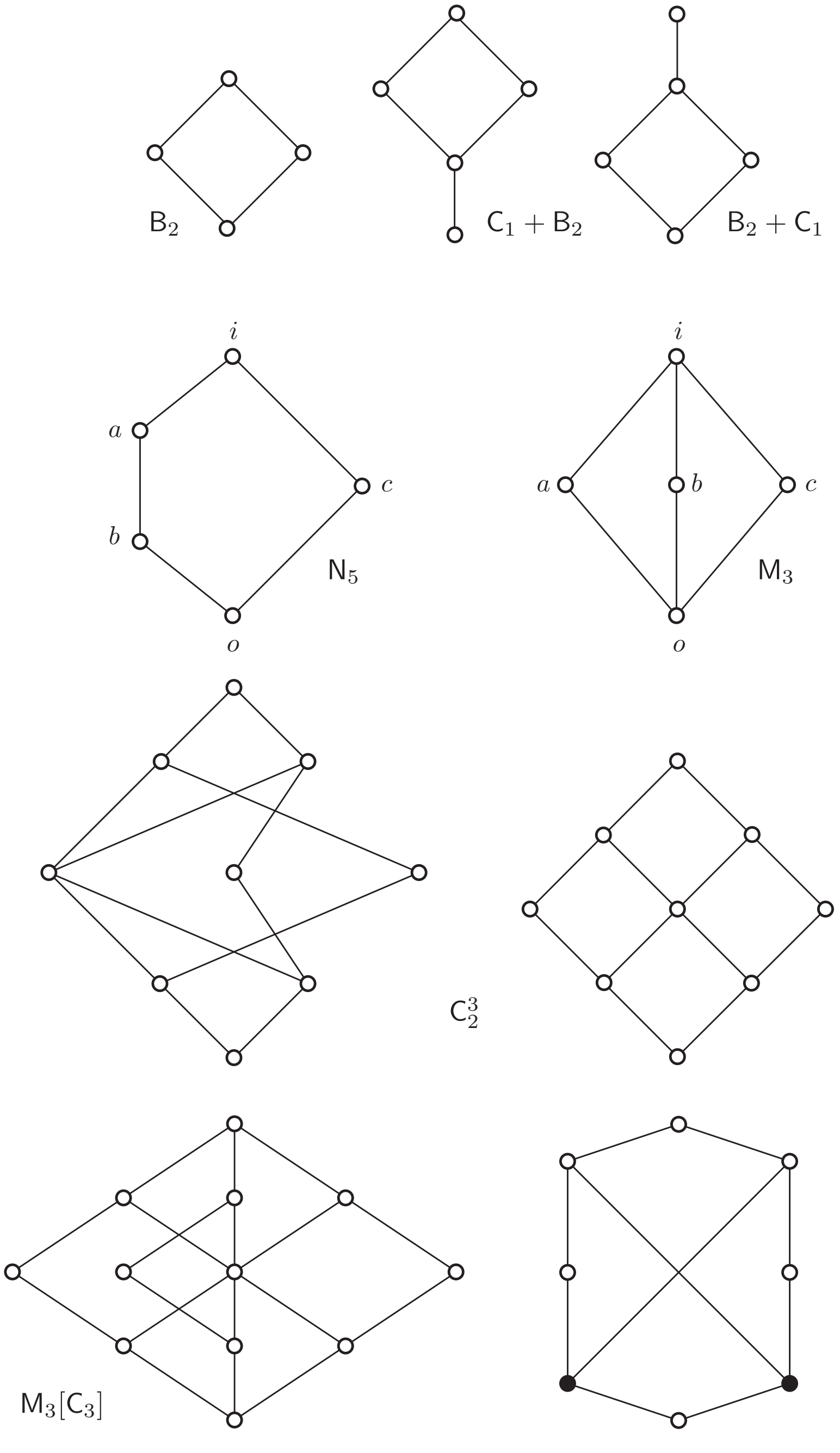

Figure 4.1: More diagrams

**Lemma 4.1.** *A diagram of a planar ordered set $P$ is the diagram of a planar lattice iff it does not have an X-configuration.*

This is almost a tautology.

Let us say that a finite ordered set $P$ with (strict) ordering relation $<$ has a complementary (strict) order relation $\lambda$, if any two distinct elements of $P$ are comparable by exactly one of $<$ and $\lambda$.

**Theorem 4.2** (Zilber's Theorem)**.** *A finite lattice $L$ is planar iff it has a complementary order relation.*

This result has immediate interesting ramifications, for instance, planar canonical forms of diagrams and their uniqueness (see my joint paper [14] with K. A. Baker).

A planar lattice $L$ has a *left boundary chain*, $\mathrm{C_l}(L)$, and a *right boundary chain*, $\mathrm{C_r}(L)$. A *left corner* lc(L) (or right corner rc(L)) of $L$ is a doubly-irreducible element in $L - \{0, 1\}$ on the left (or right) boundary of $L$. We define a *covering square* $S$ of a lattice $L$ as a four-element Boolean sublattice such that the coverings in $S$ are also coverings in $L$.

If $C$ and $D$ are maximal chains in the interval $[a, b]$, and there is no element of $L$ between $C$ and $D$, then we call $C \cup D$ a *cell*. A four-element cell is a 4-*cell*. A diagram of $\mathsf{M}_3$ has exactly two 4-cells and three covering squares. A planar lattice is called a 4-*cell lattice* if all of its cells are 4-cells. For example, $\mathsf{M}_3$ is a 4-cell lattice but $\mathsf{N}_5$ is not.

## 4.2. Two acronyms: SPS and SR

There are two classes of planar semimodular lattices we want to designate by acronyms.

As defined in my joint paper [134] with E. Knapp, a semimodular lattice $L$ is *slim*, if it contains no cover-preserving $\mathsf{M}_3$ sublattice. An *SPS lattice* is a slim, planar, semimodular lattice.

Following my joint paper [134] with E. Knapp, a semimodular lattice $L$ is *rectangular* if its left boundary chain, $\mathrm{C_l}(L)$, has exactly one doubly-irreducible element, lc(L), and its right boundary chain, $\mathrm{C_r}(D)$, has exactly one doubly-irreducible element, rc(D), and these elements are complementary, that is,

$$\mathrm{lc(D)} \vee \mathrm{rc(D)} = 1,$$
$$\mathrm{lc(D)} \wedge \mathrm{rc(D)} = 0.$$

An *SR lattice* is a slim rectangular lattice. By G. Czédli and E. T. Schmidt [74], whether a lattice is an SR lattice does not depend on the planar diagram chosen.

## 4.3. SPS lattices

Let $L$ be a planar semimodular lattice. An internal element of a cover-preserving $\mathsf{M}_3$ sublattice is called an *eye*. Removing all eyes one at a time, we get a planar semimodular lattice, called the *slimming* of $L$. We can also *add an eye*—in the obvious sense—to a 4-cell.

In an SPS lattice $L$, an element $u$ covering $n \geq 3$ elements together generate, up to isomorphism, a unique sublattice. For $n = 3$, it is $\mathsf{S}_7$, for $n = 4$, we get the lattice of Figure 4.2 (see G. Czédli [29]).

The following lemma is from my joint paper [134] with E. Knapp.

**Lemma 4.3.** *An SPS lattice $L$ is distributive iff $\mathsf{S}_7$ (see the Picture Gallery) is not a cover-preserving sublattice of $L$.*

Let $L$ be an SPS lattice. Two prime intervals of $L$ are *consecutive* if they are opposite sides of a 4-cell. As in G. Czédli and E. T. Schmidt [71], maximal sequences of consecutive prime intervals form a *trajectory*. So a trajectory is an equivalence class of the transitive reflexive closure of the consecutive relation (see Figure 4.3 for two examples). We denote by $\operatorname{Traj} L$ the set of all trajectories of $L$.

Let $C = \{o, a_l, a_r, i\}$ be a 4-cell of $L$ and let $E = [a_r, i]$ be the prime interval on the upper right of $C$. We go down by consecutive prime intervals until we reach the left boundary of $L$. These intervals form the *left wing* of $E$ and of the 4-cell $C$. We define the *right wing* symmetrically.

A simple trajectory goes from an upper boundary to a lower boundary; a more interesting trajectory has a *top edge*, which has a left and a right wing (see the first diagram of Figure 4.3).

**Lemma 4.4.** *Let $L$ be an SR lattice. Then for a join-irreducible congruence $\boldsymbol{\alpha}$ exactly one of the following conditions holds.*

(i) *There is a unique edge on the lower left boundary of $L$ generating $\boldsymbol{\alpha}$, but no such edge on the lower right boundary of $L$.*

(ii) *There is a unique edge on the lower right boundary of $L$ generating $\boldsymbol{\alpha}$, but no such edge on the lower left boundary of $L$.*

(iii) *There is a unique edge on the lower left boundary and a unique edge on the lower right boundary of $L$ both generating $\boldsymbol{\alpha}$.*

*For a rectangular lattice $L$, in general, the same statement holds except for the uniqueness.*

The first and third statement of the next lemma can be found in the literature (see my joint papers [134]–[138] with E. Knapp, and also G. Czédli and E. T. Schmidt [73]–[74]).

**Lemma 4.5.** *Let $L$ be an SPS lattice.*

(i) *An element of $L$ has at most two covers.*

(ii) *If the elements $a, u, v, w \in L$ satisfy $a = u \wedge v = v \wedge w = w \wedge u$, then two of them are comparable.*

(iii) *Let $x \in L$ cover three distinct elements $u, v$, and $w$. Then the set $\{u, v, w\}$ generates an $\mathsf{S}_7$ sublattice.*

Lemma 4.5(i) and (ii) state in different ways that there are only two directions "to go up" from an element. The next lemma states this in one more way (see G. Czédli and E. T. Schmidt [71]).

**Lemma 4.6.** *Let $L$ be an SPS lattice. Let $\mathfrak{q}, \mathfrak{q}_1, \mathfrak{q}_2$ be pairwise distinct prime intervals of $L$ satisfying $\mathfrak{q}_1 \overset{\text{dn}}{\sim} \mathfrak{q}$ and $\mathfrak{q}_2 \overset{\text{dn}}{\sim} \mathfrak{q}$. Then $\mathfrak{q}_1 \sim \mathfrak{q}_2$.*

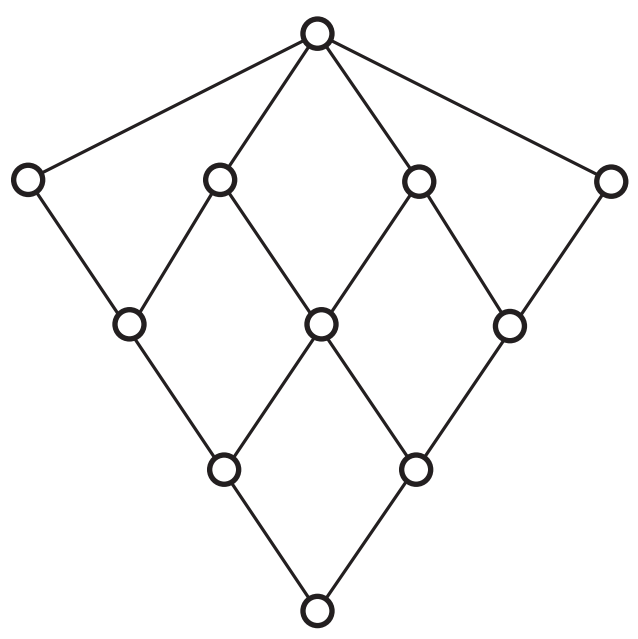

Figure 4.2: A special semimodular lattice

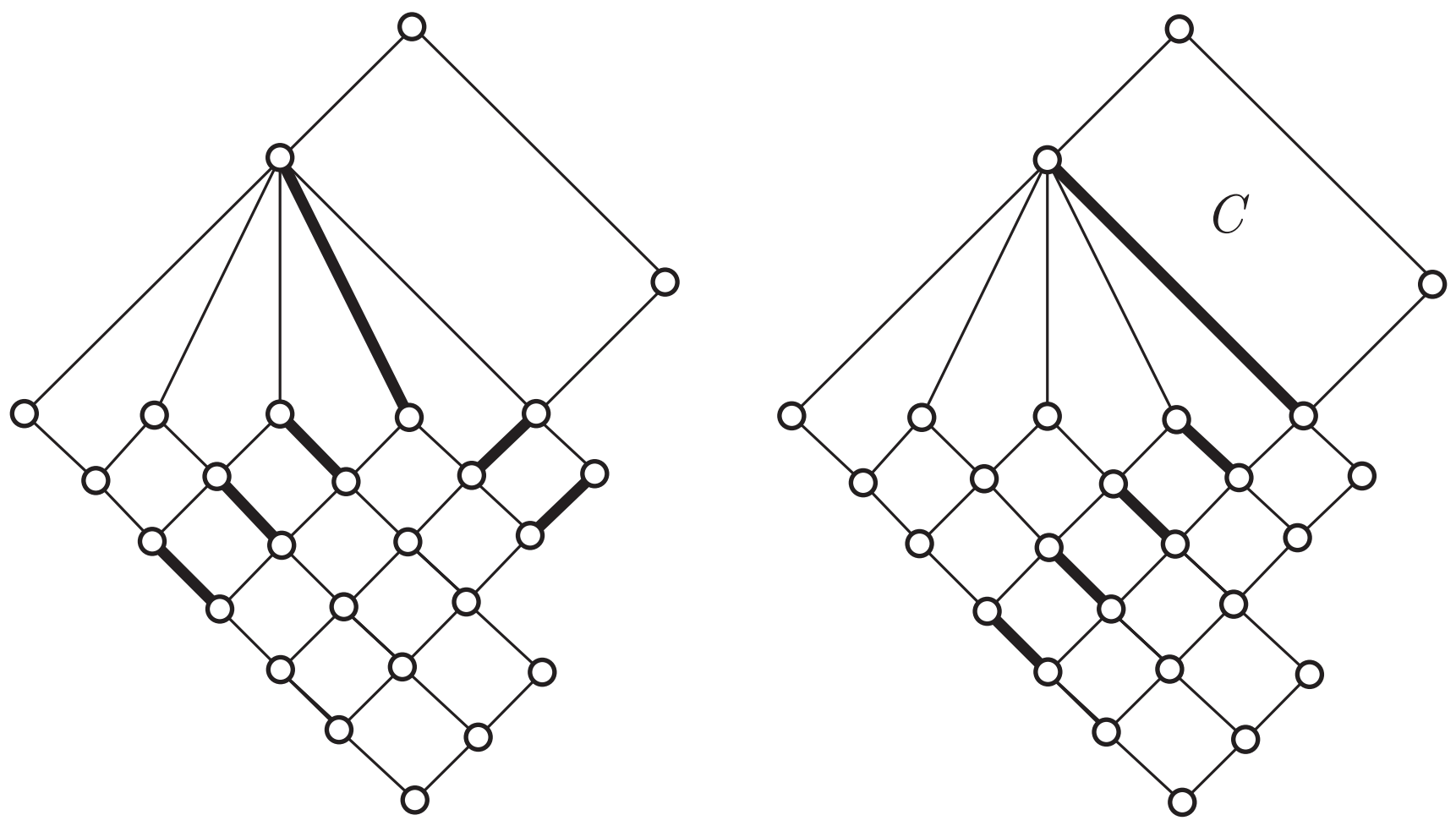

Figure 4.3: A trajectory and a left wing of a 4-cell

Call a lattice $L$ *meet-semidistributive*, if the condition

(SD$_\wedge$) $x \wedge y = x \wedge z$ implies that $x \wedge y = x \wedge (y \vee z)$

holds. Dually, we define (SD$_\vee$) and *join-semidistributive.*

SPS lattices are meet-semidistributive (see LTF-[105] for more discussion, especially Section VI.2.7 by K. Adaricheva). For an overview of semidistributivity, see Chapters 3–6 ([8]–[11]) by K. Adaricheva and J. B. Nation in LTS2-[207].

## 4.4. Forks

Our goal is to construct all planar semimodular lattices from planar distributive lattices, based on G. Czédli and E. T. Schmidt [73].

Let $L$ be an SPS lattice and let $C$ be a 4-cell in $L$. We construct a lattice extension $L[C]$ of $L$ as follows.

We observe that the lower left edge of $C$ is in an interval $A = \mathsf{C}_2 \times \mathsf{C}_n$, stretching from the lower left edge of $C$ to the left boundary of $L$, see the gray-filled elements in Figure 4.4. We replace the interval $A$ with the interval $B = \mathsf{C}_3 \times \mathsf{C}_n$ (as in the third diagram of Figure 4.4, the new elements are black-filled). We do the same on the right side of $C$ and add one more element in $C$, the join of the two top black-filled elements.

Let $L[C]$ denote the lattice we obtain. We say that $L[C]$ is obtained from $L$ by *inserting a fork* at the 4-cell $C$. It is easy to see that $L[C]$ is an SPS lattice. For the new elements we shall use the notation of Figure 4.5.

The direct product of two nontrivial chains is a *grid.*

Now we can state the main result of G. Czédli and E. T. Schmidt [73].

**Theorem 4.7.** *An SPS lattice with at least three elements can be constructed from a grid by the following two steps.*

(i) *Inserting forks.*

(ii) *Followed by removing corners.*

My joint paper [57] with G. Czédli (see also Chapter 3 in LTS1-[206]) presents a twin of this construction, called *resection.*

Sometimes, we can *delete* a fork (see G. Czédli and E. T. Schmidt [73], as illustrated by Figure 4.6).

**Lemma 4.8.** *Let $L$ be an SPS lattice and let $S$ be a cover-preserving* $\mathsf{S}_7$ *in $L$. Let us assume that the top element $t$ of $S$ is* minimal, *that is, there is no $S'$ a covering* $\mathsf{S}_7$ *with top element $t'$ satisfying that $t' < t$. Then $L$ has a sublattice $L^-$ with* 4*-cell $C = S - \{m, b_l, b_r\}$ such that $L = L^-[C]$.*

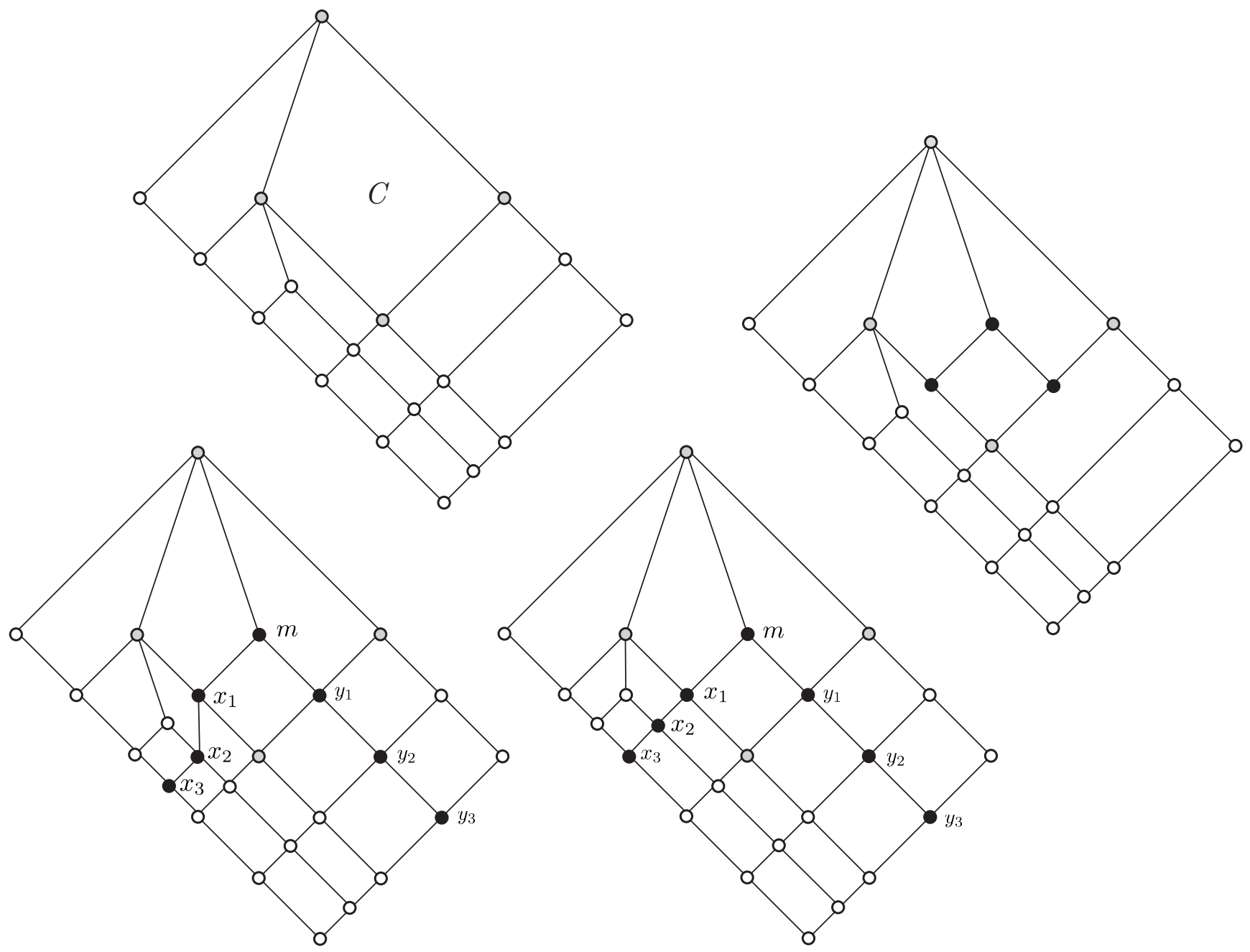

Figure 4.4: Inserting a fork into $L$ at $C$ (the third and fourth diagrams represent the same lattice)

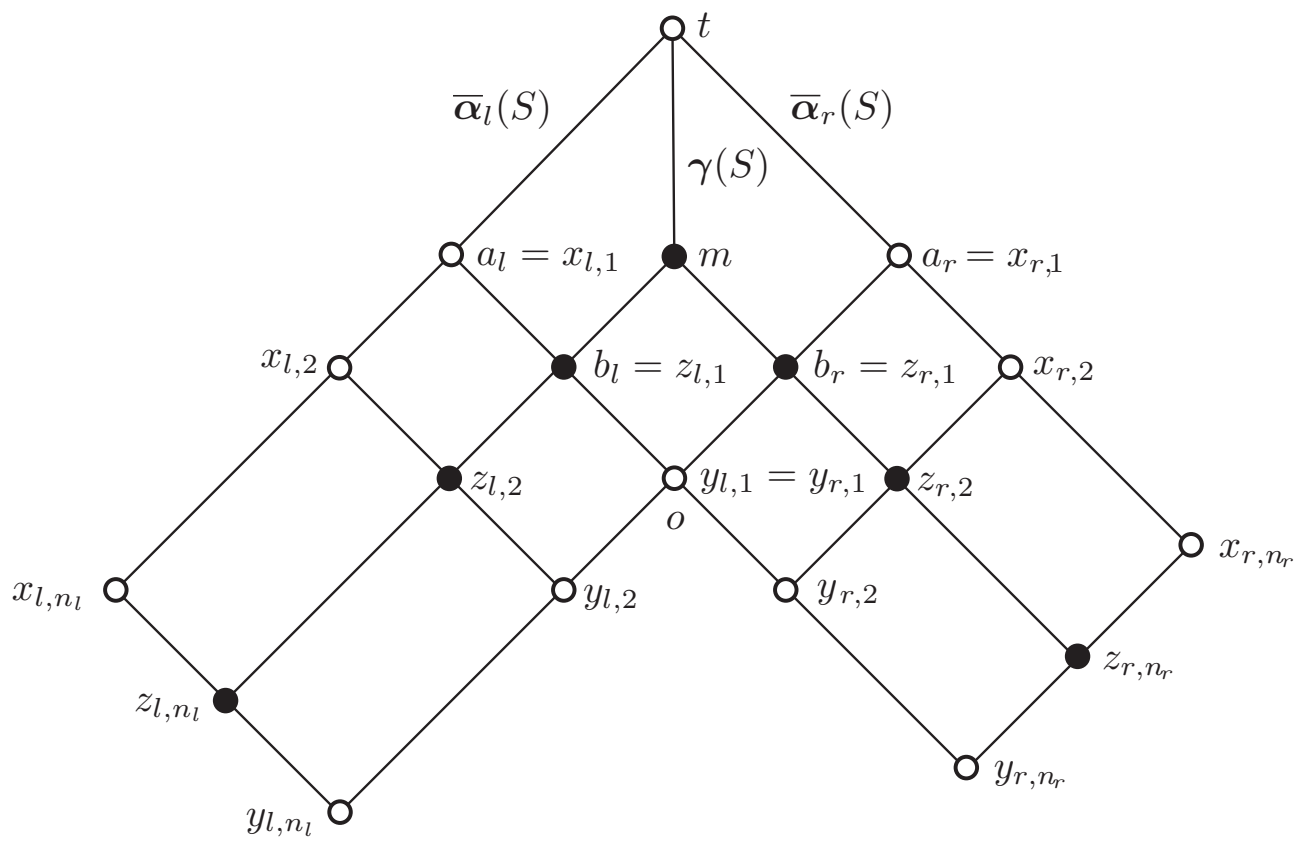

Figure 4.5: Notation for the fork construction

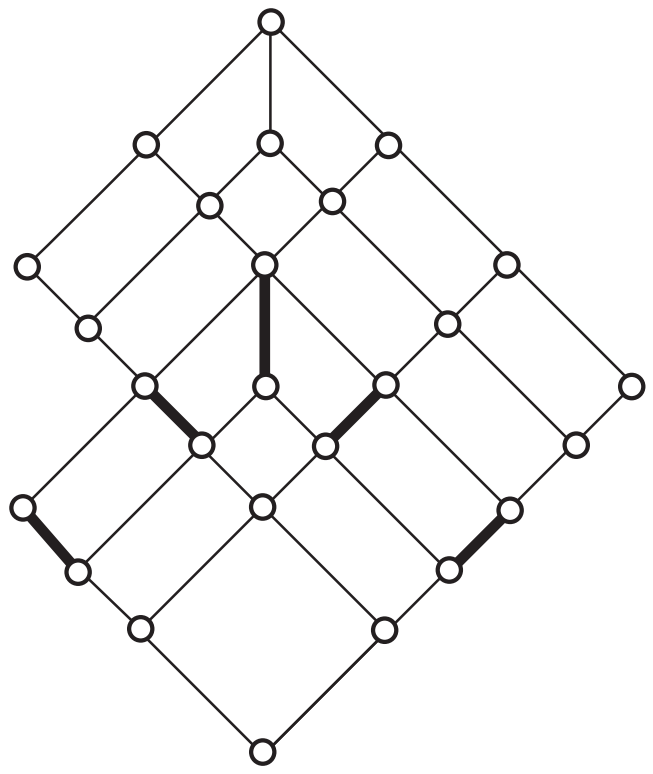

Figure 4.6: Deleting a fork (a trajectory)

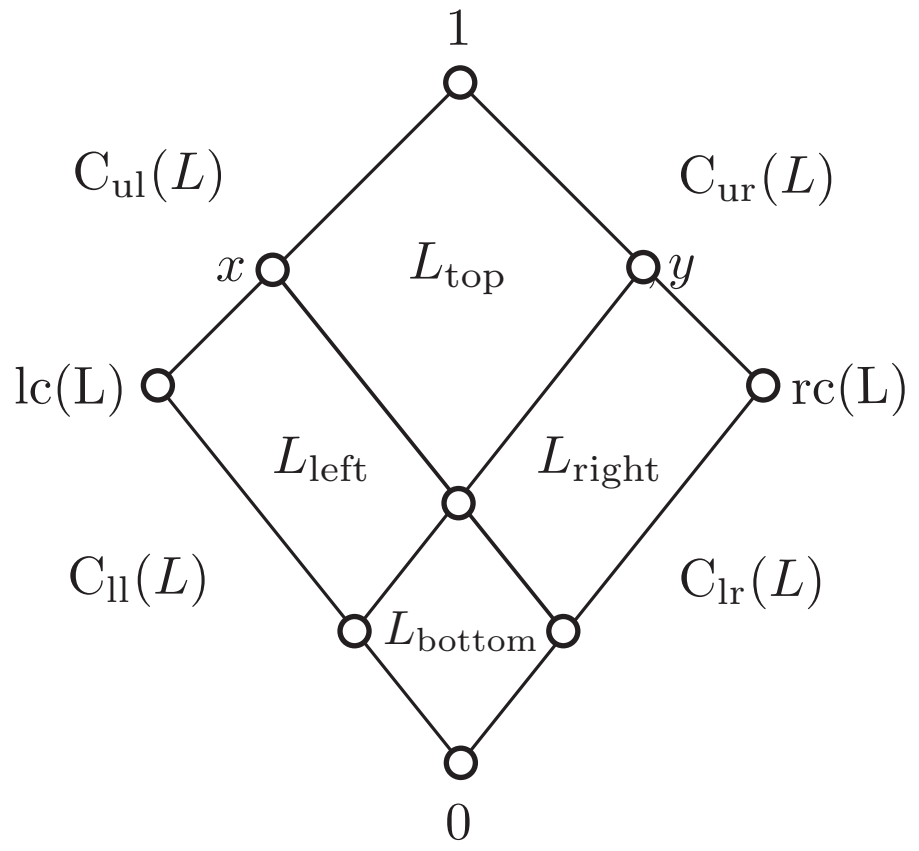

Figure 4.7: Decomposing an SR lattice

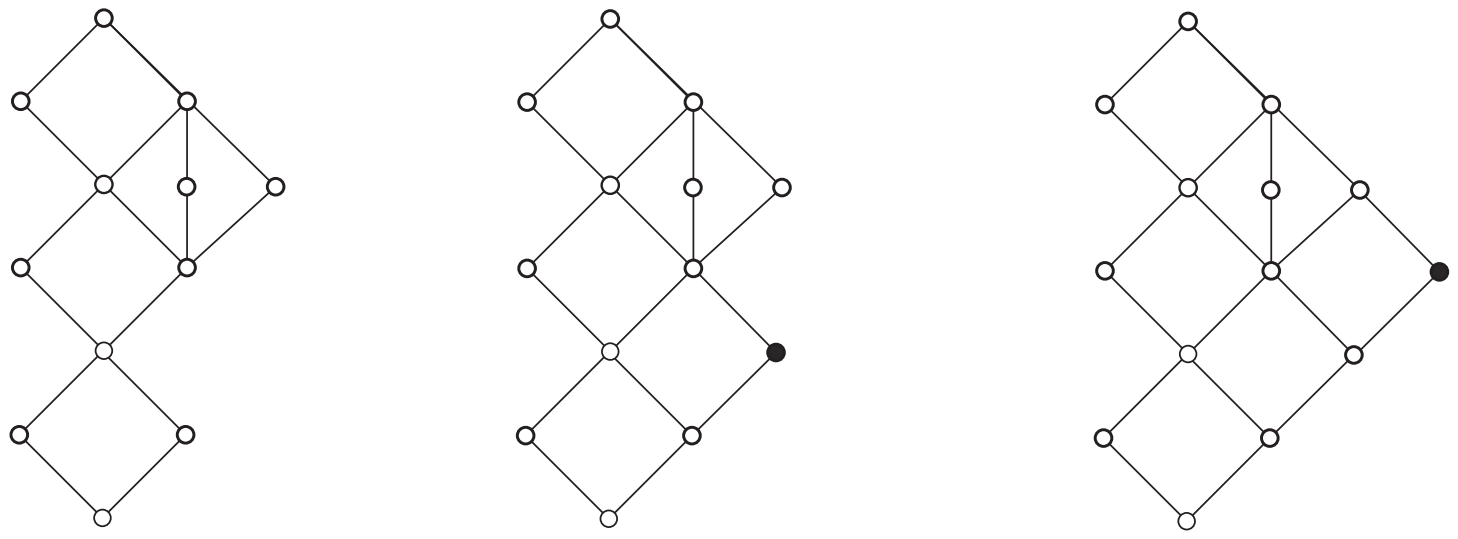

Figure 4.8: Two steps to a congruence-preserving rectangular extension

There are other constructions that yield all planar semimodular lattices. SPS lattices can also be described by *Jordan–Hölder permutations* (see the papers G. Czédli and E. T. Schmidt [75], G. Czédli, L. Ozsvárt, and B. Udvari [67], and G. Czédli [34]). These descriptions are generalized to a larger class of lattices in G. Czédli [31] and K. Adaricheva and G. Czédli [4] (see also K. Adaricheva, R. Freese, and J. B. Nation [5]).

There is a description with matrices in G. Czédli [27] and smaller diagrams in G. Czédli [37]; they can also be built from smaller building blocks called rectangular lattices, which are discussed in the next section.

## 4.5. Rectangular lattices

The following two lemmas are from my joint paper [138] with E. Knapp.

**Lemma 4.9.** *Let $L$ be a rectangular lattice. Then the intervals* $[0, \mathrm{lc(L)}]$, $[\mathrm{lc(L)}, 1]$, $[0, \mathrm{rc(L)}]$, *and* $[\mathrm{rc(L)}, 1]$ *are chains.*

So the chains $\mathrm{C_l}(L)$ and $\mathrm{C_r}(L)$ are split into two, a lower and an upper part: $\mathrm{C_{ll}}(L) = [0, \mathrm{lc(L)}]$, $\mathrm{C_{ul}}(L) = [\mathrm{lc(L)}, 1]$, $\mathrm{C_{lr}}(L) = [0, \mathrm{rc(L)}]$, and $\mathrm{C_{ur}}(L) = [\mathrm{rc(L)}, 1]$ (see Figure 4.7).

Lemma 4.9 has an interesting form for SR lattices.

**Lemma 4.10.** *Let $L$ be an SR lattice. Then for any $x \in L$, the following equation holds:*

$$x = (x \wedge \mathrm{lc(L)}) \vee (x \wedge \mathrm{rc(L)}).$$

*Moreover, the map* $\varphi\colon x \to (x \wedge \mathrm{lc(L)}, x \wedge \mathrm{rc(L)})$ *embeds $L$ into the grid* $\mathrm{C_l}(L) \times \mathrm{C_r}(L)$.

The structure of rectangular lattices is easily described utilizing Theorem 4.7 (see G. Czédli and E. T. Schmidt [73]).

**Theorem 4.11** (Structure Theorem for SR Lattices)**.** *$L$ is an SR lattice iff it can be obtained from a grid by inserting forks.*

There is a slightly stronger version of this result, implicit in G. Czédli and E. T. Schmidt [69].

**Theorem 4.12** (Structure Theorem, Strong Version)**.** *For every SR lattice $K$, there is a grid $G$ and sequences*

$$G = K_1, K_2, \dots, K_{n-1}, K_n = K$$

*of SR lattices and*

$$C_1 = \{o_1, c_1, d_1, i_1\}, C_2 = \{o_2, c_2, d_2, i_2\}, \dots, C_{n-1} = \{o_{n-1}, c_{n-1}, d_{n-1}, i_{n-1}\}$$

*of* 4*-cells in the appropriate lattices such that*

$$G = K_1, K_1[C_1] = K_2, \dots, K_{n-1}[C_{n-1}] = K_n = K.$$

*Moreover, the principal ideals* $\mathrm{id}(c_{n-1})$ *and* $\mathrm{id}(d_{n-1})$ *are distributive.*

*Proof.* We prove this result by induction on the number $n$ of cover-preserving $\mathsf{S}_7$-s in $K$. If $n = 0$, then $K$ is distributive by Lemma 4.3, so the statement is trivial. Now let us assume that the statement holds for $n - 1$. Let $K$ be an SR lattice with $n$ covering $\mathsf{S}_7$-s. As in Lemma 4.8, we take $S$, a minimal covering $\mathsf{S}_7$ in $K$. Then we form the sublattice $K^-$ by deleting the fork at $S$. So we get a 4-cell $C = C_{n-1} = \{o_{n-1}, c_{n-1}, d_{n-1}, i_{n-1}\}$ of $K^-$ such that $K = K^-[C]$. Since $K^-$ has $n - 1$ covering $\mathsf{S}_7$-s, we get the sequence

$$G = K_1, K_1[C_1] = K_2, \ldots, K_{n-2}[C_{n-2}] = K_{n-1} = K^-,$$

which, along with $K = K^-[C]$, proves the statement for $K$.

By the minimality of $S$, it follows from Lemma 4.3 that the principal ideals $\mathrm{id}(c_{n-1})$ and $\mathrm{id}(d_{n-1})$ are distributive. □

The following result of G. Czédli [29] is a stronger version of the Structure Theorems.

**Theorem 4.13.** *Every SR lattice is obtained from a grid by a sequence of fork insertions at distributive* 4*-cells, and every diagram obtained this way is an SR diagram.*

A rectangular lattice $L$ is a *patch lattice*, if the corners are dual atoms. Let $L$ be a nontrivial lattice. If $L$ cannot be obtained as a gluing of two nontrivial lattices, we call $L$ *gluing indecomposable.*

The following is a result of G. Grätzer and E. Knapp [138] (see also G. Czédli and E. T. Schmidt [73]); it is an easy consequence of Theorem 4.17.

**Theorem 4.14.** *Let $L$ be an SPS lattice with at least four elements. Then $L$ is a patch lattice iff it is gluing indecomposable.*

The next result of G. Czédli and E. T. Schmidt [74] follows easily from Theorem 4.11.

**Theorem 4.15.** *A patch lattice $L$ can be obtained from the four-element Boolean lattice by inserting forks.*

My paper [106] provides an alternative approach to the results of G. Czédli and E. T. Schmidt [74].

Finally, start with a planar semimodular lattice $L$. Can we add a corner? Yes, unless $L$ is rectangular. So we get the following result, The following two lemmas are from my joint paper [138] with E. Knapp (illustrated by Figure 4.8). Three more steps are needed to get a rectangular extension.

**Lemma 4.16.** *Let $K$ be a planar semimodular lattice. Then $K$ has a rectangular extension $L$. In fact, we can construct $L$ as a congruence-preserving extension of $K$.*

Note that under reasonable additional assumptions, the lattice $K$ (and even its diagram) is unique, as formulated in G. Czédli [40].

**Triple gluing**

For an SR lattice $L$, let

$$x \in \mathrm{C_{ul}}(L) - \{1, \mathrm{lc(L)}\},$$
$$y \in \mathrm{C_{ur}}(L) - \{1, \mathrm{rc(L)}\}.$$

We introduce some notation (see Figure 4.7).

$$\begin{aligned} L_{\mathrm{top}}(x,y) &= [x \wedge y, 1], \\ L_{\mathrm{left}}(x,y) &= [\mathrm{lc(L)} \wedge y, x], \\ L_{\mathrm{right}}(x,y) &= [x \wedge \mathrm{rc(L)}, y], \\ L_{\mathrm{bottom}}(x,y) &= [0, (\mathrm{lc(L)} \wedge y) \vee (x \wedge \mathrm{rc(L)})] \end{aligned} \tag{1}$$

The following decomposition theorem is from my joint paper [138] with E. Knapp.

**Theorem 4.17** (Decomposition Theorem for SR Lattices)**.** *Let $L$ be an SR lattice, and let*

$$x \in \mathrm{C_{ul}}(L) - \{1, \mathrm{lc(L)}\}, \tag{2}$$
$$y \in \mathrm{C_{ur}}(L) - \{1, \mathrm{rc(L)}\}. \tag{3}$$

*Then $L$ can be decomposed into four SR lattices*

$$L_{\mathrm{top}}(x,y), L_{\mathrm{left}}(x,y), L_{\mathrm{right}}(x,y), L_{\mathrm{bottom}}(x,y),$$

*and the lattice $L$ can be reconstructed from these by repeated gluing.*

We now reformulate Theorem 4.17. As in my joint paper [152] with H. Lakser, let $G, Y, Z$, and $U$ be rectangular lattices, arranged in the plane as in Figure 4.9, such that the facing boundary chains have the same number of elements. First, we glue $U$ and $Y$ together over the boundaries, to obtain the rectangular lattice $X$, then we glue $Z$ and $G$ together over the boundaries, to obtain the rectangular lattice $W$. Finally, we glue $X$ and $W$ together over the boundaries, to obtain the rectangular lattice $V$, which we call the *triple gluing* of $G, Y, Z, U$. Now we restate Theorem 4.17 using triple gluing, using the notation (1).

**Theorem 4.18.** *Let $L$ be an SR lattice and let*

$$x \in \mathrm{C_{ul}}(L) - \{1, \mathrm{lc(L)}\}, \tag{4}$$
$$y \in \mathrm{C_{ur}}(L) - \{1, \mathrm{rc(L)}\}. \tag{5}$$

*Then $L$ is the triple gluing of the four intervals*

$$[x \wedge y, 1], [\mathrm{lc(L)} \wedge y, x], [x \wedge \mathrm{rc(L)}, y], [0, (\mathrm{lc(L)} \wedge y) \vee (x \wedge \mathrm{rc(L)})].$$

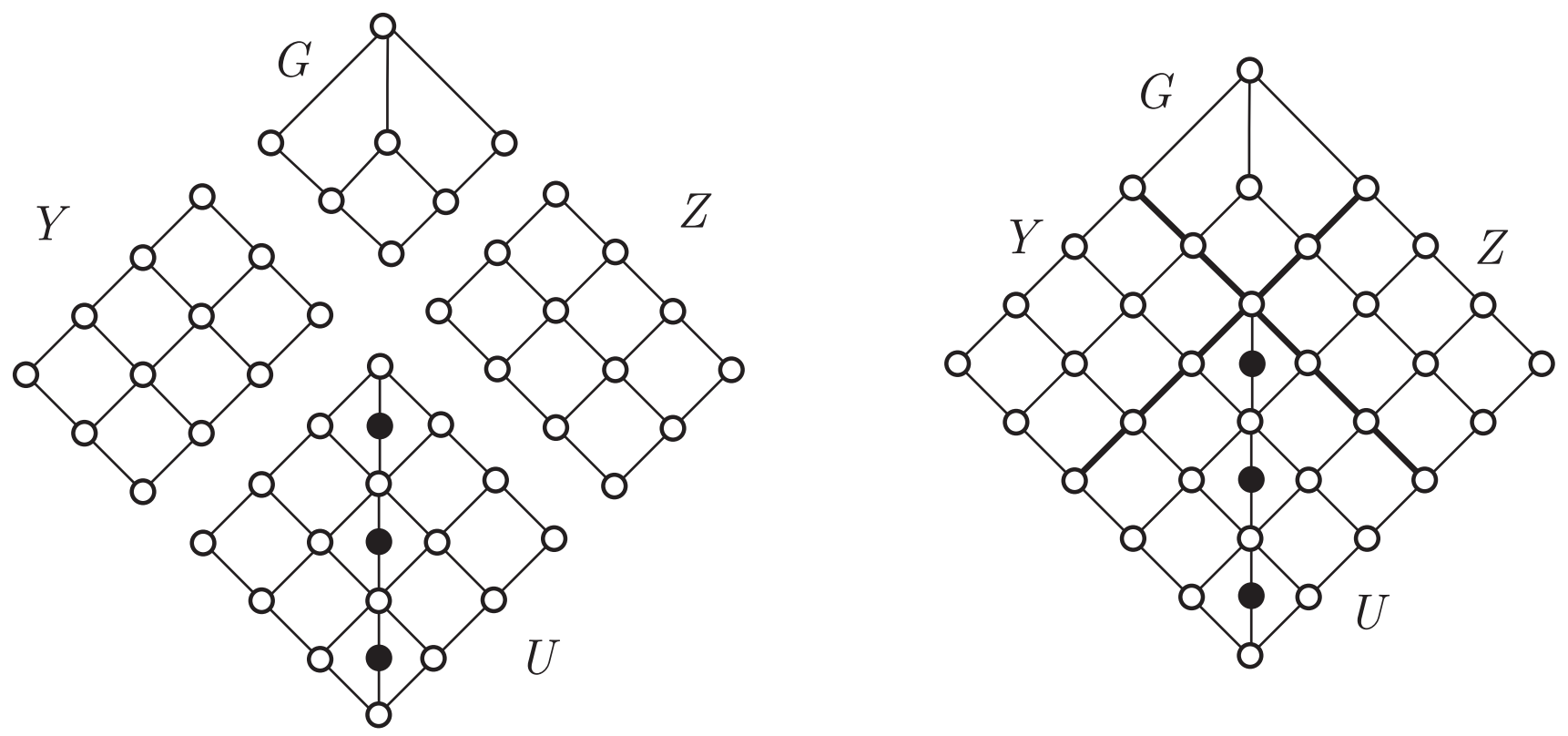

Figure 4.9: Triple gluing

The basic results on gluing (see Section 2.4) extend to triple gluings. We will use Lemma 2.9 in the following form (see my joint paper [152] with H. Lakser).

**Lemma 4.19.** *A congruence* $\boldsymbol{\alpha}$ *of* $V$ *is uniquely associated with the four congruences* $\boldsymbol{\alpha}_G$ *of* $G$, $\boldsymbol{\alpha}_Y$ *of* $Y$, $\boldsymbol{\alpha}_Z$ *of* $Z$, *and* $\boldsymbol{\alpha}_U$ *of* $U$, *satisfying the condition that* $\boldsymbol{\alpha}_G$ *and* $\boldsymbol{\alpha}_Y$ *agree on the facing boundaries, and the same for* $\{G, Z\}$, $\{Y, U\}$, *and* $\{Z, U\}$.

### Congruences of rectangular lattices

There is a "coordinatization" of congruences of rectangular lattices.

**Theorem 4.20.** *Let* $L$ *be a rectangular lattice and let* $\boldsymbol{\alpha}$ *be a congruence of* $L$. *Let* $\boldsymbol{\alpha}^l$ *denote the restriction of* $\boldsymbol{\alpha}$ *to* $C_{\mathrm{ll}}$. *Let* $\boldsymbol{\alpha}^r$ *denote the restriction of* $\boldsymbol{\alpha}$ *to* $C_{\mathrm{lr}}$. *Then the congruence* $\boldsymbol{\alpha}$ *is determined by the pair* $(\boldsymbol{\alpha}^l, \boldsymbol{\alpha}^r)$. *In fact,*

$$\boldsymbol{\alpha} = \mathrm{con}(\boldsymbol{\alpha}^l \cup \boldsymbol{\alpha}^r).$$

*Proof.* Since $\boldsymbol{\alpha} \geq \mathrm{con}(\boldsymbol{\alpha}^l \cup \boldsymbol{\alpha}^r)$, it is sufficient to prove that

(P) if the prime interval $\mathfrak{p}$ of $A$ is collapsed by the congruence $\boldsymbol{\alpha}$, then it is collapsed by the congruence $\mathrm{con}(\boldsymbol{\alpha}^l \cup \boldsymbol{\alpha}^r)$.

First, let $L$ be a slim patch lattice. By Theorem 4.15, we obtain $L$ from the square, $\mathsf{C}_2^2$, with a sequence of $n$ fork insertions. We induct on $n$.

If $n = 0$, then $L = \mathsf{C}_2^2$, and the statement is trivial.

Let the statement hold for $n - 1$ and let $K$ be the patch lattice we obtain by $n - 1$ fork insertions into $\mathsf{C}_2^2$, so that we obtain $L$ from $K$ by one fork insertion at the covering square $S$. We have three cases to consider.

**Case 1.** $\mathfrak{p}$ is a prime interval of $K$. Then the statement holds for $\mathfrak{p}$ and $\boldsymbol{\alpha}\rceil_K$, the restriction of $\boldsymbol{\alpha}$ to $K$ by induction. So $\mathfrak{p}$ is collapsed by $\operatorname{con}((\boldsymbol{\alpha}\rceil_K)^l \cup (\boldsymbol{\alpha}\rceil_K)^r)$ in $K$. Therefore, (P) holds in $L$.

In the next two cases, we assume that $\mathfrak{p}$ is not in $K$.

**Case 2.** $\mathfrak{p}$ is perspective to a prime interval of $K$. Same proof as in Case 1. This case includes $\mathfrak{p} = [o, a]$ and any one of the new intervals up-perspective with $[o, a]$.

**Case 3.** $\mathfrak{p} = [a, c]$ and any one of the new intervals is up-perspective with $[a, c]$. Then the fork extension defines the terminating prime interval $\mathfrak{q} = [y, z]$ on the boundary of $L$, which is up-perspective with $\mathfrak{p}$, verifying (P).

Second, let $L$ be a patch lattice, not necessarily slim. This case is obvious because (P) is preserved when inserting an eye.

Finally, if $L$ is not a patch lattice, we induct on $|L|$. By Theorem 4.14, $L$ is the rectangular gluing of the rectangular lattices $A$ and $B$ over the ideal $I$ and filter $J$. Let $\mathfrak{p}$ be a prime interval of $L$. Then $\mathfrak{p}$ is a prime interval of $A$ or $B$, say, of $A$. (If $\mathfrak{p}$ is a prime interval of $B$, then the argument is easier.) By induction, $\mathfrak{p}$ is collapsed by $\operatorname{con}(\boldsymbol{\alpha}\rceil_{C_{\mathrm{ll}}(A)} \cup \boldsymbol{\alpha}\rceil_{C_{\mathrm{lr}}(A)})$, so it is collapsed by $\operatorname{con}(\boldsymbol{\alpha}\rceil_{C_{\mathrm{ll}}(L)} \cup \boldsymbol{\alpha}\rceil_{C_{\mathrm{lr}}(L)}) = \operatorname{con}(\boldsymbol{\alpha}^l \cup \boldsymbol{\alpha}^r)$. □

Congruences of SR lattices are not that different from congruences of SPS lattices.

Finally, start with a planar semimodular lattice $L$. Can we add a corner? Yes, unless $L$ is rectangular. So we get the following result, illustrated by Figure 4.8. Three more steps are needed there to get a rectangular extension.

**Lemma 4.21.** *Let $K$ be a planar semimodular lattice. Then $K$ has a rectangular extension $L$. In fact, $L$ is a congruence-preserving extension of $K$.*

## 4.6. Rectangular intervals

We call the interval $I = [o, i]$ of an SPS lattice $L$ *rectangular*, if there are complementary $a, b \in I$ such that the element $a$ is to the left of the element $b$.

The next result is a new way to construct an SR lattice from an SPS lattice (see my paper [126]).

**Theorem 4.22.** *Let $L$ be an SPS lattice and let $I$ be a rectangular interval of $L$. Then the lattice $I$ is also an SR lattice.*

*Proof.* Theorem 20.3 obviously holds for grids.

Otherwise, we can assume that the SR lattice $K$ is not a grid, so $n = \mathrm{Rank(K)} > 1$. Let $K^-$ be the lattice we obtain by deleting a minimal fork in $K^-$ at the covering square

$$C_{n-1} = \{o_{n-1}, c_{n-1}, d_{n-1}, i_{n-1}\}.$$

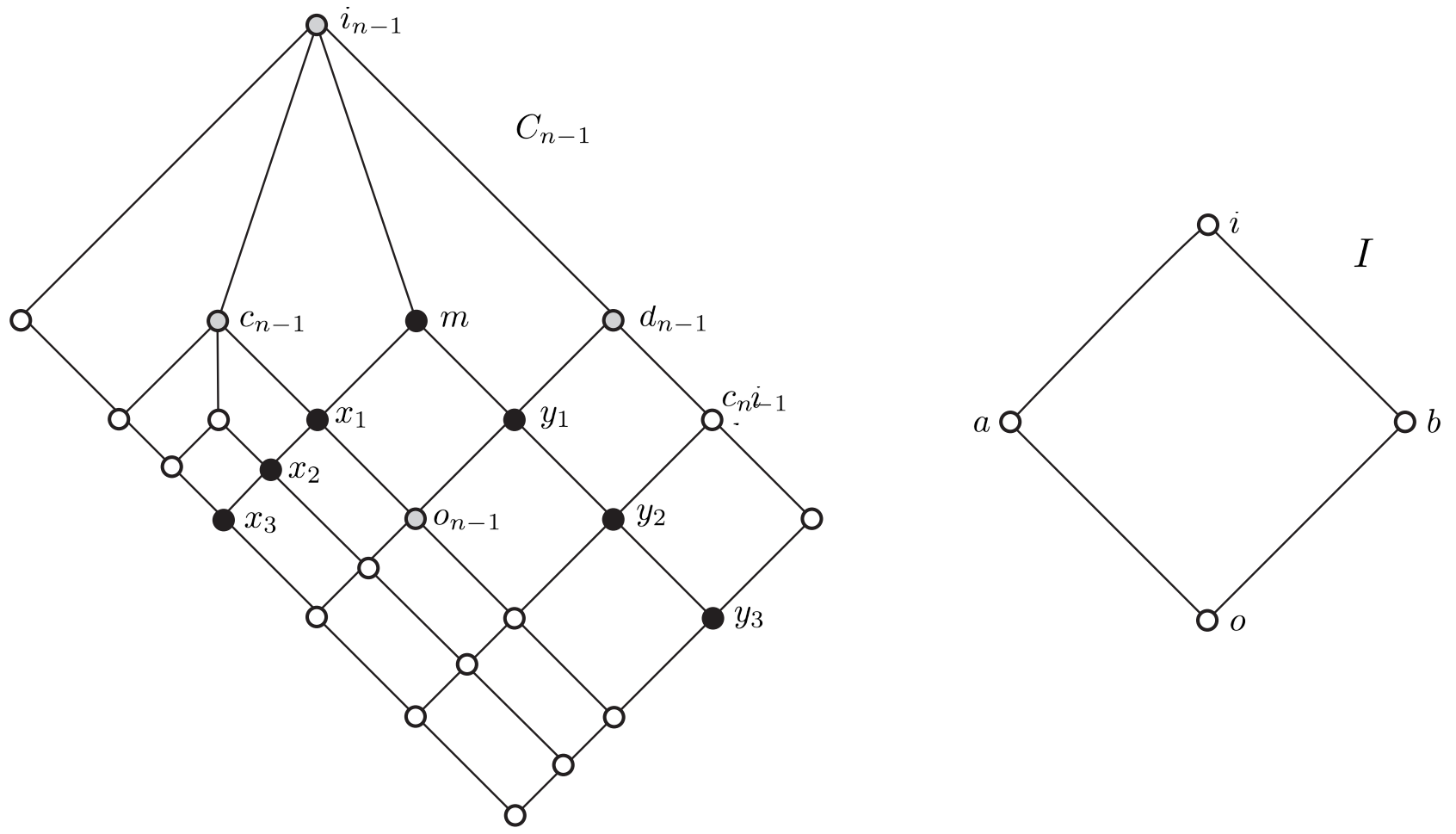

Figure 4.10: Proving Theorem 4.22: Case 1

We obtain $K$ from $K^-$ by inserting a fork at $C_{n-1}$. We add the element $m$ in the middle of $C_{n-1}$, and add the sequences of elements $x_1, \ldots$ on the left going down and $y_1, \ldots$ on the right going down as in Figure 4.4.

Let $I$ be a rectangular interval in $K$ with corners $a, b$, where $a$ is to the left of $b$. We want to prove that $I$ is an SR lattice. Of course, the lattice $I$ is slim.

We induct on $n = \text{Rank(K)}$. There are three subcases.

**Case 1.** $I$ is disjoint to $\text{id}(m)$, as illustrated in Figure 4.10. Then the interval $I$ is not changed as we add the fork to $K^-$. By induction, $I$ is rectangular in $K^-$, therefore, $I$ is also rectangular in $K$.

**Case 2.** In Figure 4.11 (and Figure 4.12), the bold lines form the boundary of the rectangular sublattice $I$ in $K^-$, the elements of $C_{n-1}$ are gray-filled, and the elements $m$, $x_1$, ..., $y_1$, ... are black-filled. The element $m$ is internal in $I$, so the element $a$ is $c_{n-1}$ or it is to the left of $c_{n-1}$, and symmetrically (see Figure 4.11). Therefore, $C_{n-1} = [o_{n-1}, i_{n-1}]_{K^-}$ is a covering square in $K^-$ and we obtain the interval $[o_{n-1}, i_{n-1}]_K$ of $K$ by adding a fork to $C_{n-1}$ at $[o_{n-1}, i_{n-1}]_{K^-}$. A fork extension of an SR lattice is also an SR lattice, so we conclude that $I$ is an SR lattice.

**Case 3.** $m$ is not an internal element of $I$ but some $x_i$ or $y_i$ is (see Figure 4.12), where $y_2$ is an internal element of $I$. By utilizing that $\text{id}(d_{n-1})$ is distributive, we conclude that we obtain $I$ from $[o, i]_{K^-}$ by replacing a cover preserving $\mathsf{C}_m \times \mathsf{C}_2$ by $\mathsf{C}_m \times \mathsf{C}_3$, and so $I$ remains rectangular. □

We proceed with some applications of this result. The next statement follows directly.

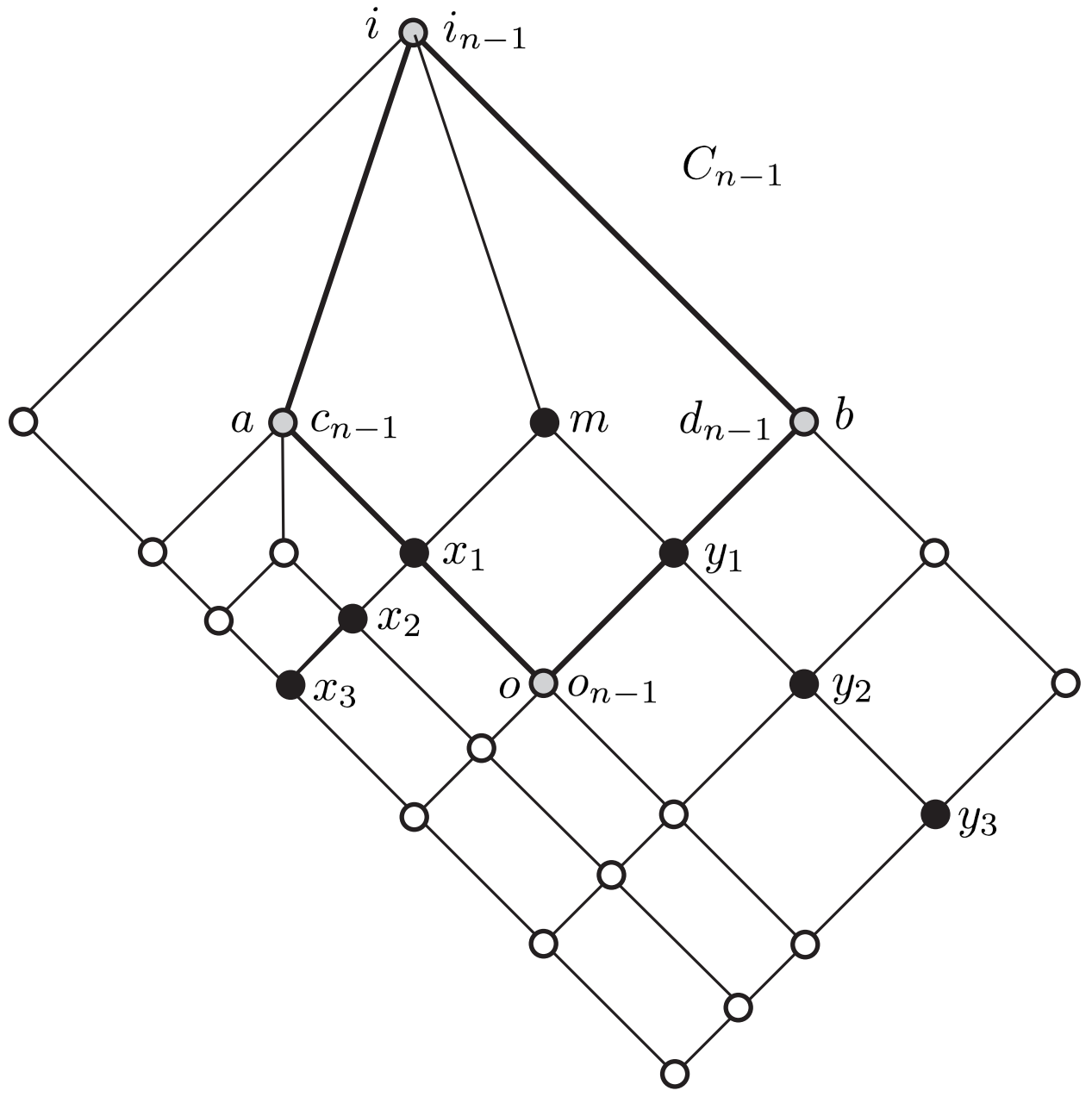

Figure 4.11: Proving Theorem 4.22: Case 2

**Corollary 4.23.** *Let $L$ be an SPS lattice and let $I$ be a rectangular interval of $L$. Let $(P)$ be any property of slim rectangular lattices. Then the property $(P)$ holds for the lattice $I$.*

For instance, let (P) be the property: the intervals $[o, a]$ and $[o, b]$ are chains and all elements of the lower boundary of $I$ are meet-reducible, except for $a, b$. Then we get the main result of G. Czédli [52].

**Corollary 4.24.** *Let $L$ be an SPS lattice and let $I$ be a rectangular interval of $L$ with corners $a, b$. Then $[o, a]$ and $[o, b]$ are chains and all the elements of the lower boundary of $I$ except for $a, b$ are meet-reducible.*

Another nice application is the following.

**Corollary 4.25.** *Let $L$ be an SPS lattice and let $I$ be a rectangular interval of $L$ with corners $a, b$. Then for any $x \in I$, the following equation holds:*

$$x = (x \wedge a) \vee (x \wedge b).$$

There is a more elegant way to formulate the last result.

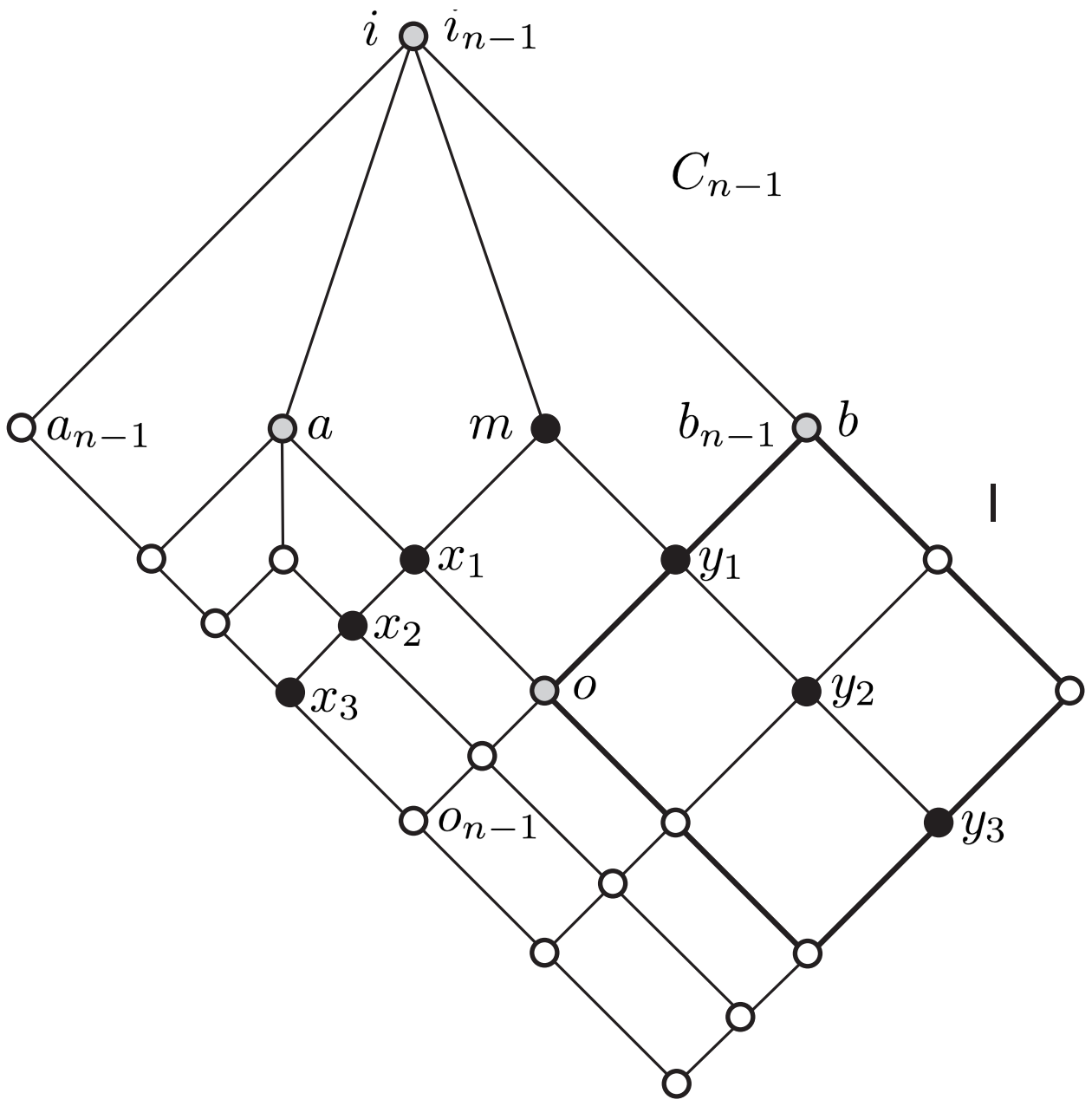

Figure 4.12: Proving Theorem 4.22: Case 3

**Corollary 4.26.** *Let $L$ be an SPS lattice and let $a, b$, and $c$ be pairwise incomparable elements of $L$. If $a$ is to the left of $b$, and $b$ is to the left of $c$, then*

$$b = (b \wedge a) \vee (b \wedge c).$$

## 4.7. Special diagrams for SR lattices

### Two approaches

Let $L$ be an SR lattice. When is a diagram $L$ "nice?" Clearly, if the diagram has some special properties that allow us to easily visualize the properties of $L$ and may even help us with some proofs.

Our "nice" is not the same as the "nice" of our graphic artist friends. They may draw $\mathsf{S}_7$ with the golden ratio, as in the second diagram of Figure 4.13, which may be nicer but not practical. We work with vector graphics apps and grids. The first diagram of Figure 4.13 is easy to draw, but the second diagram is not. Try to draw Figure 29.3 with the golden ratio.

The first attempt to define "nice" was in my joint paper [138] with E. Knapp, published in 2008; we named it *natural.* The third diagram of Figure 4.13 illustrates how we obtain the natural diagram of $\mathsf{S}_7$ by deleting two elements from the diagram of $\mathsf{C}_3^2$.

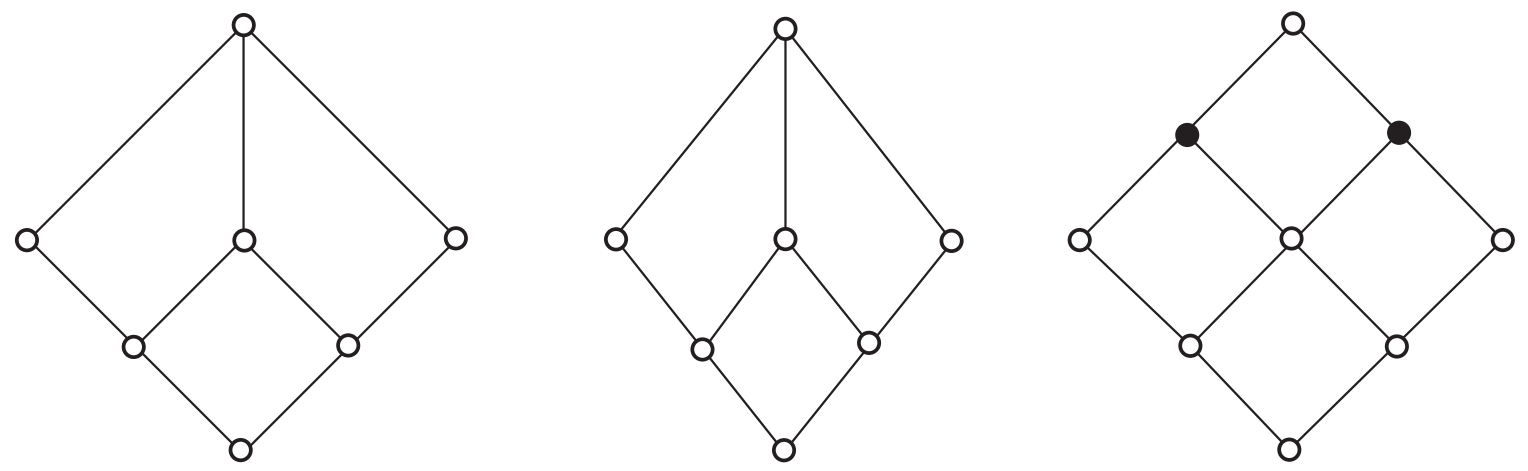

Figure 4.13: 1. $\mathsf{S}_7$; 2. $\mathsf{S}_7$ with golden ratio; 3. a natural diagram of $\mathsf{S}_7$

A decade later (in 2017), the second attempt to define "nice' diagrams was presented in G. Czédli [40].

A $\mathcal{C}_1$-diagram is defined as a diagram with some very nice properties; the angles of the edges are very restrictive. On the other hand, their existence is not trivial.

### Natural diagrams

In this section, we follow my joint paper with E. Knapp [138].

For an SR lattice $L$, let $\mathrm{C_l}(L)$ be the lower left and $\mathrm{C_r}(L)$ the lower right boundary chain of $L$, respectively, and let lc(L) be the left and rc(L) the right corner of $L$, respectively.

We regard the grid $G = \mathrm{C_l}(L) \times \mathrm{C_r}(L)$ as a planar lattice, with $\mathrm{C_l}(L) = \mathrm{C_l}(G)$ and $\mathrm{C_r}(L) = \mathrm{C_r}(G)$. Then the map

$$\psi\colon x \mapsto (x \wedge \mathrm{lc(L)}, x \wedge \mathrm{rc(L)}) \tag{6}$$

is a meet-embedding of $L$ into $G$; the map $\psi$ also preserves the bounds. Therefore, the image of $L$ under $\psi$ in $G$ is a diagram of $L$; we call it the *natural diagram* representing $L$. For instance, the third diagram of Figure 4.13 shows the natural diagram representing $\mathsf{S}_7$.

Let $L$ be an SR lattice. By the Structure Theorem, Strong Version, we can represent $L$ in the form $L = K[C]$, where $K$ is an SR lattice and $C = \{o, c, d, i\}$ is a distributive 4-cell of $K$. Let $\mathcal{D}$ be a diagram of $K$. We form the diagram $\mathcal{D}[C]$ by adding the element $m$ in (the centre of) $C$, the elements $m, x_1, \ldots$, and $m, y_1, \ldots$, as in the last diagram of Figure 4.4, so that the lines spanned by the elements $m, x_1, \ldots$ and m, $y_1, \ldots$ are both normal.

**Lemma 4.27.** *Let $L$, $C$, $K$, $\mathcal{D}$, and $\mathcal{D}[C]$ be as in the previous paragraph. Then $\mathcal{D}[C]$ is a diagram of $L$.*

*Proof.* This is obvious. ☐

**Lemma 4.28.** *Let us make the assumptions of Lemma 4.27. If $\mathcal{D}$ is a natural diagram of $K$, then $\mathcal{D}[C]$ is a natural diagram of $L$.*

*Proof.* Let $\mathcal{D}$ be a natural diagram of $K$. Let the line $m, x_1, \dots$ terminate with $x_{k_l}$ and the line $m, y_1, \dots$ with $y_{k_r}$. We have to show that all the new elements in $L$ can be represented as a join $u_l \vee u_r$, where $u_l \in \mathrm{C_l}(L)$ and $u_r \in \mathrm{C_r}(L)$. Indeed, $m = x_{k_l} \vee x_{k_r}$. The others follow from the distributivity assumptions. □

### $\mathcal{C}_1$-diagrams

This research tool, introduced by G. Czédli, has been playing an important role in some recent papers (see G. Czédli [40]–[52], my joint paper [59] with G. Czédli, and G. Grätzer [124]; for the definition, see G. Czédli [40] and G. Grätzer [124]).

In the diagram of an SR lattice $K$, a *normal edge* (*line*) has a slope of 45° or 135°. Any edge (line) of slope strictly between 45° and 135° is *steep.*

Figure 4.14 depicts the lattice $\mathsf{S}_7$. A *peak sublattice* $\mathsf{S}_7$ (*peak sublattice*, for short) of a lattice $L$ is a sublattice isomorphic to $\mathsf{S}_7$ such that the three edges at the top are covers in the lattice $L$.

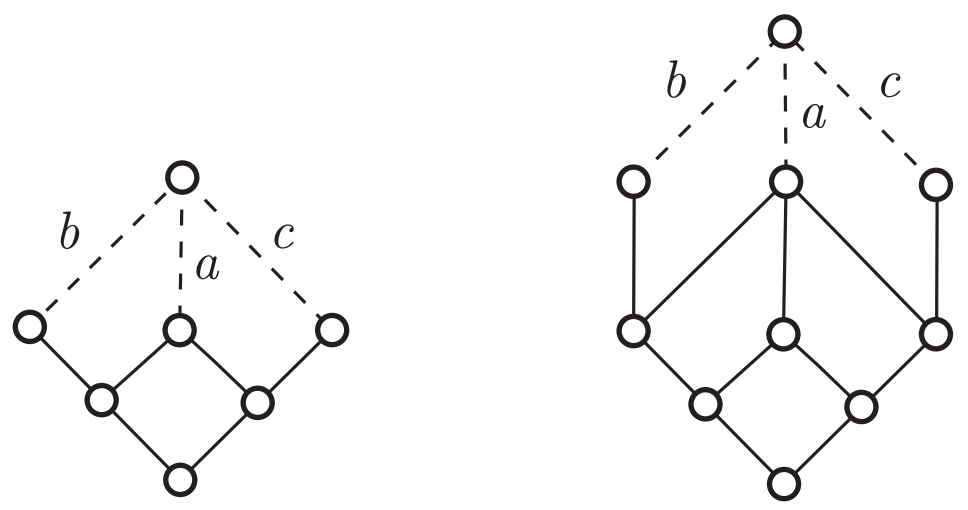

Figure 4.14: Two diagrams of the lattice $\mathsf{S}_7$ as a peak sublattice

**Definition 4.29.** A diagram of an SR lattice $L$ is a $\mathcal{C}_1$-*diagram*, if the middle edge of a peek sublattice is steep and all other edges are normal.

**Theorem 4.30.** *Every slim rectangular lattice $L$ has a $\mathcal{C}_1$-diagram.*

This was proved in G. Czédli [40]. My note [124] presents a short and direct proof. This result also follows from Theorem 4.32.

## 4.8. Natural diagrams and $\mathcal{C}_1$-diagrams

We start with a trivial statement.

**Lemma 4.31.** *Let us make the assumptions of Lemma 4.27. If $\mathcal{D}$ is a $\mathcal{C}_1$-diagram of $K$, then $\mathcal{D}[C]$ is a $\mathcal{C}_1$-diagram of $L$.*

Now we state our second result on SR lattices.

**Theorem 4.32.** *Let $L$ be a SR lattice. Then a natural diagram of $L$ is a $\mathcal{C}_1$-diagram. Conversely, every $\mathcal{C}_1$-diagram is natural.*

*Proof.* Let us assume that the SR lattice $L$ can be obtained from a grid $G$ by inserting forks $n$-times, where $n = \mathrm{Rank(L)}$. We induct on $n$. The case $n = 0$ is trivial because then $L$ is a grid. So let us assume that the theorem holds for $n - 1$.

By the Structure Theorem, Strong Version (Theorem 4.12), there is an SR lattice $K$ and a distributive 4-cell $C = \{o, a, b, i\}$ of $K$ such that $K$ can be obtained from the grid $G$ by inserting forks $(n-1)$-times and also $L = K[C]$ holds.

Now form the natural diagram $\mathcal{D}$ of $K$. By induction, it is a $\mathcal{C}_1$-diagram. By Lemma 4.27, the diagram $\mathcal{D}[C]$ is both natural and $\mathcal{C}_1$.

We prove the converse the same way. ☐

Natural diagrams exist by definition. So as we noted before, Theorem 4.30 also follows from Theorem 4.32.

G. Czédli [40] also defined $\mathcal{C}_2$-*diagrams.* A $\mathcal{C}_1$-diagram is $\mathcal{C}_2$, if any two edges on the lower boundary are of the same length.

We use Theorem 4.32 to prove the following two results of G. Czédli [40].

**Theorem 4.33.** *Let $L$ be a SR lattice. Then $L$ has a $\mathcal{C}_2$-diagram.*

*Proof.* Let $C_l$ and $C_r$ be chains of the same length as $\mathrm{C_l}(L)$ and $\mathrm{C_r}(L)$, respectively. Then $\mathrm{C_l}(L) \times \mathrm{C_r}(L)$ and $C_l \times C_r$ are isomorphic, so we can regard the map $\psi$ (see equation (6)) as a map from $L$ into $C_l \times C_r$, a bounded and meet-preserving map. So the natural diagram it defines is the diagram of the lattice $L$.

If we choose $C_l$ and $C_r$ so that the edges are of the same size, we obtain a $\mathcal{C}_2$-diagram of the SR lattice $L$. ☐

Natural diagrams have a left-right symmetry. The symmetric diagram is obtained with the map

$$\widetilde{\psi} \colon x \mapsto (x \wedge \mathrm{rc(L)}, x \wedge \mathrm{lc(L)}), \tag{7}$$

replacing equation (6).

**Theorem 4.34** (Uniqueness Theorem)**.** *Let $L$ be a slim rectangular lattice. Then the $\mathcal{C}_1$-diagram of $L$ is unique up to left-right symmetry.*

## 4.9. Discussion

Semidistributivity is an important lattice property. Can we prove an RT for finite join-semidistributive lattices?

The three-element chain cannot be represented as the congruence lattice of a finite semidistributive lattice. In fact, consider the following result.

**Theorem 4.35.** *A finite distributive lattice $D$ with a three-element filter cannot be represented as the congruence lattice of a finite semidistributive lattice.*

We start with an easy statement.

**Lemma 4.36.** *Let $a$ be an atom in a finite meet-semidistributive lattice $L$. Define*

$$P_a = \mathrm{id}(a^*) = \{\, x \in L \mid x \leq a^* \,\}.$$

*Then $P_a$ is a prime ideal.*

*Proof.* Indeed, if $P_a$ is not a prime ideal, then there are elements $b, c \in L$ such that $b, c \notin P_a$ but $b \wedge c \in P_a$. From $b \notin P_a$, we get that $a \wedge b \neq 0$. Since $a$ is an atom, we conclude that $a \leq b$. Similarly, $a \leq c$. Therefore, $a \leq b \wedge c$, contradicting that $b \wedge c \in P_a$. □

With a prime ideal $P$, we associate the congruence $\boldsymbol{\varepsilon}_P$ with two congruence classes: $P$ and $L - P$.

By way of contradiction, let the finite meet-semidistributive lattice $L$ represent $D$, a finite distributive lattice with a there-element filter $F$. Let $\boldsymbol{\alpha}$ generate $F$ as a filter. Then $\mathrm{Con}(L/\boldsymbol{\alpha})$ represents $F$. Since the homomorphic image of a finite meet-semidistributive lattice is itself a finite meet-semidistributive lattice (see Corollary 3-1.21 of [8]), we obtain the following: there is a finite meet-semidistributive lattice $K$ representing the three-element chain. Clearly, $K$ has exactly one nontrivial congruence.

If $K$ has exactly one atom $a$, then $K$ has at least two nontrivial congruences, namely $\mathrm{con}(0, a)$ and $\mathrm{con}(a, 1)$, a contradiction.

If $K$ has at least two distinct atoms $a$ and $b$, then $K$ has at least two nontrivial congruences, namely $\boldsymbol{\varepsilon}_{P_a}$ and $\boldsymbol{\varepsilon}_{P_b}$, again a contradiction. This completes the proof of Theorem 4.35.

So congruence lattices of finite semidistributive lattices form a proper subclass of the class of finite distributive lattices.

Theorem 12 in K. Adaricheva, R. Freese, and J. B. Nation [5] characterize the class of finite semidistributive lattices that can be represented as congruence lattices of finite semidistributive lattices.

R. P. Dilworth [78] first studied lattices with unique irreducible decomposition. P. Crawley and R. P. Dilworth [24] present a deeper account. In recent research, *finite convex geometries* are finite lattices with unique irreducible decomposition.

It turns out that a finite lattice is join-distributive iff it has unique irreducible decomposition (see Theorem 3-1.4, K. Adaricheva and J. B. Nation [8]).

There is a lot of research on finite semidistributive lattices (see Chapters 3–6 by K. Adaricheva and J. B. Nation in LTS2-[207] for a survey).

G. Czédli et al. (see [3], [4], [39], [43], [64], and [221]) studied the connection between SPS lattices and finite convex geometries.

**Problem 1.** Can every finite distributive lattice be represented as the congruence lattice of an infinite semidistributive lattice?

R. Freese and J. B. Nation [89] prove that $\mathsf{C}_2$ can be represented as the congruence lattice of an infinite semidistributive lattice.

# Part II

# Some Special Techniques

# Chapter 5

# *Chopped Lattices*

The first basic technique is the use of a chopped lattice, a finite meet-semilattice $(M, \wedge)$ regarded as a *partial algebra* $(M, \wedge, \vee)$, where $\vee$ is a partial operation. It turns out that the ideals of a chopped lattice form a lattice with the same congruence lattice as that of the chopped lattice. So to construct a finite lattice with a given congruence lattice it is enough to construct such a chopped lattice. The problem is how to ensure that the ideal lattice of the chopped lattice has some given properties. As an example, we will look at sectionally complemented lattices.

Chopped lattices were introduced by myself and H. Lakser (published in GLT-[99]). They were named in my joint paper [180] with E. T. Schmidt and generalized to the infinite case in another joint paper with E. T. Schmidt [181].

## 5.1. Basic definitions

An ($n$-ary) *partial operation* on a nonempty set $A$ is a map from a subset of $A^n$ to $A$. For $n = 2$, we call the partial operation *binary*. A *partial algebra* is a nonempty set $A$ with partial operations defined on $A$.

A finite meet-semilattice $(M, \wedge)$ may be regarded as a partial algebra, $(M, \wedge, \vee)$, called a *chopped lattice*, where $\wedge$ is an operation and $\vee$ is a partial operation: $a \vee b$ is the least upper bound of $a$ and $b$, provided that it exists.

We can obtain an example of a chopped lattice by taking a finite lattice with unit, 1, and defining $M = L - \{1\}$. The converse also holds: by adding a new unit 1 to a chopped lattice $M$, we obtain a finite lattice $L$, and chopping off the unit element, we get $M$ back.

A more useful example is obtained with *merging.* Let $C$ and $D$ be lattices such that $J = C \cap D$ is an ideal in both $C$ and $D$. Then with the natural ordering, $\mathrm{Merge}(C, D) = C \cup D$, called the *merging* of $C$ and $D$, is a chopped lattice. Note that whenever $a \vee b = c$ in $\mathrm{Merge}(C, D)$, then we have $a, b, c \in C$ and $a \vee b = c$ in $C$ or $a, b, c \in D$ and $a \vee b = c$ in $D$.

Among finite ordered sets, chopped lattices are easy to spot.

**Lemma 5.1.** *Let $M$ be a finite ordered set and let* Max *be the set of maximal elements of $M$. If $\downarrow m$ is a lattice, for each $m \in \mathrm{Max}$, and if $m \wedge n$ exists, for all $m$, $n \in \mathrm{Max}$, then $M$ is a meet-semilattice.*

*Proof.* Indeed, for $x, y \in M$,

$$x \wedge y = (x \wedge_m a) \wedge_a (y \wedge_n a),$$

where $m, n \in \mathrm{Max}$, $x \leq m$, $y \leq n$, $a = m \wedge n$, and $\wedge_m$, $\wedge_n$, $\wedge_a$ denotes the meet in the lattice $\downarrow m$, $\downarrow n$, and $\downarrow a$, respectively. It is easy to see that $x \wedge y$ is the greatest lower bound of $x$ and $y$ in $M$. □

A meet-subsemilattice $A$ of a chopped lattice $M$ is a *sublattice*, if whenever $c = a \vee b$ in $A$, then we have $c = a \vee b$ in $M$.

We define an equivalence relation $\boldsymbol{\alpha}$ to be a *congruence* of a chopped lattice $M$ as we defined it for lattices in Section 1.3.3: we require that $(\mathrm{SP}_\wedge)$ and $(\mathrm{SP}_\vee)$ hold, the latter with the proviso: whenever $a \vee c$ and $b \vee d$ exist. The set $\mathrm{Con}\, M$ of all congruence relations of $M$ ordered by set inclusion is a lattice.

**Lemma 5.2.** *Let $M$ be a chopped lattice and let $\boldsymbol{\alpha}$ be an equivalence relation on $M$ satisfying the following two conditions for $x, y, z \in M$:*

(1) *If $x \equiv y \pmod{\boldsymbol{\alpha}}$; then $x \wedge z \equiv y \wedge z \pmod{\boldsymbol{\alpha}}$.*

(2) *If $x \equiv y \pmod{\boldsymbol{\alpha}}$ and $x \vee z$ and $y \vee z$ exist, then $x \vee z \equiv y \vee z \pmod{\boldsymbol{\alpha}}$.*

*Then $\boldsymbol{\alpha}$ is a congruence relation on $M$.*

*Proof.* Condition (1) states that $\boldsymbol{\alpha}$ preserves $\wedge$.

Now let $x, y, u, v \in S$ with $x \equiv y \pmod{\boldsymbol{\alpha}}$ and $u \equiv v \pmod{\boldsymbol{\alpha}}$; let $x \vee u$ and $y \vee v$ exist. Then $x \equiv x \wedge y \equiv y \pmod{\boldsymbol{\alpha}}$ and $(x \wedge y) \vee u$ and $(x \wedge y) \vee v$ exist. Thus, by condition (2),

$$x \vee u \equiv (x \wedge y) \vee u \equiv (x \wedge y) \vee v \equiv y \vee v \pmod{\boldsymbol{\alpha}}.$$ □

A nonempty subset $I$ of the chopped lattice $M$ is an *ideal* iff it is a down set with the property,

(Id) $a, b \in I$ implies that $a \vee b \in I$, provided that $a \vee b$ exists in $M$.

The set $\operatorname{Id} M$ of all ideals of $M$ ordered by set inclusion is a lattice. For $I, J \in \operatorname{Id} M$, the meet is $I \cap J$, but the join is a bit more complicated to describe.

**Lemma 5.3.** *Let $I$ and $J$ be ideals of the chopped lattice $M$. Define*

$$U(I,J)_0 = I \cup J,$$
$$U(I,J)_i = \{\, x \mid x \leq u \vee v,\ u,\ v \in U(I,J)_{i-1} \,\} \text{ for } 0 < i < \omega.$$

*Then*

$$I \vee J = \bigcup(\, U(I,J)_i \mid i < \omega\,).$$

*Proof.* Define $U = \bigcup(\, U(I,J)_i \mid i < \omega\,)$. If $K$ is an ideal of $M$, then $I \subseteq K$ and $J \subseteq K$ imply—by induction—that $U \subseteq K$. So it is sufficient to prove that $U$ is an ideal of $M$.

Obviously, $U$ is a down set. Also, $U$ has property (Id), since if $a, b \in U$ and $a \vee b$ exists in $M$, then $a, b \in U(I,J)_n$, for some $0 < n < \omega$, and therefor $a \vee b \in U(I,J)_{n+1} \subseteq U$. ☐

Most lattice concepts and notations for them will be used for chopped lattices without further explanation.

Infinite chopped lattices are briefly mentioned on page 113.

## 5.2. Compatible vectors of elements

Let $M$ be a chopped lattice, and let $\operatorname{Max}(M)$ (Max if $M$ is understood) be the set of maximal elements of $M$. Then $M = \bigcup(\, \operatorname{id}(m) \mid m \in \operatorname{Max}\,)$ and each $\operatorname{id}(m)$ is a (finite) lattice. A *vector* (associated with $M$) is of the form $(i_m \mid m \in \operatorname{Max})$, where $i_m \leq m$ for all $m \in M$. We order the vectors componentwise.

With every ideal $I$ of $M$, we can associate the vector $(i_m \mid m \in \operatorname{Max})$ defined by $I \cap \operatorname{id}(m) = \operatorname{id}(i_m)$. Clearly, $I = \bigcup(\, \operatorname{id}(i_m) \mid m \in M\,)$. Such vectors are easy to characterize. Let us call the vector $(j_m \mid m \in \operatorname{Max})$ *compatible* if $j_m \wedge n = j_n \wedge m$ for all $m, n \in \operatorname{Max}$.

**Lemma 5.4.** *Let $M$ be a chopped lattice.*

(i) *There is a one-to-one correspondence between ideals and compatible vectors of $M$.*

(ii) *Given any vector $\mathbf{g} = (g_m \mid m \in \operatorname{Max})$, there is a smallest compatible vector $\overline{\mathbf{g}} = (i_m \mid m \in \operatorname{Max})$ containing $\mathbf{g}$.*

(iii) *Let $I$ and $J$ be ideals of $M$, with corresponding compatible vectors*

$$(i_m \mid m \in \operatorname{Max}) \text{ and } (j_m \mid m \in \operatorname{Max}).$$

*Then*

(a) *$I \leq J$ in* $\operatorname{Id} M$ *iff* $i_m \leq j_m$ *for all* $m \in \operatorname{Max}$.

(b) *The compatible vector corresponding to* $I \wedge J$ *is* $(i_m \wedge j_m \mid m \in \operatorname{Max})$.

(c) *Let* $\mathbf{a} = (i_m \vee j_m \mid m \in \operatorname{Max})$. *Then the compatible vector corresponding to* $I \vee J$ *is* $\overline{\mathbf{a}}$.

*Proof.*

(i) Let $I$ be an ideal of $M$. Then $(i_m \mid m \in \operatorname{Max})$ is compatible since $i_m \wedge m \wedge n$ and $i_n \wedge m \wedge n$ both generate the principal ideal $I \cap \operatorname{id}(m) \cap \operatorname{id}(n)$ for all $m, n \in \operatorname{Max}$.

Conversely, let $(j_m \mid m \in \operatorname{Max})$ be compatible, and define

$$I = \bigcup( \operatorname{id}(j_m) \mid m \in \operatorname{Max} ).$$

Observe that

$$I \cap \operatorname{id}(m) = \operatorname{id}(j_m), \tag{1}$$

for $m \in \operatorname{Max}$. Indeed, if $x \in I \cap \operatorname{id}(m)$ and $x \in \operatorname{id}(j_n)$, for $n \in \operatorname{Max}$, then $x \leq m \wedge j_n = n \wedge j_m$ (since $(j_m \mid m \in \operatorname{Max})$ is compatible), so $x \leq j_m$, that is, $x \in \operatorname{id}(j_m)$. The reverse inclusion is obvious.

$I$ is obviously a down set. To verify property (Id) for $I$, let $a, b \in I$ and let us assume that $a \vee b$ exists in $M$. Then $a \vee b \leq m$, for some $m \in \operatorname{Max}$, so $a \leq m$ and $b \leq m$. By (1), we get $a \leq j_m$ and $b \leq j_m$, therefore, $a \vee b \leq j_m \in I$. Since $I$ is a down set, it follows that $a \vee b \in I$, verifying property (Id).

(ii) Obviously, the vector $(m \mid m \in \operatorname{Max})$ contains all other vectors and it is compatible. Since the componentwise meet of compatible vectors is compatible, the statement follows.

(iii) This is obvious since, by (ii), we are dealing with a closure system (see Section 1.2.2). ☐

## 5.3. Compatible vectors of congruences

Let $M$ be a chopped lattice. With any congruence $\boldsymbol{\alpha}$ of $M$, we can associate the *restriction vector* $(\boldsymbol{\alpha}\rceil_m \mid m \in \operatorname{Max})$, where $\boldsymbol{\alpha}\rceil_m$ is the restriction of $\boldsymbol{\alpha}$ to $\operatorname{id}(m)$. The restriction $\boldsymbol{\alpha}\rceil_m$ is a congruence of the lattice $\operatorname{id}(m)$.

Let $\boldsymbol{\gamma}_m$ be a congruence of $\operatorname{id}(m)$ for all $m \in \operatorname{Max}$. The *congruence vector* $(\boldsymbol{\gamma}_m \mid m \in \operatorname{Max})$ is called *compatible* if $\boldsymbol{\gamma}_m$ restricted to $\operatorname{id}(m \wedge n)$ is the same as $\boldsymbol{\gamma}_n$ restricted to $\operatorname{id}(m \wedge n)$ for $m, n \in \operatorname{Max}$. Obviously, a restriction vector is compatible. The converse also holds.

**Lemma 5.5.** *Let* $(\boldsymbol{\gamma}_m \mid m \in \operatorname{Max})$ *be a compatible congruence vector of a chopped lattice* $M$. *Then there is a unique congruence* $\boldsymbol{\alpha}$ *of* $M$ *such that the restriction vector of* $\boldsymbol{\alpha}$ *agrees with* $(\boldsymbol{\gamma}_m \mid m \in \operatorname{Max})$.

*Proof.* Let $(\boldsymbol{\gamma}_m \mid m \in \mathrm{Max})$ be a compatible congruence vector. We define a binary relation $\boldsymbol{\alpha}$ on $M$ as follows.

Let $m, n \in \mathrm{Max}$. For $x \in \mathrm{id}(m)$ and $y \in \mathrm{id}(n)$, let $x \equiv y \pmod{\boldsymbol{\alpha}}$ iff $x \equiv x \wedge y \pmod{\boldsymbol{\gamma}_m}$ and $y \equiv x \wedge y \pmod{\boldsymbol{\gamma}_n}$.

Obviously, $\boldsymbol{\alpha}$ is reflexive and symmetric. To prove transitivity, let $m, n, k \in \mathrm{Max}$, and let $x \in \mathrm{id}(m)$, $y \in \mathrm{id}(n)$, $z \in \mathrm{id}(k)$; let $x \equiv y \pmod{\boldsymbol{\alpha}}$ and $y \equiv z \pmod{\boldsymbol{\alpha}}$, that is,

$$x \equiv x \wedge y \pmod{\boldsymbol{\gamma}_m}, \tag{2}$$

$$y \equiv x \wedge y \pmod{\boldsymbol{\gamma}_n}, \tag{3}$$

$$y \equiv y \wedge z \pmod{\boldsymbol{\gamma}_n}, \tag{4}$$

$$z \equiv y \wedge z \pmod{\boldsymbol{\gamma}_k}. \tag{5}$$

Then meeting the congruence (4) with $x$ (in the lattice $\mathrm{id}(n)$), we get

$$x \wedge y \equiv x \wedge y \wedge z \pmod{\boldsymbol{\gamma}_n}, \tag{6}$$

and from (3), by meeting with $z$, we obtain

$$y \wedge z \equiv x \wedge y \wedge z \pmod{\boldsymbol{\gamma}_n}. \tag{7}$$

Since $x \wedge y$ and $x \wedge y \wedge z \in \mathrm{id}(m)$, by compatibility, (6) implies that

$$x \wedge y \equiv x \wedge y \wedge z \pmod{\boldsymbol{\gamma}_m}. \tag{8}$$

Now (2) and (8) yield

$$x \equiv x \wedge y \wedge z \pmod{\boldsymbol{\gamma}_m}. \tag{9}$$

Similarly,

$$z \equiv x \wedge y \wedge z \pmod{\boldsymbol{\gamma}_k}. \tag{10}$$

$\boldsymbol{\gamma}_m$ is a lattice congruence on $\mathrm{id}(m)$ and $x \wedge y \wedge z \leq x \wedge z \leq x$, so

$$x \equiv x \wedge z \pmod{\boldsymbol{\gamma}_m}. \tag{11}$$

Similarly,

$$z \equiv x \wedge z \pmod{\boldsymbol{\gamma}_k}. \tag{12}$$

Equations (11) and (12) yield that $x \equiv z \pmod{\boldsymbol{\alpha}}$, proving transitivity.

$(\mathrm{SP}_\wedge)$ is easy: let $x \in \mathrm{id}(m)$, $y \in \mathrm{id}(n)$, $z \in M$; if $x \equiv y \pmod{\boldsymbol{\alpha}}$, then $x \wedge z \equiv y \wedge z \pmod{\boldsymbol{\alpha}}$ because $x \wedge z \equiv x \wedge y \wedge z \pmod{\boldsymbol{\gamma}_m}$ and $y \wedge z \equiv x \wedge y \wedge z \pmod{\boldsymbol{\gamma}_n}$.

Finally, we verify $(\mathrm{SP}_\vee)$. Let $x \equiv y \pmod{\boldsymbol{\alpha}}$ and $z \in M$, and let us assume that $x \vee z$ and $y \vee z$ exist. Then there are $p, q \in \mathrm{Max}$ such that $x \vee z \in \mathrm{id}(p)$ and $y \vee z \in \mathrm{id}(q)$. By compatibility, $x \equiv x \wedge y \pmod{\boldsymbol{\gamma}_p}$, so $x \vee z \equiv (x \wedge y) \vee z \pmod{\boldsymbol{\gamma}_p}$. Since $(x \wedge y) \vee z \leq (x \vee z) \wedge (y \vee z) \leq x \vee z$, we also have

$$x \vee z \equiv (x \vee z) \wedge (y \vee z) \pmod{\boldsymbol{\gamma}_p}.$$

Similarly,

$$y \vee z \equiv (x \vee z) \wedge (y \vee z) \pmod{\boldsymbol{\gamma}_q}.$$

The last two displayed equations show that $x \vee z \equiv y \vee z \pmod{\boldsymbol{\alpha}}$. □

## 5.4. From the chopped lattice to the ideal lattice

The map $m \mapsto \mathrm{id}(m)$ embeds the chopped lattice $M$ with zero into the lattice $\mathrm{Id}\, M$, so we can regard $\mathrm{Id}\, M$ as an extension. It is, in fact, a congruence-preserving extension (see my joint paper with H. Lakser [139], proof first published in LTFC-[98]).

**Theorem 5.6.** *Let $M$ be a chopped lattice. Then $\mathrm{Id}\, M$ is a congruence-preserving extension of $M$.*

*Proof.* Let $\boldsymbol{\alpha}$ be a congruence relation of $M$. If $I, J \in \mathrm{Id}\, M$, define

$$I \equiv J \pmod{\overline{\boldsymbol{\alpha}}} \quad \text{iff} \quad I/\boldsymbol{\alpha} = J/\boldsymbol{\alpha},$$

where $I/\boldsymbol{\alpha}$ denotes $I$ modulo $\boldsymbol{\alpha}]I$ and the same for $J/\boldsymbol{\alpha}$. Obviously, $\overline{\boldsymbol{\alpha}}$ is an equivalence relation. Let $I \equiv J \pmod{\overline{\boldsymbol{\alpha}}}$, $N \in \mathrm{Id}\, M$, and $x \in I \cap N$. Then $x \equiv y \pmod{\boldsymbol{\alpha}}$, for some $y \in J$, and so $x \equiv x \wedge y \pmod{\boldsymbol{\alpha}}$ and $x \wedge y \in J \cap N$. This shows that $(I \cap N)/\boldsymbol{\alpha} \subseteq (J \cap N)/\boldsymbol{\alpha}$. Similarly, $(J \cap N)/\boldsymbol{\alpha} \subseteq (I \cap N)/\boldsymbol{\alpha}$, so $I \cap N \equiv J \cap N \pmod{\overline{\boldsymbol{\alpha}}}$.

To prove $I \vee N \equiv J \vee N \pmod{\overline{\boldsymbol{\alpha}}}$, by symmetry, it is sufficient to verify that $I \vee N \subseteq (J \vee N)/\boldsymbol{\alpha}$. By Lemma 5.3, this is equivalent to proving that $U_n \subseteq (J \vee N)/\boldsymbol{\alpha}$ for $n < \omega$. This is obvious for $n = 0$.

Now assume that $U_{n-1} \subseteq (J \vee N)/\boldsymbol{\alpha}$ and let $x \in U_n$. Then $x \leq t_1 \vee t_2$ for some $t_1, t_2 \in U_{n-1}$. Thus $t_1 \equiv u_1 \pmod{\boldsymbol{\alpha}}$ and $t_2 \equiv u_2 \pmod{\boldsymbol{\alpha}}$, for some $u_1, u_2 \in J \vee N$, and so $t_1 \equiv t_1 \wedge u_1 \pmod{\boldsymbol{\alpha}}$ and $t_2 \equiv t_2 \wedge u_2 \pmod{\boldsymbol{\alpha}}$. Observe that $t_1 \vee t_2$ is an upper bound for $\{t_1 \wedge u_1, t_2 \wedge u_2\}$; consequently, $(t_1 \wedge u_1) \vee (t_2 \wedge u_2)$ exists. Therefore,

$$t_1 \vee t_2 \equiv (t_1 \wedge u_1) \vee (t_2 \wedge u_2) \pmod{\boldsymbol{\alpha}}.$$

Finally,

$$x \equiv x \wedge (t_1 \vee t_2) = x \wedge ((t_1 \wedge u_1) \vee (t_2 \wedge u_2)) \pmod{\boldsymbol{\alpha}},$$

and

$$x \wedge ((t_1 \wedge u_1) \vee (t_2 \wedge u_2)) \in J \vee N.$$

Thus $x \in (J \vee N)/\boldsymbol{\alpha}$, completing the induction, verifying that $\overline{\boldsymbol{\alpha}}$ is a congruence relation of $\operatorname{Id} M$.

If $a \equiv b \pmod{\boldsymbol{\alpha}}$ and $x \in \operatorname{id}(a)$, then $x \equiv x \wedge b \pmod{\boldsymbol{\alpha}}$. Thus $\operatorname{id}(a) \subseteq \operatorname{id}(b)/\boldsymbol{\alpha}$. Similarly, $\operatorname{id}(b) \subseteq \operatorname{id}(a)/\boldsymbol{\alpha}$, and so $\operatorname{id}(a) \equiv \operatorname{id}(b) \pmod{\overline{\boldsymbol{\alpha}}}$. Conversely, if $\operatorname{id}(a) \equiv \operatorname{id}(b) \pmod{\overline{\boldsymbol{\alpha}}}$, then $a \equiv b_1 \pmod{\boldsymbol{\alpha}}$ and $a_1 \equiv b \pmod{\boldsymbol{\alpha}}$ for some $a_1 \leq a$ and $b_1 \leq b$. Forming the join of these two congruences, we get $a \equiv b \pmod{\overline{\boldsymbol{\alpha}}}$. Thus $\overline{\boldsymbol{\alpha}}$ has all the properties required by Lemma 5.3.

To show the uniqueness, let $\boldsymbol{\gamma}$ be a congruence relation of $\operatorname{Id} M$ satisfying $\operatorname{id}(a) \equiv \operatorname{id}(b) \pmod{\boldsymbol{\gamma}}$ iff $a \equiv b \pmod{\boldsymbol{\alpha}}$. Let $I, J \in \operatorname{Id} M$, $I \equiv J \pmod{\boldsymbol{\gamma}}$, and $x \in I$. Then

$$\begin{aligned} \operatorname{id}(x) \cap I &\equiv \operatorname{id}(x) \cap J \pmod{\boldsymbol{\gamma}},\\ \operatorname{id}(x) \cap I &= \operatorname{id}(x),\\ \operatorname{id}(x) \cap J &= \operatorname{id}(y) \end{aligned}$$

for some $y \in J$. Thus $\operatorname{id}(x) \equiv \operatorname{id}(y) \pmod{\boldsymbol{\gamma}}$, and so $x \equiv y \pmod{\boldsymbol{\alpha}}$, proving that $I \subseteq J/\boldsymbol{\alpha}$. Similarly, $J \subseteq I/\boldsymbol{\alpha}$, and so $I \equiv J \pmod{\overline{\boldsymbol{\alpha}}}$. Conversely, if $I \equiv J \pmod{\overline{\boldsymbol{\alpha}}}$, then take all congruences of the form $x \equiv y \pmod{\boldsymbol{\alpha}}$, $x \in I$, $y \in J$. By our assumption regarding $\boldsymbol{\gamma}$, we get the congruence $\operatorname{id}(x) \equiv \operatorname{id}(y) \pmod{\boldsymbol{\gamma}}$, and by our definition of $\overline{\boldsymbol{\alpha}}$, the join of all these congruences yields $I \equiv J \pmod{\boldsymbol{\alpha}}$. Thus $\boldsymbol{\gamma} = \overline{\boldsymbol{\alpha}}$. □

This result is very useful. It means that in order to construct a finite lattice $L$ to represent a given finite distributive lattice $D$ as a congruence lattice, it is sufficient to construct a chopped lattice $M$ with $\operatorname{Con} M \cong D$, since $\operatorname{Con} M \cong \operatorname{Con}(\operatorname{Id} M) = \operatorname{Con} L$, where $L = \operatorname{Id} M$, and $L$ is a finite lattice.

This result also allows us to construct congruence-preserving extensions.

**Corollary 5.7.** *Let $M = \operatorname{Merge}(A, B)$ be a chopped lattice with $A = \operatorname{id}(a)$ and $B = \operatorname{id}(b)$. If $a \wedge b > 0$ and $B$ is simple, then $\operatorname{Id} M$ is a congruence-preserving extension of $A$.*

*Proof.* Let $\boldsymbol{\alpha}$ be a congruence of $A$. Then $(\boldsymbol{\alpha}, \boldsymbol{\gamma})$ is a compatible congruence vector iff

$$\boldsymbol{\gamma} = \begin{cases} \mathbf{0} & \text{if } \boldsymbol{\alpha} \text{ is discrete on } [0, a \wedge b];\\ \mathbf{1}, & \text{otherwise.} \end{cases}$$

So $\boldsymbol{\gamma}$ is determined by $\boldsymbol{\alpha}$ and the statement follows. □

## 5.5. Sectional complementation

We introduce *sectionally complemented chopped lattices* as we did for lattices in Section 2.1.

We illustrate the use of compatible vectors with two results on sectionally complemented chopped lattices. The first result is from my joint paper [184] with E. T. Schmidt.

**Lemma 5.8** (Atom Lemma). *Let $M$ be a chopped lattice with two maximal elements $m_1$ and $m_2$. We assume that* $\mathrm{id}(m_1)$ *and* $\mathrm{id}(m_2)$ *are sectionally complemented lattices. If $p = m_1 \wedge m_2$ is an atom, then* $\mathrm{Id}\, M$ *is sectionally complemented.*

*Proof.* Figure 5.1 illustrates the setup.

To show that $\mathrm{Id}\, M$ is sectionally complemented, let $I \subseteq J$ be two ideals of $M$, represented by the compatible vectors $(i_1, i_2)$ and $(j_1, j_2)$, respectively. Let $s_1$ be the sectional complement of $i_1$ in $j_1$ and let $s_2$ be the sectional complement of $i_2$ in $j_2$. If $p \wedge s_1 = p \wedge s_2$, then $(s_1, s_2)$ is a compatible vector, representing an ideal $S$ that is a sectional complement of $I$ in $J$. Otherwise, without loss of generality, we can assume that $p \wedge s_1 = 0$ and $p \wedge s_2 = p$. Since $\mathrm{id}(m_2)$ is sectionally complemented, there is a sectional complement $s_2'$ of $p$ in $[0, s_2]$. Then $(s_1, s_2')$ satisfies $p \wedge s_1 = p \wedge s_2'$ $(= 0)$, and so it is compatible; therefore, $(s_1, s_2')$ represents an ideal $S$ of $M$. Obviously, $I \wedge S = \{0\}$.

From $p \wedge s_2 = p$, it follows that $p \leq s_2 \leq j_2$. Since $J$ is an ideal and $j_2 \wedge p = p$, it follows that $j_1 \wedge p = p$, that is, $p \leq j_1$. Obviously, $I \vee S \subseteq J$. So to show that $I \vee S = J$, it is sufficient to verify that $j_1, j_2 \in I \vee S$. Evidently, $j_1 = i_1 \vee s_1 \in I \vee S$. Note that $p \leq j_1 = i_1 \vee s_1 \in I \vee S$. Thus $p, s_2', i_2 \in I \vee S$, and therefore

$$p \vee s_2' \vee i_2 = (p \vee s_2') \vee i_2 = s_2 \vee i_2 = j_2 \in I \vee S. \qquad \square$$

The second result (G. Grätzer, H. Lakser, and M. Roddy [155]) shows that the ideal lattice of a sectionally complemented chopped lattice is not always sectionally complemented.

**Theorem 5.9.** *There is a sectionally complemented chopped lattice $M$ whose ideal lattice* $\mathrm{Id}\, M$ *is not sectionally complemented.*

*Proof.* Let $M$ be the chopped lattice of Figure 5.2, where $L_1 = \mathrm{id}(m_1)$ and $L_2 = \mathrm{id}(m_2)$. Note that $p$ is meet-irreducible in $\mathrm{id}(m_2)$ and $q$ is meet-irreducible in $\mathrm{id}(m_1)$.

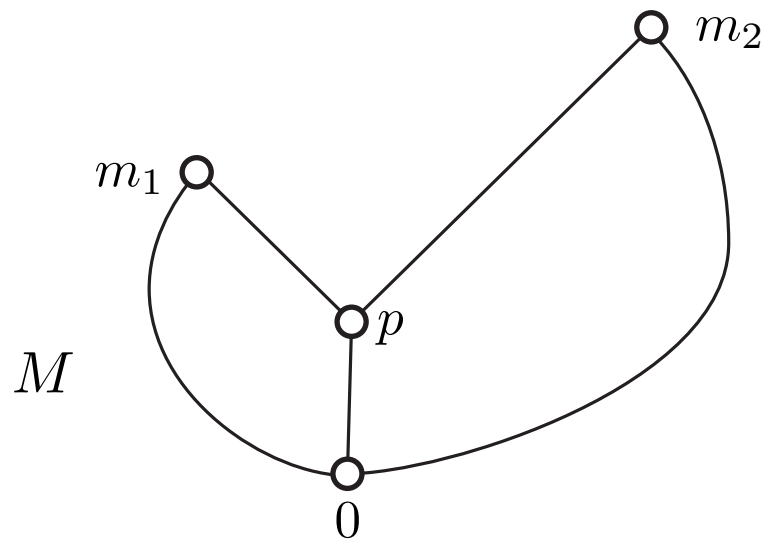

Figure 5.1: The Atom Lemma illustrated

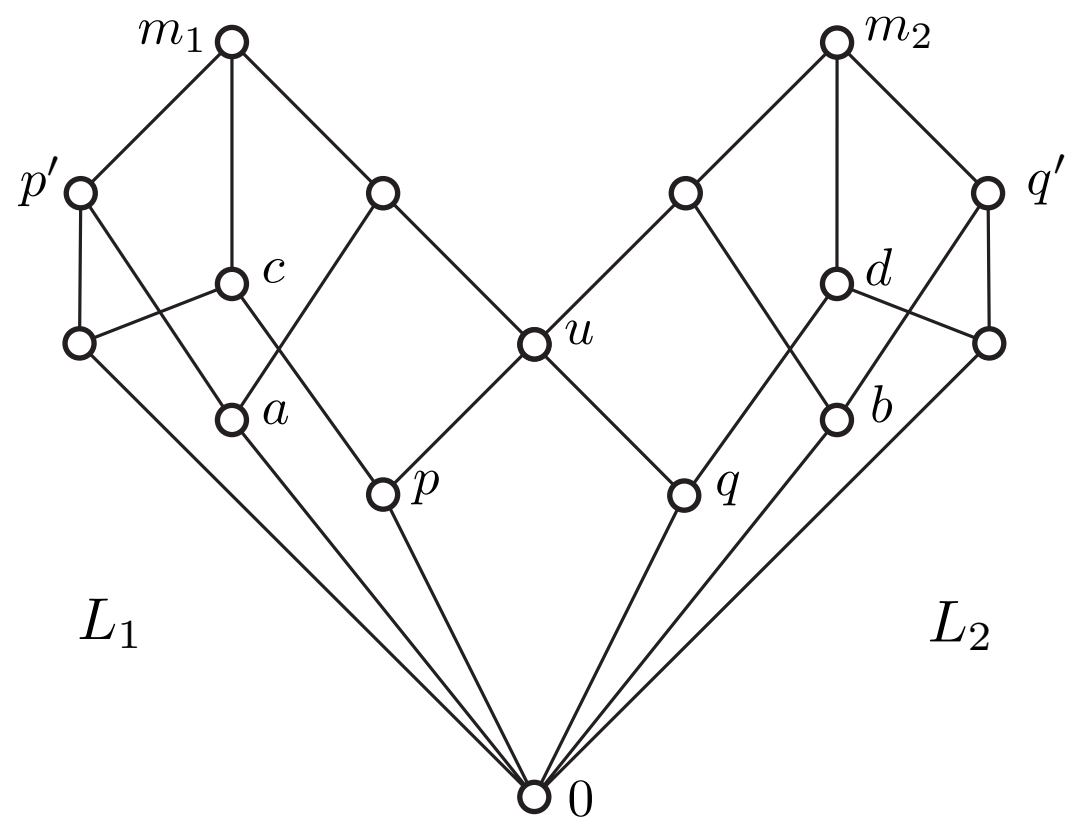

Figure 5.2: The chopped lattice $M$

The unit element of the ideal lattice of $M$ is the compatible vector $(m_1, m_2)$. We show that the compatible vector $(a, b)$ has no complement in the ideal lattice of $M$.

Assume, to the contrary, that the compatible vector $(s, t)$ is a complement of $(a, b)$. Since $(a, b) \leq (a \vee u, m_2)$, a compatible vector, $(s, t) \nleq (a \vee u, m_2)$, that is,

$$s \nleq a \vee u. \tag{13}$$

Similarly, by considering $(m_1, b \vee u)$, we conclude that

$$t \nleq b \vee u. \tag{14}$$

Now $(a, b) \leq (p', q')$, which is a compatible vector. Thus $(s, t) \nleq (p', q')$, and so either $s \nleq p'$ or $t \nleq q'$. Without loss of generality, we may assume that $s \nleq p'$. It then follows by (13) that $s$ can be only $c$ or $m_1$. Then since $s \wedge a = 0$, we conclude that $s = c$. Thus $s \wedge u = p$, and so $t \wedge u = p$. But $p$ is meet-irreducible in $L_2$. Thus $t = p \leq b \vee u$, contradicting (14). $\square$

This result illustrates that the Atom Lemma (Lemma 5.8) cannot be extended to the case where $[0, m_1 \wedge m_2]$ is a four-element Boolean lattice.

# Chapter 6

# *Boolean Triples*

In Part IV, we construct congruence-preserving extensions of finite lattices, extensions with special properties, such as sectionally complemented, semimodular, and so on.

It is easy to construct a proper congruence-preserving extension of a nontrivial finite lattice. In the early 1990s, in my joint paper [180] with E. T. Schmidt, we raised the question whether *every* nontrivial lattice has a proper congruence-preserving extension. (See another of my joint papers with E. T. Schmidt [181].)

It took almost a decade for the answer to appear in my joint paper [199] with F. Wehrung. For infinite lattices, the affirmative answer was provided by the Boolean triples construction, which is described in this chapter. It is interesting that Boolean triples also provide an important tool for finite lattices.

## 6.1. The general construction

In this section, I describe a congruence-preserving extension of a (finite) lattice $L$, introduced in my joint paper with F. Wehrung [199]. We will see that this generalizes a construction of E. T. Schmidt [238] for bounded distributive lattices.

For a lattice $L$, let us call the triple $(x, y, z) \in L^3$ *Boolean* iff

$$\text{(F)} \qquad \begin{aligned} x &= (x \vee y) \wedge (x \vee z),\\ y &= (y \vee x) \wedge (y \vee z),\\ z &= (z \vee x) \wedge (z \vee y), \end{aligned}$$

where F stands for "Fixed point definition."

Note that by Lemma 2.6, if (F) holds, then $\mathrm{sub}(\{x, y, z\})$ is Boolean.

(F) is a "Fixed point definition" because the triple $(x, y, z)$ satisfies (F) iff $p(x, y, z) = (x, y, z)$, where

$$p(x, y, z) = ((x \vee y) \wedge (x \vee z), (y \vee x) \wedge (y \vee z), (z \vee x) \wedge (z \vee y)).$$

We denote by $\mathsf{M}_3[L] \subseteq L^3$ the ordered set of Boolean triples of $L$ (ordered as an ordered subset of $L^3$, that is, componentwise). If we apply the construction to an interval $[a, b]$ of $L$, we write $\mathsf{M}_3[a, b]$ for $\mathsf{M}_3[[a, b]]$.

Observe that any Boolean triple $(x, y, z) \in L^3$ satisfies

$$\text{(B)} \qquad x \wedge y = y \wedge z = z \wedge x.$$

where B stands for "balanced." Indeed, if $(x, y, z)$ is Boolean, then

$$x \wedge y = y \wedge z = z \wedge x = (x \vee y) \wedge (y \vee z) \wedge (z \vee x).$$

We call such triples *balanced*.

Here are some of the basic properties of Boolean triples.

**Lemma 6.1.** *Let $L$ be a lattice.*

(i) *$(x, y, z) \in L^3$ is Boolean iff there is a triple $(u, v, w) \in L^3$ satisfying*

$$\text{(E)} \qquad \begin{aligned} x &= u \wedge v, \\ y &= u \wedge w, \\ z &= v \wedge w. \end{aligned}$$

(ii) *$\mathsf{M}_3[L]$ is a closure system in $L^3$. For $(x, y, z) \in L^3$, the closure is*

$$\overline{(x, y, z)} = ((x \vee y) \wedge (x \vee z), (y \vee x) \wedge (y \vee z), (z \vee x) \wedge (z \vee y)).$$

(iii) *If $L$ has $0$, then the ordered subset $\{\, (x, 0, 0) \mid x \in L \,\}$ is a sublattice of $\mathsf{M}_3[L]$ and $\gamma \colon x \mapsto (x, 0, 0)$ is an isomorphism between $L$ and this sublattice.*

(iv) *If $L$ is bounded, then $\mathsf{M}_3[L]$ has a spanning $\mathsf{M}_3$, that is, a bounded sublattice isomorphic to $\mathsf{M}_3$, namely,*

$$\{(0, 0, 0), (1, 0, 0), (0, 1, 0), (0, 0, 1), (1, 1, 1)\}.$$

*Remark.* (E) is the "Existential definition" because $(x, y, z)$ is Boolean iff there *exists* a triple $(u, v, w)$ satisfying (E).

*Proof.*

(i) If $(x, y, z)$ is Boolean, then $u = x \vee y$, $v = x \vee z$, and $w = y \vee z$ satisfy (E). Conversely, if there is a triple $(u, v, w) \in L^3$ satisfying (E), then by Lemma 2.6, the sublattice generated by $x, y, z$ is isomorphic to a quotient of $\mathsf{B}_3$ and $x, y, z$ are the images of the three atoms of $\mathsf{B}_3$. Thus $(x \vee y) \wedge (x \vee z) = x$, the first equation in (F). The other two equations are proved similarly.

(ii) $\mathsf{M}_3[L] \neq \varnothing$; for instance, for all $x \in L$, the diagonal element $(x, x, x) \in \mathsf{M}_3[L]$.

For $(x, y, z) \in L^3$, define $u = x \vee y$, $v = x \vee z$, and $w = y \vee z$. Set $x_1 = u \wedge v$, $y_1 = u \wedge w$, and $z_1 = v \wedge w$. Then $(x_1, y_1, z_1)$ is Boolean by (i) and $(x, y, z) \leq (x_1, y_1, z_1)$ in $L^3$. Now if $(x, y, z) \leq (x_2, y_2, z_2)$ in $L^3$ and $(x_2, y_2, z_2)$ is Boolean, then

$$\begin{aligned} x_2 &= (x_2 \vee y_2) \wedge (x_2 \vee z_2) && \text{(by (F))} \\ &\geq (x \vee y) \wedge (x \vee z) && \text{(by } (x, y, z) \leq (x_2, y_2, z_2)\text{)} \\ &= u \wedge v = x_1, \end{aligned}$$

and similarly, $y_2 \geq y_1$, $z_2 \geq z_1$. Thus $(x_2, y_2, z_2) \geq (x_1, y_1, z_1)$, and so $(x_1, y_1, z_1)$ is the smallest Boolean triple containing $(x, y, z)$.

(iii) and (iv) are obvious. ☐

$\mathsf{M}_3[L]$ is difficult to draw, in general. Figure 6.1 shows the diagram of $\mathsf{M}_3[\mathsf{C}_3]$ with the three-element chain $\mathsf{C}_3 = \{0, a, 1\}$.

If $C$ is an arbitrary bounded chain, with bounds 0 and 1, it is easy to picture $\mathsf{M}_3[C]$, as sketched in Figure 6.1. The element $(x, y, z) \in C^3$ is Boolean iff it is of the form $(x, y, y)$, or $(y, x, y)$, or $(y, y, x)$, where $y \leq x$ in $C$. So the diagram is made up of three isomorphic "flaps" overlapping on the diagonal. Two of the flaps form the "base," a planar lattice: $\mathsf{C}_3^2$, the third one (shaded) comes up out of the plane pointing in the direction of the viewer.

We get some more examples of $\mathsf{M}_3[L]$ from the following observation.

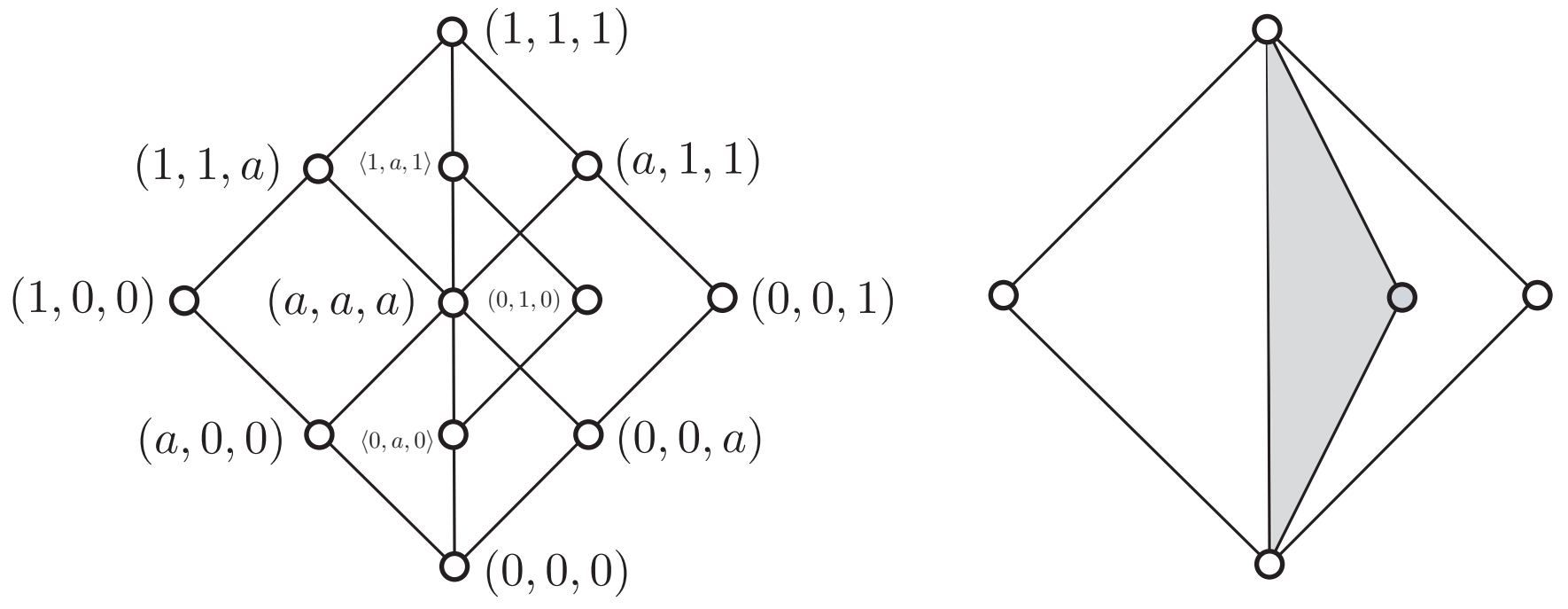

Figure 6.1: The lattice $\mathsf{M}_3[\mathsf{C}_3]$ with a sketch

**Lemma 6.2.** *Let the lattice $L$ have a direct product decomposition: $L = L_1 \times L_2$. Then $\mathsf{M}_3[L] \cong \mathsf{M}_3[L_1] \times \mathsf{M}_3[L_2]$.*

*Proof.* This is obvious, since $((x_1, y_1), (x_2, y_2), (x_3, y_3))$ is a Boolean triple iff $(x_1, x_2, x_3)$ and $(y_1, y_2, y_3)$ are both Boolean triples ($x_i \in L_1$ and $y_i \in L_2$ for $i = 1, 2, 3$). □

## 6.2. The congruence-preserving extension property

Let $L$ be a nontrivial lattice with zero and let

$$\gamma\colon x \mapsto (x, 0, 0) \in \mathsf{M}_3[L]$$

be an embedding of $L$ into $\mathsf{M}_3[L]$.

Here is the main result of my joint paper with F. Wehrung [199].

**Theorem 6.3.** *$\mathsf{M}_3[L]$ is a congruence-preserving extension of $\gamma L$.*

The next two lemmas prove this theorem.

For a congruence $\boldsymbol{\alpha}$ of $L$, let $\boldsymbol{\alpha}^3$ denote the congruence of $L^3$ defined componentwise. Let $\mathsf{M}_3[\boldsymbol{\alpha}]$ be the restriction of $\boldsymbol{\alpha}^3$ to $\mathsf{M}_3[L]$.

**Lemma 6.4.** *$\mathsf{M}_3[\boldsymbol{\alpha}]$ is a congruence relation of $\mathsf{M}_3[L]$.*

*Proof.* $\mathsf{M}_3[\boldsymbol{\alpha}]$ is obviously an equivalence relation on $\mathsf{M}_3[L]$. Since $\mathsf{M}_3[L]$ is a meet-subsemilattice of $L^3$, it is clear that $\mathsf{M}_3[\boldsymbol{\alpha}]$ satisfies ($\mathrm{SP}_\wedge$). To verify ($\mathrm{SP}_\vee$) for $\mathsf{M}_3[\boldsymbol{\alpha}]$, let $(x_1, y_1, z_1), (x_2, y_2, z_2) \in \mathsf{M}_3[L]$, let

$$(x_1, y_1, z_1) \equiv (x_2, y_2, z_2) \pmod{\mathsf{M}_3[\boldsymbol{\alpha}]},$$

and let $(u, v, w) \in \mathsf{M}_3[L]$. Set

$$(x_i', y_i', z_i') = (x_i, y_i, z_i) \vee (u, v, w)$$

(the join formed in $\mathsf{M}_3[L]$) for $i = 1, 2$.

Then using Lemma 6.1.(ii) for $x_1 \vee u$, $y_1 \vee v$, and $z_1 \vee w$, we obtain that

$$\begin{aligned} x_1' &= (x_1 \vee u \vee y_1 \vee v) \wedge (x_1 \vee u \vee z_1 \vee w) \\ &\equiv (x_2 \vee u \vee y_2 \vee v) \wedge (x_2 \vee u \vee z_2 \vee w) = x_2' \pmod{\mathsf{M}_3[\boldsymbol{\alpha}]}, \end{aligned}$$

and similarly, $y_1' \equiv y_2' \pmod{\mathsf{M}_3[\boldsymbol{\alpha}]}$, $z_1' \equiv z_2' \pmod{\mathsf{M}_3[\boldsymbol{\alpha}]}$, hence,

$$(x_1, y_1, z_1) \vee (u, v, w) \equiv (x_2, y_2, z_2) \vee (u, v, w) \pmod{\mathsf{M}_3[\boldsymbol{\alpha}]}.$$ □

It is obvious that $\mathsf{M}_3[\boldsymbol{\alpha}]$ restricted to $\gamma L$ is $\gamma\boldsymbol{\alpha}$.

**Lemma 6.5.** *Every congruence of $\mathsf{M}_3[L]$ is of the form $\mathsf{M}_3[\boldsymbol{\alpha}]$ for a suitable congruence $\boldsymbol{\alpha}$ of $L$.*

*Proof.* Let $\boldsymbol{\gamma}$ be a congruence of $\mathsf{M}_3[L]$, and let $\boldsymbol{\alpha}$ denote the congruence of $L$ obtained by restricting $\boldsymbol{\gamma}$ to the sublattice $L' = \{\, (x, 0, 0) \mid x \in L \,\}$ of $\mathsf{M}_3[L]$, that is, for $x, y \in L$, $x \equiv y \pmod{\boldsymbol{\alpha}}$ iff $(x, 0, 0) \equiv (y, 0, 0) \pmod{\boldsymbol{\gamma}}$. We prove that $\boldsymbol{\gamma} = \mathsf{M}_3[\boldsymbol{\alpha}]$.

To show that $\boldsymbol{\gamma} \subseteq \mathsf{M}_3[\boldsymbol{\alpha}]$, let

$$(x_1, y_1, z_1) \equiv (x_2, y_2, z_2) \pmod{\boldsymbol{\gamma}}. \tag{1}$$

Meeting the congruence (1) with $(1, 0, 0)$ yields

$$(x_1, 0, 0) \equiv (x_2, 0, 0) \pmod{\boldsymbol{\gamma}},$$

and so

$$(x_1, 0, 0) \equiv (x_2, 0, 0) \pmod{\mathsf{M}_3[\boldsymbol{\alpha}]}. \tag{2}$$

Meeting the congruence (1) with $(0, 1, 0)$ yields

$$(0, y_1, 0) \equiv (0, y_2, 0) \pmod{\boldsymbol{\gamma}}. \tag{3}$$

Since

$$(0, y_1, 0) \vee (0, 0, 1) = \overline{(0, y_1, 1)} = (y_1, y_1, 1),$$

and

$$(y_1, y_1, 1) \wedge (1, 0, 0) = (y_1, 0, 0),$$

and similarly for $(0, y_2, 0)$, joining the congruence (3) with $(0, 1, 0)$ and then meeting with $(1, 0, 0)$, yields

$$(y_1, 0, 0) \equiv (y_2, 0, 0) \pmod{\boldsymbol{\gamma}},$$

and so

$$(0, y_1, 0) \equiv (0, y_2, 0) \pmod{\mathsf{M}_3[\boldsymbol{\alpha}]}. \tag{4}$$

Similarly,

$$(0, 0, z_1) \equiv (0, 0, z_2) \pmod{\mathsf{M}_3[\boldsymbol{\alpha}]}. \tag{5}$$

Joining the congruences (2), (4), and (5), we obtain

$$(x_1, y_1, z_1) \equiv (x_2, y_2, z_2) \pmod{\mathsf{M}_3[\boldsymbol{\alpha}]}, \tag{6}$$

proving that $\boldsymbol{\gamma} \subseteq \mathsf{M}_3[\boldsymbol{\alpha}]$.

To prove the converse, $\mathsf{M}_3[\boldsymbol{\alpha}] \subseteq \boldsymbol{\gamma}$, take

$$(x_1, y_1, z_1) \equiv (x_2, y_2, z_2) \pmod{\mathsf{M}_3[\boldsymbol{\alpha}]} \tag{7}$$

in $\mathsf{M}_3[L]$; equivalently,

$$(x_1, 0, 0) \equiv (x_2, 0, 0) \pmod{\boldsymbol{\gamma}}, \tag{8}$$
$$(y_1, 0, 0) \equiv (y_2, 0, 0) \pmod{\boldsymbol{\gamma}}, \tag{9}$$
$$(z_1, 0, 0) \equiv (z_2, 0, 0) \pmod{\boldsymbol{\gamma}} \tag{10}$$

in $\mathsf{M}_3[L]$.

Joining the congruence (9) with $(0, 0, 1)$ and then meeting the result with $(0, 1, 0)$, we get (as in the computation following (3)):

$$(0, y_1, 0) \equiv (0, y_2, 0) \pmod{\boldsymbol{\gamma}}. \tag{11}$$

Similarly, from (10), we conclude that

$$(0, 0, z_1) \equiv (0, 0, z_2) \pmod{\boldsymbol{\gamma}}. \tag{12}$$

Finally, joining the congruences (8), (11), and (12), we get

$$(x_1, y_1, z_1) \equiv (x_2, y_2, z_2) \pmod{\boldsymbol{\gamma}}, \tag{13}$$

that is, $\mathsf{M}_3[\boldsymbol{\alpha}] \subseteq \boldsymbol{\gamma}$. This completes the proof of this lemma. ☐

## 6.3. The distributive case

Let $D$ be a bounded distributive lattice. In 1974 (25 years before the publication of my joint paper with F. Wehrung [199]), E. T. Schmidt [238] defined $\mathsf{M}_3[D]$ as the set of balanced triples $(x, y, z) \in D^3$ (defined in Section 6.1), regarded as an ordered subset of $D^3$. Schmidt proved the following result.

**Theorem 6.6.** *Let $D$ be a bounded distributive lattice. Then $\mathsf{M}_3[D]$ (the set of balanced triples $(x, y, z) \in D^3$) is a modular lattice. The map*

$$\gamma \colon x \mapsto (x, 0, 0) \in \mathsf{M}_3[D]$$

*is an embedding of $D$ into $\mathsf{M}_3[D]$, and $\mathsf{M}_3[D]$ is a congruence-preserving extension of $\gamma D$.*

*Proof.* Examining Figure 2.6 (page 24), we immediately see that in a distributive lattice $D$, conditions (B) and (F) are equivalent, so $\mathsf{M}_3[D]$ is the Boolean triple construction. Therefore, this result follows from the results in Sections 6.1 and 6.2, except for the modularity. A direct computation of this is not so easy—although entertaining. However, we can do it without computation. Observe that it is enough to prove modularity for a finite $D$. Now if $D$ is finite, then by Lemma 6.2, $\mathsf{M}_3[D]$ can be embedded into $\mathsf{M}_3[\mathsf{B}_n] \cong (\mathsf{M}_3[\mathsf{C}_2])^n \cong \mathsf{M}_3^n$, a modular lattice; hence, $\mathsf{M}_3[D]$ is modular. ☐

## 6.4. Two interesting intervals

Let $L$ be a bounded lattice. For an arbitrary $a \in L$, consider the interval $\mathsf{M}_3[L, a]$ of $\mathsf{M}_3[L]$ (illustrated in Figure 6.2):

$$\mathsf{M}_3[L, a] = [(0, a, 0), (1, 1, 1)] \subseteq \mathsf{M}_3[L].$$

Define the map (with image $B = \gamma_a L$ in Figure 6.2):

$$\gamma_a \colon x \mapsto (x, a, x \wedge a).$$

**Lemma 6.7.** *$\gamma_a$ is an embedding of $L$ into $\mathsf{M}_3[L, a]$.*

*Proof.* $\gamma_a$ is obviously one-to-one and meet-preserving. It also preserves the join because

$$\overline{(x \vee y, a, (x \wedge a) \vee (y \wedge a))} = (x \vee y, a, (x \vee y) \wedge a)$$

for $x, y \in L$. ☐

The following result of my joint paper with E. T. Schmidt [186] is a generalization of Theorem 6.3.

**Theorem 6.8.** $\mathsf{M}_3[L, a]$ *is a congruence-preserving extension of* $\gamma_a L$.

The proof is presented in the following lemmas.

**Lemma 6.9.** *Let $(x, y, z) \in \mathsf{M}_3[L, a]$. Then $a \leq y$ and*

$$x \wedge a = z \wedge a. \tag{14}$$

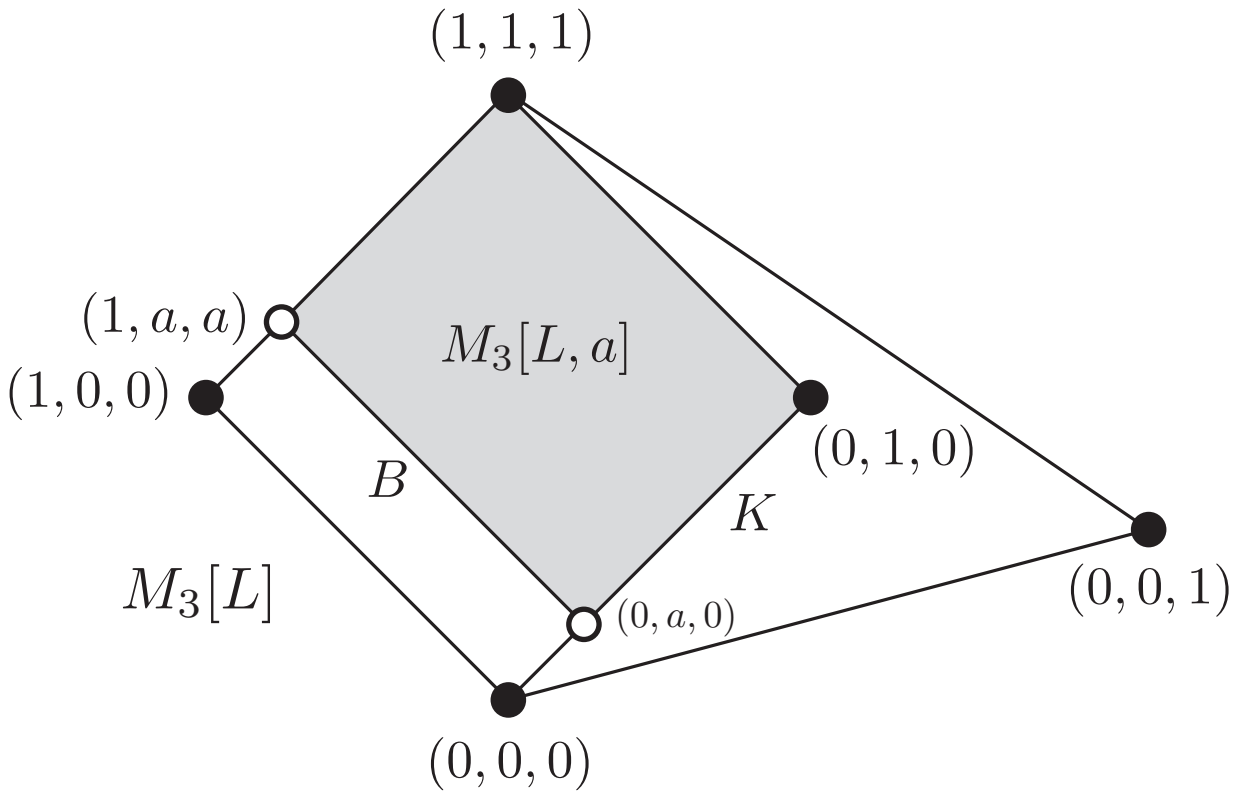

Figure 6.2: The shaded area: the lattice $\mathsf{M}_3[L, a]$

*Proof.* Since $(x, y, z)$ is Boolean, it is balanced, so $x \wedge y = z \wedge y$. Therefore, $x \wedge a = (x \wedge y) \wedge a = (z \wedge y) \wedge a = z \wedge a$, as claimed. $\square$

We need an easy decomposition statement for the elements of $\mathsf{M}_3[L, a]$. Let us use the notation ($B$ and $K$ are marked in Figure 6.2):

$$\begin{aligned} B &= \{\, (x, a, x \wedge a) \mid x \in L \,\} \ (= \gamma_a L), \\ K &= \{\, (0, x, 0) \mid x \in L,\ x \geq a \,\}, \\ J &= \{\, (x \wedge a, a, x) \mid x \in L \,\}. \end{aligned}$$

**Lemma 6.10.** *Let* $\mathbf{v} = (x, y, z) \in \mathsf{M}_3[L, a]$. *Then* $\mathbf{v}$ *has a decomposition in* $\mathsf{M}_3[L, a]$*:*

$$\mathbf{v} = \mathbf{v}_B \vee \mathbf{v}_K \vee \mathbf{v}_J, \tag{15}$$

*where*

$$\mathbf{v}_B = (x, y, z) \wedge (1, a, a) = (x, a, x \wedge a) \in B, \tag{16}$$
$$\mathbf{v}_K = (x, y, z) \wedge (0, 1, 0) = (0, y, 0) \in K, \tag{17}$$
$$\mathbf{v}_J = (x, y, z) \wedge (a, a, 1) = (z \wedge a, a, z) \in J. \tag{18}$$

*Proof.* (16) follows from (14). By symmetry, (18) follows, and (17) is trivial. Finally, the right side of (15) componentwise joins to the left side—in view of (14). $\square$

We can now describe the congruences of $\mathsf{M}_3[L, a]$.

**Lemma 6.11.** *Let* $\boldsymbol{\gamma}$ *be a congruence of* $\mathsf{M}_3[L, a]$ *and let* $\mathbf{v}, \mathbf{w} \in \mathsf{M}_3[L, a]$. *Then*

$$\mathbf{v} \equiv \mathbf{w} \pmod{\boldsymbol{\gamma}}, \tag{19}$$

*iff*

$$\mathbf{v}_B \equiv \mathbf{w}_B \pmod{\boldsymbol{\gamma}}, \tag{20}$$
$$\mathbf{v}_K \equiv \mathbf{w}_K \pmod{\boldsymbol{\gamma}}, \tag{21}$$
$$\mathbf{v}_J \equiv \mathbf{w}_J \pmod{\boldsymbol{\gamma}}. \tag{22}$$

*Proof.* (19) implies (20) by (16). Similarly, for (21) and (22). Conversely, (20)–(22) imply (19) by (15). $\square$

**Lemma 6.12.** *For a congruence* $\boldsymbol{\alpha}$ *of* $L$, *let* $\mathsf{M}_3[\boldsymbol{\alpha}, a]$ *be the restriction of* $\boldsymbol{\alpha}^3$ *to* $\mathsf{M}_3[L, a]$. *Then* $\mathsf{M}_3[\boldsymbol{\alpha}, a]$ *is a congruence of* $\mathsf{M}_3[L, a]$, *and every congruence of* $\mathsf{M}_3[L, a]$ *is of the form* $\mathsf{M}_3[\boldsymbol{\alpha}, a]$, *for a unique congruence* $\boldsymbol{\alpha}$ *of* $L$.

*Proof.* It follows from Lemma 6.4 that $\mathsf{M}_3[\boldsymbol{\alpha}, a]$ is a congruence of $\mathsf{M}_3[L, a]$. Let

$$\mathbf{v} = (x, y, z), \mathbf{v}' = (x', y', z') \in \mathsf{M}_3[L, a].$$

Let $\boldsymbol{\gamma}$ be a congruence of $\mathsf{M}_3[L, a]$, and let $\boldsymbol{\alpha}$ be the restriction of $\boldsymbol{\gamma}$ to $L$ with respect to the embedding $\gamma_a$. By Lemma 6.11,

$$\mathbf{v} \equiv \mathbf{v}' \pmod{\boldsymbol{\gamma}}$$

iff (20)–(22) hold. Note that $\mathbf{v}_B, \mathbf{v}'_B \in \gamma_a L$, so (20) is equivalent to $\mathbf{v}_B \equiv \mathbf{v}'_B \pmod{\boldsymbol{\alpha}}$.

Now consider

$$p(\mathbf{x}) = (\mathbf{x} \vee (0, 1, 0)) \wedge (1, a, a).$$

Then

$$p((x \wedge a, a, x)) = (x, 1, x) \wedge (1, a, a) = (x, a, x \wedge a).$$

So (22) implies that $p(\mathbf{v}_J) \equiv p(\mathbf{v}'_J) \pmod{\boldsymbol{\gamma}}$, and symmetrically. Thus (22) is equivalent to

$$p(\mathbf{v}_J) \equiv p(\mathbf{v}'_J) \pmod{\boldsymbol{\gamma}},$$

that is, to

$$p(\mathbf{v}_J) \equiv p(\mathbf{v}'_J) \pmod{\boldsymbol{\alpha}},$$

since $p(\mathbf{v}_J), p(\mathbf{v}'_J) \in \gamma_a L$.

Now consider

$$q(\mathbf{x}) = (\mathbf{x} \vee (a, a, 1)) \wedge (1, a, a).$$

Then $q((0, x, 0)) = (x, a, x \wedge a)$ for $x \geq a$. So $q(\mathbf{v}_K) \equiv q(\mathbf{v}'_K) \pmod{\boldsymbol{\gamma}}$, that is, $q(\mathbf{v}_K) \equiv q(\mathbf{v}'_K) \pmod{\boldsymbol{\alpha}}$.

Finally, define

$$r(\mathbf{x}) = (\mathbf{x} \vee (a, a, 1)) \wedge (0, 1, 0).$$

Then $q((x, x \wedge a, a)) = (0, x, 0)$. So $q(\mathbf{v}_B) \equiv q(\mathbf{v}'_B) \pmod{\boldsymbol{\gamma}}$. From these it follows that $\mathbf{v}_K \equiv \mathbf{v}'_K \pmod{\boldsymbol{\gamma}}$ is equivalent to $q(\mathbf{v}_K) \equiv q(\mathbf{v}'_K) \pmod{\boldsymbol{\gamma}}$ and $q(\mathbf{v}_K), q(\mathbf{v}'_K) \in \gamma_a L$, so the latter is equivalent to $q(\mathbf{v}_K) \equiv q(\mathbf{v}'_K) \pmod{\boldsymbol{\alpha}}$.

We conclude that the congruence $\mathbf{v} \equiv \mathbf{v}' \pmod{\boldsymbol{\gamma}}$ in $\mathsf{M}_3[L, a]$ is equivalent to the following three congruences in $L$:

$$\begin{aligned} \mathbf{v}_B &\equiv \mathbf{v}'_B && \pmod{\boldsymbol{\alpha}}, \\ p(\mathbf{v}_J) &\equiv p(\mathbf{v}'_J) && \pmod{\boldsymbol{\alpha}}, \\ q(\mathbf{v}_K) &\equiv q(\mathbf{v}'_K) && \pmod{\boldsymbol{\alpha}}, \end{aligned}$$

concluding the proof of the lemma. ☐

In Chapter 17, we need a smaller interval introduced in my joint paper [190] with E. T. Schmidt; namely, for $a, b \in L$ with $a < b$, we introduce the interval $\mathsf{M}_3[L, a, b]$ of $\mathsf{M}_3[L]$:

$$\mathsf{M}_3[L, a, b] = [(0, a, 0), (1, b, b)] \subseteq \mathsf{M}_3[L].$$

Again,

$$\gamma_a : x \mapsto (x, a, x \wedge a)$$

is a (convex) embedding of $L$ into $\mathsf{M}_3[L, a, b]$. (Note that if $L$ is bounded, then $\mathsf{M}_3[L, a] = \mathsf{M}_3[L, a, 1]$.)

Using the notation (illustrated in Figure 6.4; the black-filled elements form $J$):

$$\begin{aligned} B &= \{\, (x, a, x \wedge a) \mid x \in L \,\} \ (= \gamma_a L), \\ I_{a,b} &= [(0, a, 0), (0, b, 0)], \\ J &= \{\, (x \wedge a, a, x) \mid x \leq b \,\}. \end{aligned}$$

We can now generalize Lemmas 6.10–6.12.

**Lemma 6.13.** *Let* $\mathbf{v} = (x, y, z) \in \mathsf{M}_3[L, a, b]$. *Then* $\mathbf{v}$ *has a decomposition in* $\mathsf{M}_3[L, a, b]$*:*

$$\mathbf{v} = \mathbf{v}_B \vee \mathbf{v}_{I_{a,b}} \vee \mathbf{v}_J,$$

*where*

$$\begin{aligned} \mathbf{v}_B &= (x, y, z) \wedge (1, a, a) = (x, a, x \wedge a) \in B, \\ \mathbf{v}_{I_{a,b}} &= (x, y, z) \wedge (0, b, 0) \ = (0, y, 0) \in I_{a,b}, \\ \mathbf{v}_J &= (x, y, z) \wedge (a, a, b) \ = (z \wedge a, a, z) \in J. \end{aligned}$$

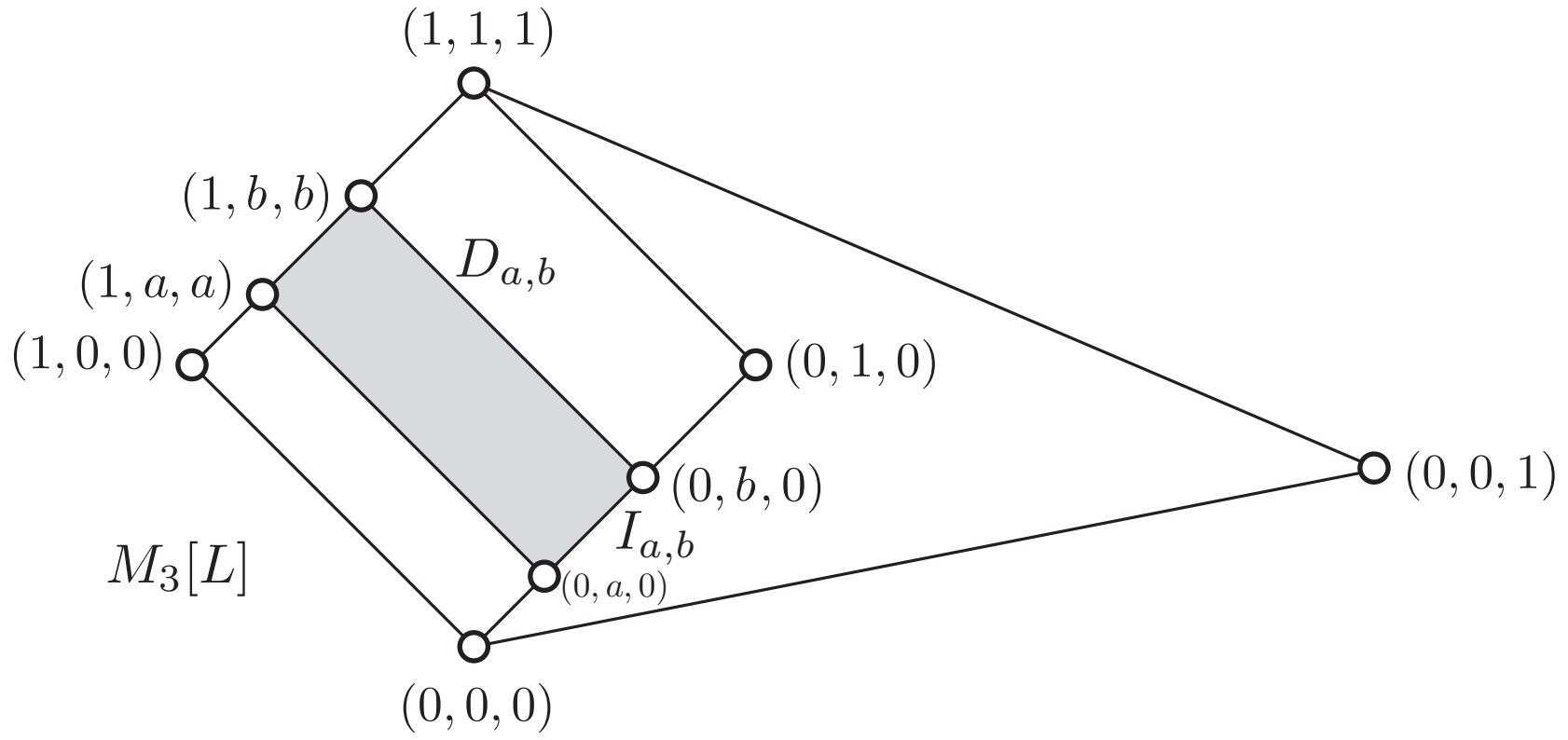

Figure 6.3: The shaded area: the lattice $\mathsf{M}_3[L, a, b]$

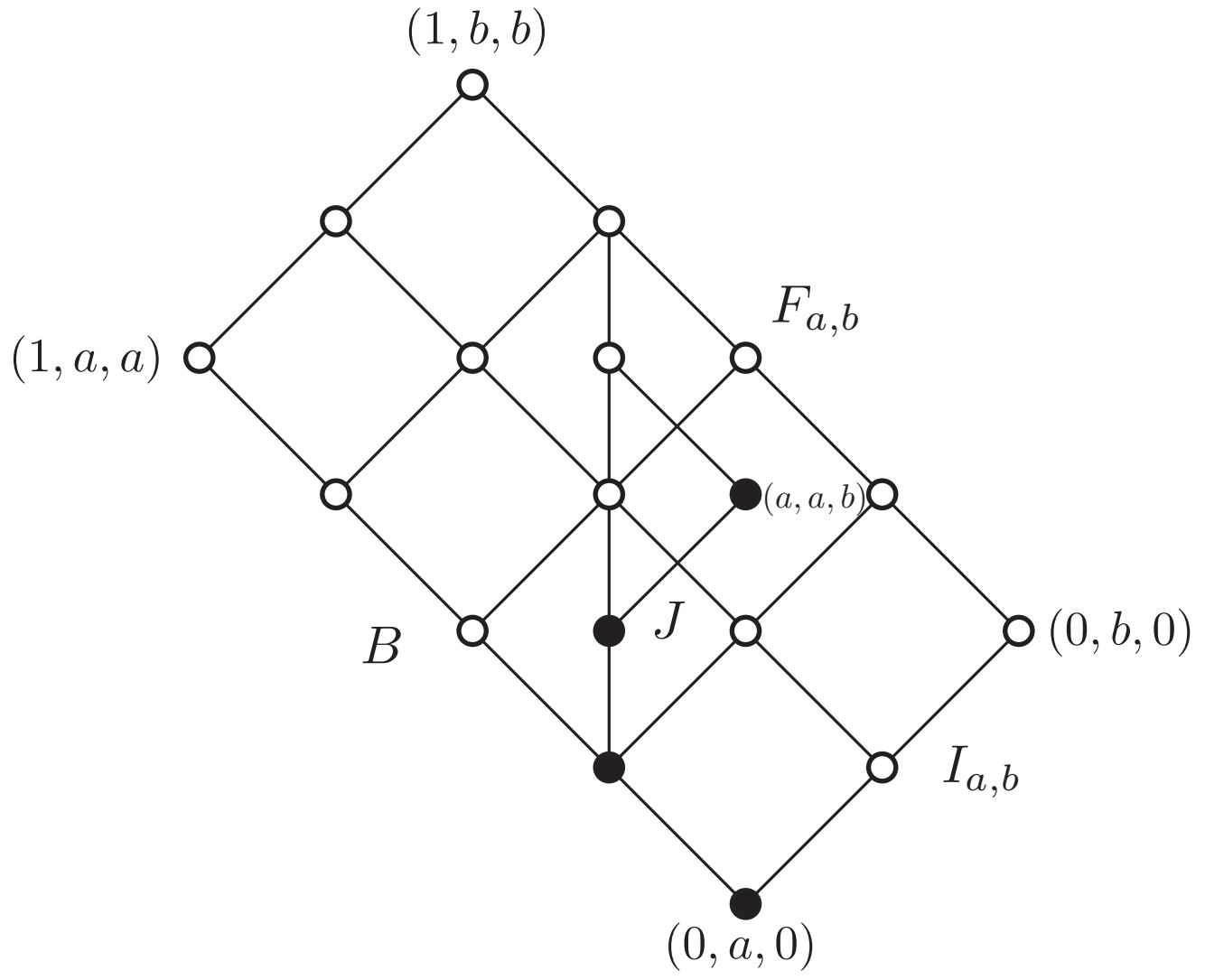

Figure 6.4: The lattice $\mathsf{M}_3[L, a, b]$

**Lemma 6.14.** *Let* $\boldsymbol{\gamma}$ *be a congruence of* $\mathsf{M}_3[L, a, b]$ *and let* $\mathbf{v}, \mathbf{w} \in \mathsf{M}_3[L, a, b]$. *Then*

$$\mathbf{v} \equiv \mathbf{w} \pmod{\boldsymbol{\gamma}}$$

*iff*

$$\mathbf{v}_B \equiv \mathbf{w}_B \pmod{\boldsymbol{\gamma}}, \tag{23}$$

$$\mathbf{v}_{I_{a,b}} \equiv \mathbf{w}_{I_{a,b}} \pmod{\boldsymbol{\gamma}}, \tag{24}$$

$$\mathbf{v}_J \equiv \mathbf{w}_J \pmod{\boldsymbol{\gamma}}. \tag{25}$$

**Lemma 6.15.** *For a congruence* $\boldsymbol{\alpha}$ *of* $L$, *let* $\mathsf{M}_3[\boldsymbol{\alpha}, a, b]$ *be the restriction of* $\boldsymbol{\alpha}^3$ *to* $\mathsf{M}_3[L, a, b]$. *Then* $\mathsf{M}_3[\boldsymbol{\alpha}, a, b]$ *is a congruence of* $\mathsf{M}_3[L, a, b]$, *and every congruence of* $\mathsf{M}_3[L, a, b]$ *is of the form* $\mathsf{M}_3[\boldsymbol{\alpha}, a, b]$ *for a unique congruence* $\boldsymbol{\alpha}$ *of* $L$. *It follows that* $\gamma_a$ *is a congruence-preserving convex embedding of* $L$ *into* $\mathsf{M}_3[L, a, b]$.

We will use the notation

$$F_{a,b} = [(0, b, 0), (1, b, b)] = \{\, (x, b, x \wedge b) \mid x \in L \,\}.$$

The following two observations are trivial.

**Lemma 6.16.**

(i) $I_{a,b}$ *is an ideal of* $\mathsf{M}_3[L, a, b]$ *and* $I_{a,b}$ *is isomorphic to the interval* $[a, b]$ *of* $L$.

(ii) *$F_{a,b}$ is a filter of $\mathsf{M}_3[L, a, b]$ and $F_{a,b}$ is isomorphic to $L$.*

We can say a lot more about $F_{a,b}$. But first we need another lemma.

**Lemma 6.17.** *Let $L$ be a lattice, let $[u, v]$ and $[u', v']$ be intervals of $L$, and let $\delta \colon [u, v] \mapsto [u', v']$ be an isomorphism between these two intervals. Let $\delta$ and $\delta^{-1}$ be term functions in $L$. Then $L$ is a congruence-preserving extension of $[u, v]$ iff it is a congruence-preserving extension of $[u', v']$.*

*Proof.* Let us assume that $L$ is a congruence-preserving extension of $[u', v']$. Let $\boldsymbol{\alpha}$ be a congruence relation of $[u, v]$ and let $\delta\boldsymbol{\alpha}$ be the image of $\boldsymbol{\alpha}$ under $\delta$. Since $\delta$ is an isomorphism, it follows that $\delta\boldsymbol{\alpha}$ is a congruence of $[u', v']$, and so $\delta\boldsymbol{\alpha}$ has a unique extension to a congruence $\boldsymbol{\gamma}$ of $L$.

*We claim that $\boldsymbol{\gamma}$ extends $\boldsymbol{\alpha}$ to $L$ and extends it uniquely.*

**Step 1.** *$\boldsymbol{\gamma}$ extends $\boldsymbol{\alpha}$.*

Let $x \equiv y \pmod{\boldsymbol{\alpha}}$. Then $\delta x \equiv \delta y \pmod{\delta\boldsymbol{\alpha}}$, since $\delta$ is an isomorphism. By definition, $\delta$ extends $\delta\boldsymbol{\alpha}$, so $\delta x \equiv \delta y \pmod{\boldsymbol{\alpha}}$. Since $\delta^{-1}$ is a term function, the last congruence implies that $x \equiv y \pmod{\boldsymbol{\alpha}}$. Conversely, let $x \equiv y \pmod{\boldsymbol{\gamma}}$ and $x, y \in [u, v]$. Then $\delta x, \delta y \in [u', v']$; since $\boldsymbol{\gamma}$ is a term function, it follows that $\delta x \equiv \delta y \pmod{\boldsymbol{\gamma}}$. Since $\boldsymbol{\gamma}$ extends $\delta\boldsymbol{\alpha}$, we conclude that $\delta x \equiv \delta y \pmod{\delta\boldsymbol{\alpha}}$. Using that $\delta^{-1}$ is an isomorphism, we obtain that $x \equiv y \pmod{\boldsymbol{\alpha}}$, verifying the claim.

**Step 2.** *$\boldsymbol{\gamma}$ extends $\boldsymbol{\alpha}$ uniquely.*

Let $\boldsymbol{\gamma}$ extend $\boldsymbol{\alpha}$ to $L$. As in the previous paragraph—mutatis mutandis—we conclude that $\boldsymbol{\gamma}$ extends $\delta\boldsymbol{\alpha}$ to $L$; hence, $\boldsymbol{\alpha} = \boldsymbol{\gamma}$, proving the uniqueness and the claim.

By symmetry, the lemma is proved. ☐

We have already observed in Lemma 6.12 that $\mathsf{M}_3[L, a, b]$ is a congruence-preserving convex extension of $[(0, a, 0), (1, a, a)]$. Since the isomorphism

$$x \mapsto x \vee (0, b, 0)$$

between $[(0, a, 0), (1, a, a)]$ and $[(0, b, 0), (1, b, b)] = F_{a,b}$ is a term function, and so is the inverse

$$x \mapsto x \wedge (1, a, a),$$

from Lemma 6.17 we conclude the following.

**Corollary 6.18.** *The lattice $\mathsf{M}_3[L, a, b]$ is a congruence-preserving convex extension of the filter $F_{a,b}$.*

We summarize our results (my joint paper with E. T. Schmidt [190]).

**Lemma 6.19.** *Let $L$ be a bounded lattice, and let $a, b \in L$ with $a < b$. Then there exists a bounded lattice $L_{a,b}$ (with bounds $0_{a,b}$ and $1_{a,b}$) and $u_{a,b}, v_{a,b} \in L_{a,b}$, such that the following conditions are satisfied:*

(i) *$v_{a,b}$ is a complement of $u_{a,b}$.*

(ii) *$F_{a,b} = [v_{a,b}, 1_{a,b}] \cong L$.*

(iii) *$I_{a,b} = [0_{a,b}, v_{a,b}] \cong [a, b]$.*

(iv) *$L_{a,b}$ is a congruence-preserving (convex) extension of $[0_{a,b}, u_{a,b}]$ and of $[v_{a,b}, 1_{a,b}]$.*

(v) *The congruences on $I_{a,b}$ and $F_{a,b}$ are* synchronized, *that is, if $\boldsymbol{\alpha}$ is a congruence on $L$, $\overline{\boldsymbol{\alpha}}$ is the extension of $\boldsymbol{\alpha}$ to $L_{a,b}$ (we map $\boldsymbol{\alpha}$ to $F_{a,b}$ under the isomorphism, and then by* (iv) *we uniquely extend it to $L_{a,b}$), and $x, y \in [a, b]$, then we can denote by $x_{F_{a,b}}$, $y_{F_{a,b}} \in F_{a,b}$ the images of $x, y$ in $F_{a,b}$ and by $x_{I_{a,b}}, y_{I_{a,b}} \in I_{a,b}$ the images of $x, y$ in $I_{a,b}$; synchronization means that*

$$x_{F_{a,b}} \equiv y_{F_{a,b}} \pmod{\overline{\boldsymbol{\alpha}}}$$

*is equivalent to*

$$x_{I_{a,b}} \equiv y_{I_{a,b}} \pmod{\overline{\boldsymbol{\alpha}}}.$$

*Proof.* Of course, $L_{a,b} = \mathsf{M}_3[L, a, b]$, $0_{a,b} = (0, a, 0)$, $1_{a,b} = (1, b, b)$, $u_{a,b} = (1, a, a)$, and $v_{a,b} = (0, b, 0)$. □

The lattice $L_{a,b}$ is illustrated with $L = \mathsf{C}_5$ in Figure 6.5, the five-element chain, $a$ is the atom and $b$ is the dual atom of $S$. This figure is the same as Figure 6.4, only the notation is changed.

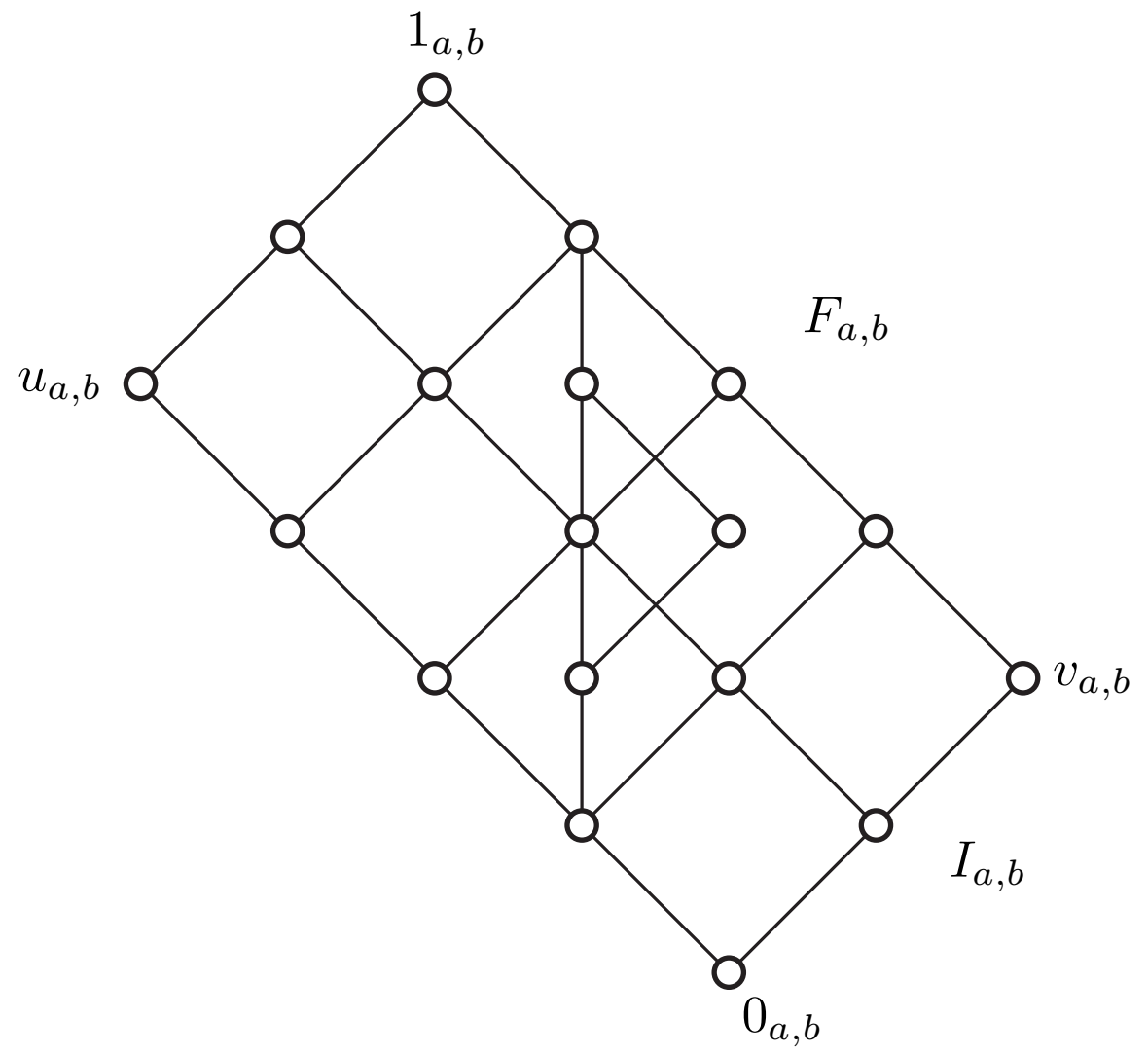

Figure 6.5: The lattice $L_{a,b}$

**Chapter**

# 7

# *Cubic Extensions*

In this chapter, for a finite lattice $K$, we introduce an extension $\operatorname{Cube} K$ with the following properties: (i) The lattice $K$ is a congruence-reflecting sublattice of $\operatorname{Cube} K$. (ii) $\operatorname{Con}(\operatorname{Cube} K)$ is Boolean. (iii) The minimal extension of the meet-irreducible congruences are the dual atoms of $\operatorname{Con}(\operatorname{Cube} K)$; their ordering is "flattened."

## 7.1. The construction

Let $K$ be a finite lattice. Following my joint paper with E. T. Schmidt [184], for every meet-irreducible congruence $\boldsymbol{\gamma}$ of $K$ (in symbols, $\boldsymbol{\gamma} \in \operatorname{Con}_{\mathrm{M}} K$), we form the quotient lattice $K/\boldsymbol{\gamma}$, and extend it to a finite *simple* lattice $\operatorname{Simp} K_{\boldsymbol{\gamma}}$ with zero $0_{\boldsymbol{\gamma}}$ and unit $1_{\boldsymbol{\gamma}}$, using Lemma 2.3.

Let $\operatorname{Cube}_{\mathrm{Simp}} K$ be the direct product of the lattices $\operatorname{Simp} K_{\boldsymbol{\gamma}}$ for $\boldsymbol{\gamma} \in \operatorname{Con}_{\mathrm{M}} K$:

$$\operatorname{Cube}_{\mathrm{Simp}} K = \prod(\operatorname{Simp} K_{\boldsymbol{\gamma}} \mid \boldsymbol{\gamma} \in \operatorname{Con}_{\mathrm{M}} K).$$

If the function Simp is understood, we write $\operatorname{Cube} K$ for $\operatorname{Cube}_{\mathrm{Simp}} K$. We regard $\operatorname{Simp} K_{\boldsymbol{\gamma}}$, for $\boldsymbol{\gamma} \in \operatorname{Con}_{\mathrm{M}} K$, an ideal of $\operatorname{Cube} K$, as in Section 2.2.

For $a \in K$, define $\operatorname{Diag}(a) \in \operatorname{Cube} K$ as follows.

$$\operatorname{Diag}(a) = (a/\boldsymbol{\gamma} \mid \boldsymbol{\gamma} \in \operatorname{Con}_{\mathrm{M}} K).$$

$K$ has a natural (diagonal) embedding into $\operatorname{Cube} K$ by

$$\gamma \colon a \mapsto \operatorname{Diag}(a) \quad \text{for } a \in K.$$

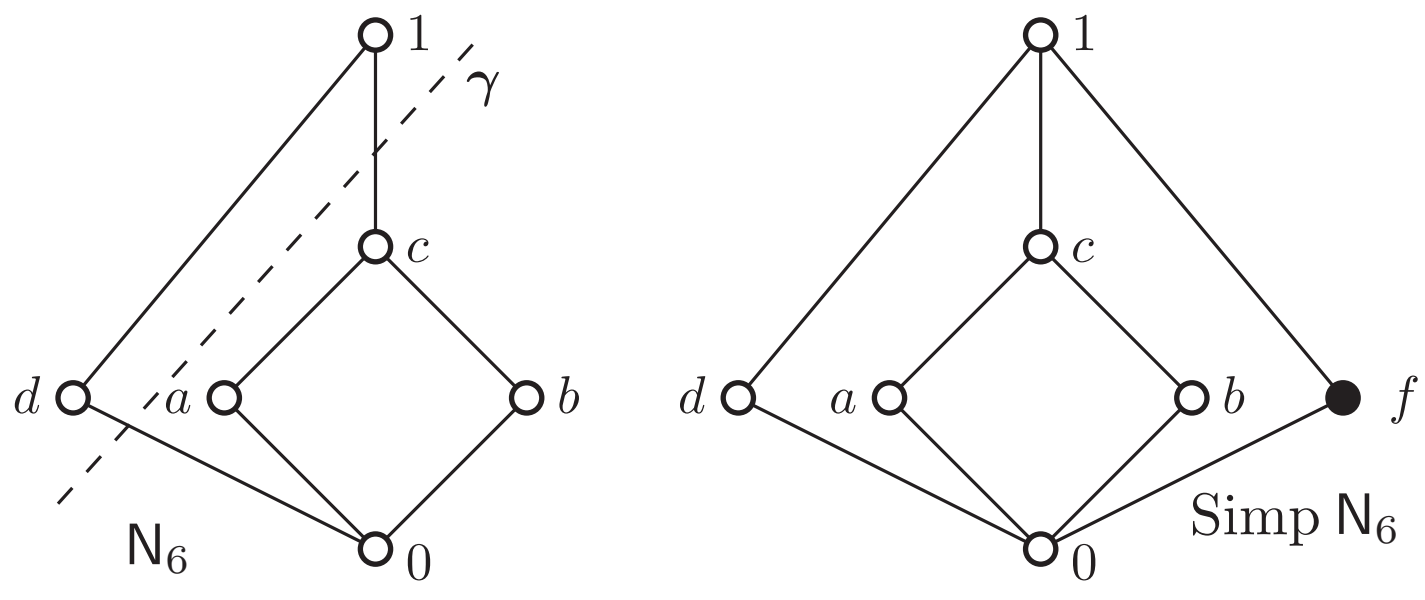

Figure 7.1: The lattices $\mathsf{N}_6$ and Simp $\mathsf{N}_6$

Let Diag$(K)$ be the sublattice $\{\,\mathrm{Diag}(a) \mid a \in K\,\}$ of Cube $K$, and for a congruence $\boldsymbol{\alpha}$ of $K$, let Diag$(\boldsymbol{\alpha})$ denote the corresponding congruence of Diag$(K)$. By identifying $a$ with Diag$(a)$, for $a \in K$, we can view Cube $K$ as an extension of $K$; we call Cube $K$ a *cubic extension* of $K$.

Cubic extensions are hard to draw because they are direct products. Here is a small example. We take the lattice $\mathsf{N}_6$ (see Figure 7.1). The lattice $\mathsf{N}_6$ is subdirectly irreducible; it has two meet-irreducible congruences, $\mathbf{0}$ and $\boldsymbol{\gamma}$. Since $\mathsf{N}_6/\boldsymbol{\gamma} = \mathsf{C}_2$, it is simple, so we can choose Simp $\mathsf{N}_6/\boldsymbol{\gamma} = \mathsf{C}_2$. The other lattice quotient is $\mathsf{N}_6/\mathbf{0} = \mathsf{N}_6$ and for this we choose a simple extension Simp $\mathsf{N}_6$, by adding an element; (see Figure 7.1). The lattice Cube $\mathsf{N}_6$ is shown in Figure 7.2 along with the embedding Diag of $\mathsf{N}_6$ into Cube $\mathsf{N}_6$. The images of elements of $\mathsf{N}_6$ under Diag are black-filled.

As an alternative, you may draw the chopped lattice $M$ whose ideal lattice is Cube $\mathsf{N}_6$ (see Figure 7.3). Unfortunately, the embedding Diag is not easy to see with this representation.

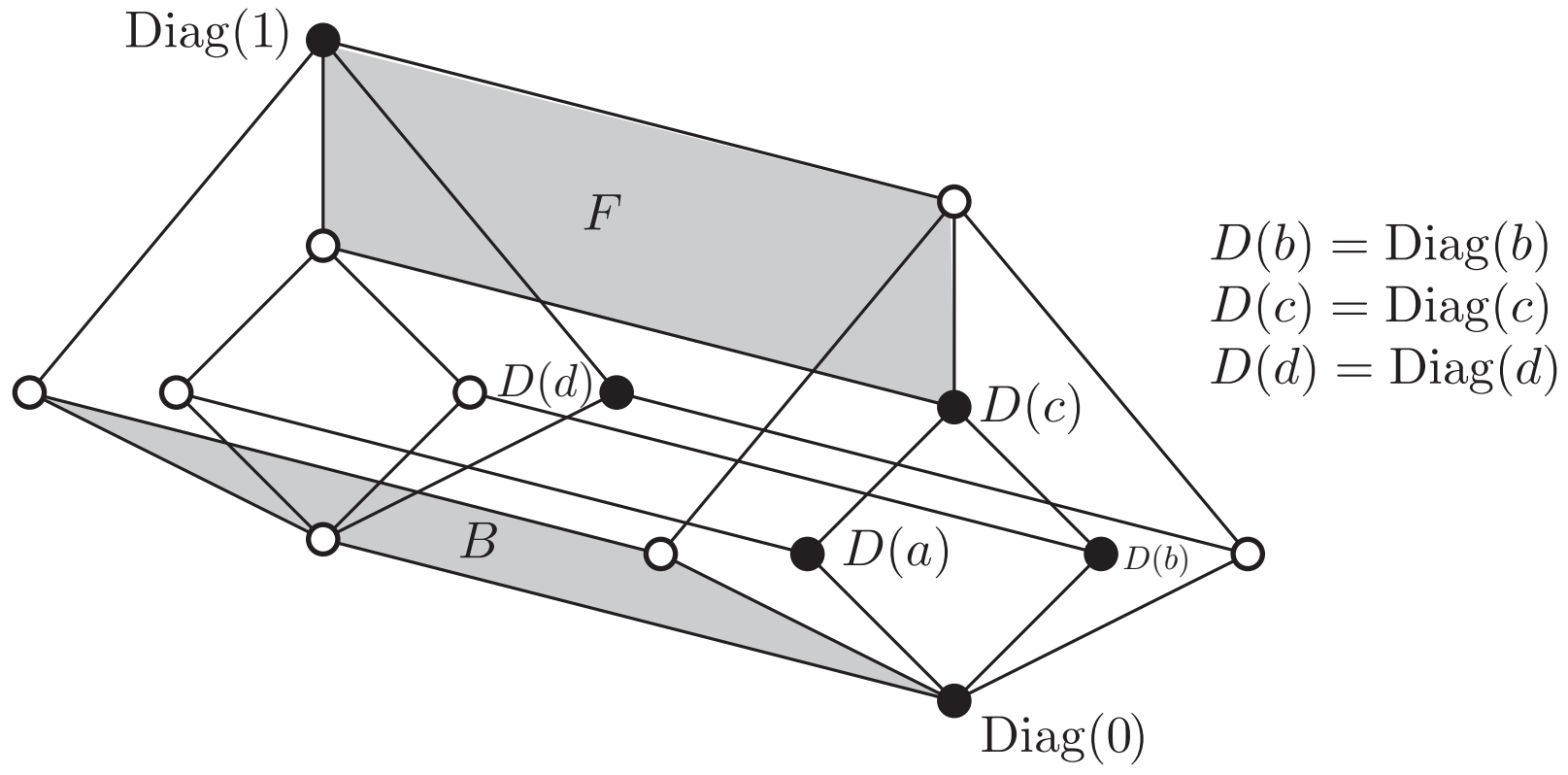

Figure 7.2: The embedding of $\mathsf{N}_6$ into Cube $\mathsf{N}_6$

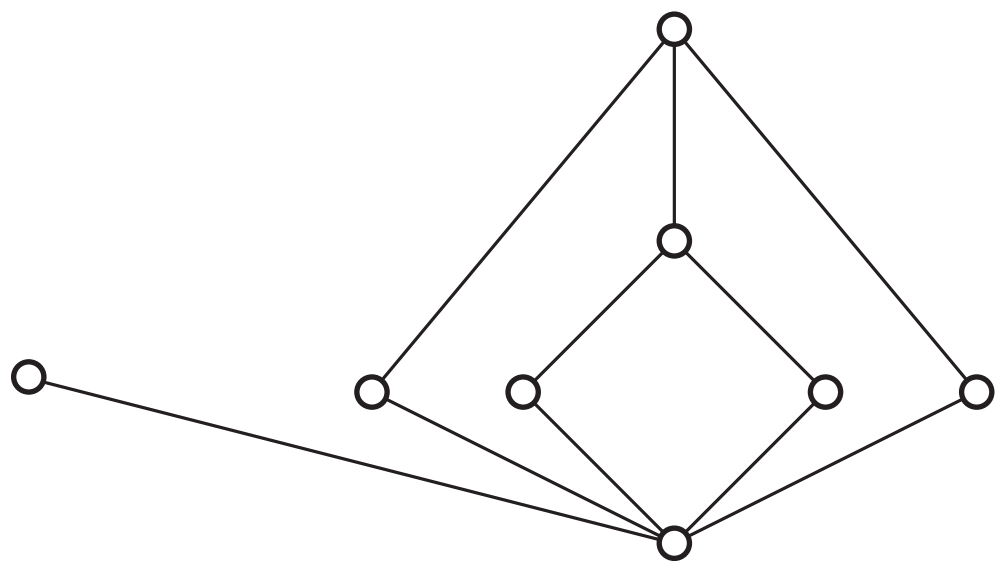

Figure 7.3: The chopped lattice $M$ whose ideals represent $\operatorname{Cube} \mathsf{N}_6$

## 7.2. The basic property

The following crucial property of cubic extensions was stated in my joint paper [184] with E. T. Schmidt.

**Theorem 7.1.** *Let $K$ be a finite lattice and let $\operatorname{Cube} K$ be a cubic extension of $K$. Then $K$ $(= \operatorname{Diag}(K))$ is a congruence-reflecting extension in $\operatorname{Cube} K$.*

*Proof.* For $\boldsymbol{\kappa} \in \operatorname{Con} K$ and $\boldsymbol{\gamma} \in \operatorname{Con_M} K$, define the congruence $\operatorname{Cube}(\boldsymbol{\kappa}, \boldsymbol{\gamma})$ on the lattice $\operatorname{Simp} K_{\boldsymbol{\gamma}}$ as follows.

$$\operatorname{Cube}(\boldsymbol{\kappa}, \boldsymbol{\gamma}) = \begin{cases} \mathbf{0} & \text{if } \boldsymbol{\kappa} \leq \boldsymbol{\gamma}, \\ \mathbf{1} & \text{if } \boldsymbol{\kappa} \nleq \boldsymbol{\gamma}; \end{cases}$$

and define

$$\operatorname{Cube}(\boldsymbol{\kappa}) = \prod(\, \operatorname{Cube}(\boldsymbol{\kappa}, \boldsymbol{\gamma}) \mid \boldsymbol{\gamma} \in \operatorname{Con_M} K\,),$$

a congruence of the lattice $\operatorname{Cube} K$.

Observe that $\operatorname{Cube}(\mathbf{0}_K) = \mathbf{0}_{\operatorname{Cube} K}$ and $\operatorname{Cube}(\mathbf{1}_K) = \mathbf{1}_{\operatorname{Cube} K}$, that is, $\operatorname{Cube}(\mathbf{0}) = \mathbf{0}$ and $\operatorname{Cube}(\mathbf{1}) = \mathbf{1}$.

We show that $\operatorname{Diag}(\boldsymbol{\kappa})$ is the restriction of $\operatorname{Cube}(\boldsymbol{\kappa})$ to $K$. First, assume that $\operatorname{Diag}(u) \equiv \operatorname{Diag}(v) \pmod{\operatorname{Diag}(\boldsymbol{\kappa})}$ in $\operatorname{Diag}(K)$. Then $u \equiv v \pmod{\boldsymbol{\kappa}}$ in $K$, by the definition of $D$. The congruence $u \equiv v \pmod{\boldsymbol{\kappa}}$ implies that $u \equiv v \pmod{\boldsymbol{\gamma}}$, for all $\boldsymbol{\gamma} \in \operatorname{Con_M} K$ satisfying $\boldsymbol{\kappa} \leq \boldsymbol{\gamma}$, that is, $u/\boldsymbol{\gamma} = v/\boldsymbol{\gamma}$, for all $\boldsymbol{\gamma} \in \operatorname{Con_M} K$ satisfying $\boldsymbol{\kappa} \leq \boldsymbol{\gamma}$. This, in turn, can be written as

$$u/\boldsymbol{\gamma} \equiv v/\boldsymbol{\gamma} \pmod{\operatorname{Cube}(\boldsymbol{\kappa}, \boldsymbol{\gamma})}, \text{ for all } \boldsymbol{\gamma} \in \operatorname{Con_M} K \text{ satisfying } \boldsymbol{\kappa} \leq \boldsymbol{\gamma},$$

since, by definition, $\boldsymbol{\kappa}_{\boldsymbol{\gamma}} = \boldsymbol{\kappa}$ for $\boldsymbol{\kappa} \leq \boldsymbol{\gamma}$.

Again, by definition, $\boldsymbol{\kappa}_{\boldsymbol{\gamma}} = \mathbf{1}$ for $\boldsymbol{\kappa} \nleq \boldsymbol{\gamma}$. Therefore, the congruence $u/\boldsymbol{\gamma} \equiv v/\boldsymbol{\gamma} \pmod{\operatorname{Cube}(\boldsymbol{\kappa}, \boldsymbol{\gamma})}$ always holds.

We conclude that $u/\boldsymbol{\gamma} \equiv v/\boldsymbol{\gamma} \pmod{\mathrm{Cube}(\boldsymbol{\kappa}, \boldsymbol{\gamma})}$ for all $\boldsymbol{\gamma} \in \mathrm{Con_M}\, K$. This congruence is equivalent to $\mathrm{Diag}(u) \equiv \mathrm{Diag}(v) \pmod{\mathrm{Cube}(\boldsymbol{\kappa})}$ in $\mathrm{Cube}\, K$, which was to be proved.

Second, assume that $\mathrm{Diag}(u) \equiv \mathrm{Diag}(v) \pmod{\mathrm{Cube}(\boldsymbol{\kappa})}$ in $\mathrm{Cube}\, K$. Then $u/\boldsymbol{\gamma} \equiv v/\boldsymbol{\gamma} \pmod{\mathrm{Cube}(\boldsymbol{\kappa}, \boldsymbol{\gamma})}$, for all $\boldsymbol{\gamma} \in \mathrm{Con_M}\, K$; in particular, for all $\boldsymbol{\gamma} \in \mathrm{Con_M}\, K$ satisfying $\boldsymbol{\kappa} \leq \boldsymbol{\gamma}$. Thus $u/\boldsymbol{\gamma} = v/\boldsymbol{\gamma}$, for all $\boldsymbol{\gamma} \in \mathrm{Con_M}\, K$ satisfying $\boldsymbol{\kappa} \leq \boldsymbol{\gamma}$, that is, $u \equiv v \pmod{\boldsymbol{\gamma}}$ for all $\boldsymbol{\gamma} \in \mathrm{Con_M}\, K$ satisfying $\boldsymbol{\kappa} \leq \boldsymbol{\gamma}$. Therefore,

$$u \equiv v \pmod{\bigwedge(\, \boldsymbol{\gamma} \in \mathrm{Con_M}\, K \mid \boldsymbol{\kappa} \leq \boldsymbol{\gamma}\,)}.$$

The lattice $\mathrm{Con}\, K$ is finite, so every congruence is a meet of meet-irreducible congruences; therefore,

$$\boldsymbol{\kappa} = \bigwedge(\, \boldsymbol{\gamma} \in \mathrm{Con_M}\, K \mid \boldsymbol{\kappa} \leq \boldsymbol{\gamma}\,),$$

and so $u \equiv v \pmod{\boldsymbol{\kappa}}$ in $K$, that is, $\mathrm{Diag}(u) \equiv \mathrm{Diag}(v) \pmod{\mathrm{Diag}(\boldsymbol{\kappa})}$ in $\mathrm{Diag}(K)$. ☐

For $\boldsymbol{\kappa} \in \mathrm{Con}\, K$, the set

$$\Delta_{\boldsymbol{\kappa}} = \{\, \boldsymbol{\gamma} \in \mathrm{Con_M}\, K \mid \boldsymbol{\kappa} \nleq \boldsymbol{\gamma}\,\}$$

is a down set of $\mathrm{Con_M}\, K$, and every down set of $\mathrm{Con_M}\, K$ is of the form $\Delta_{\boldsymbol{\kappa}}$, for a unique $\boldsymbol{\kappa} \in \mathrm{Con}\, K$. The down set $\Delta_{\boldsymbol{\kappa}}$ of $\mathrm{Con_M}\, K$ describes $\mathrm{Cube}(\boldsymbol{\kappa})$, and conversely.

In the example of Section 7.1 (see Figures 7.1–7.3), the congruence lattice $\mathrm{Con}\, \mathsf{N}_6 = \{\mathbf{0}, \boldsymbol{\gamma}, \mathbf{1}\}$. Since $\mathrm{Cube}(\mathbf{0}) = \mathbf{0}$ and $\mathrm{Cube}(\mathbf{1}) = \mathbf{1}$, we only have to compute $\mathrm{Cube}(\boldsymbol{\gamma})$. Clearly, $\mathrm{Cube}(\boldsymbol{\gamma}) = \mathbf{0} \times \mathbf{1}$ (that is, $\mathrm{Cube}(\boldsymbol{\gamma}) = \mathbf{0}_{\mathsf{C}_2} \times \mathbf{1}_{\mathrm{Simp}\, \mathsf{N}_6}$); it splits $\mathrm{Cube}\, \mathsf{N}_6$ into two parts as shown by the dashed line in Figure 7.2.

We now summarize the properties of the lattice $\mathrm{Cube}\, K$.

**Theorem 7.2.** *Let $K$ be a finite lattice with a cubic extension* $\mathrm{Cube}\, K$. *Then*

(i) $\mathrm{Cube}\, K$ *is finite.*

(ii) *Choose an atom $s_{\boldsymbol{\gamma}} \in \mathrm{Simp}\, K_{\boldsymbol{\gamma}}$, for each $\boldsymbol{\gamma} \in \mathrm{Con_M}\, K$, and let $B$ be the Boolean ideal of* $\mathrm{Cube}\, K$ *generated by these atoms. Then $B$ is a congruence-determining ideal of* $\mathrm{Cube}\, K$.

(iii) *There are one-to-one correspondences among the subsets of* $\mathrm{Con_M}\, K$, *the sets of atoms of $B$, and the congruences $\boldsymbol{\kappa}$ of* $\mathrm{Cube}\, K$; *the subset of* $\mathrm{Con_M}\, K$ *corresponding to the congruence $\boldsymbol{\kappa}$ of* $\mathrm{Cube}\, K$ *is*

$$\boldsymbol{\delta}_{\boldsymbol{\kappa}} = \{\, \boldsymbol{\gamma} \in \mathrm{Con_M}\, K \mid \boldsymbol{\kappa} \nleq \boldsymbol{\gamma}\,\},$$

*which, in turn corresponds to the set*

$$S_{\boldsymbol{\kappa}} = \{\, s_{\boldsymbol{\gamma}} \mid s_{\boldsymbol{\gamma}} \equiv 0 \pmod{\boldsymbol{\kappa}} \,\}$$

*of atoms of $B$. Hence, the congruence lattice of* Cube $K$ *is a finite Boolean lattice.*

(iv) *Every congruence $\boldsymbol{\kappa}$ of $K$ has an extension* $\mathrm{Cube}(\boldsymbol{\kappa})$ *to a congruence of* Cube $K$ *corresponding to the down set $\Delta_{\boldsymbol{\kappa}}$ of* $\mathrm{Con_M}\,K$.

In the example of Section 7.1, for $B$ we can choose the shaded ideal of Figure 7.2.

**Corollary 7.3.** *Choose a dual atom $t_{\boldsymbol{\gamma}} \in \mathrm{Simp}\,K_{\boldsymbol{\gamma}}$, for each $\boldsymbol{\gamma} \in \mathrm{Con_M}\,K$, and define $\bar{t}_{\boldsymbol{\gamma}}$ as the element of* Cube $K$ *whose $\boldsymbol{\gamma}$-component is $t_{\boldsymbol{\gamma}}$ and all other components are* 1. *The element $\bar{t}_{\boldsymbol{\gamma}}$ is a dual atom of* Cube $K$. *Let $F$ be the Boolean filter of* Cube $K$ *generated by these dual atoms. Then $F$ is a congruence-determining filter of* Cube $K$.

*There are one-to-one correspondences among the subsets of* $\mathrm{Con_M}\,K$, *the sets of dual atoms of $F$, and the congruences $\boldsymbol{\kappa}$ of* Cube $K$; *the subset of* $\mathrm{Con_M}\,K$ *corresponding to the congruence $\boldsymbol{\kappa}$ of* Cube $K$ *is*

$$D_{\boldsymbol{\kappa}} = \{\, \boldsymbol{\gamma} \in \mathrm{Con_M}\,K \mid \boldsymbol{\kappa} \nleq \boldsymbol{\gamma} \,\},$$

*which, in turn, corresponds to the set*

$$T_{\boldsymbol{\kappa}} = \{\, s_{\boldsymbol{\gamma}} \mid s_{\boldsymbol{\gamma}} \equiv 1 \pmod{\mathbf{0}} \,\}$$

*of dual atoms of $F$. Also, $|F| = |B|$.*

In the example of Section 7.1, we can choose for $F$ the shaded filter of Figure 7.2.

So a cubic extension Cube $K$ of $K$

(i) has a "cubic" congruence lattice (the Boolean lattice $\mathsf{B}_n$, an "$n$-dimensional cube");

(ii) $K$ and its cubic extension, Cube $K$, have the same number of meet-irreducible congruences;

(iii) the congruences $\boldsymbol{\kappa}$ of $K$ extend to Cube $K$. As a rule, the cubic extension has many more congruences than the $\mathrm{Cube}(\boldsymbol{\kappa})$-s.

There are other small examples in Sections 14.2 and 15.2, in particular, in Figures 14.3 and 15.2.

# Index